\documentclass[journal]{IEEEtran}
\usepackage{cite}
\ifCLASSINFOpdf
\else
\fi
\usepackage{graphicx,color}
\usepackage{setspace}
\usepackage{mwe}
\usepackage{multicol}
\usepackage[utf8]{inputenc}
\usepackage[english]{babel}

\usepackage{amssymb,amsmath,amsthm}

\usepackage{amsfonts}
\usepackage{amssymb}
\usepackage{scalefnt}
\usepackage{xr}
\usepackage{flushend}
\usepackage{enumerate}
\usepackage{bbm}
\usepackage{comment}
\newtheorem{theorem}{Theorem}
\newtheorem{proposition}{Proposition}
\newtheorem{corollary}{Corollary}%[theorem]
\newtheorem{lemma}{Lemma}
\newtheorem*{remark}{Remark}
\newtheorem{assumption}{Assumption}
\theoremstyle{definition}
\newtheorem{definition}{Definition}[section]
\theoremstyle{aside}

\newtheorem{problem}{Problem}

\let\Pr\relax
\DeclareMathOperator*{\Pr}{\mathbb{P}}
 
\usepackage{color}

\def\limsup{\mathop{\rm limsup}}
\def\1{\mathbf{1}}

\newcommand{\G}{\mathbb{G}}

\newcommand{\W}{W}

\newcommand{\maxdev}{\Delta_\text{max}(T_0,\delta)}

\def\mbbG{\mathbb{G}}
\newcommand\MM{\mathcal{M}}
\newcommand\NN{\mathcal{N}}

\newcommand\LL{\mathcal{L}}

\allowdisplaybreaks
\ifCLASSOPTIONcompsoc
  \usepackage[caption=false,font=normalsize,labelfont=sf,textfont=sf,labelformat=empty]{subfig}
\else
  \usepackage[caption=false,font=footnotesize,labelformat=empty]{subfig}
\begin{document}
\title{Characterizing Trust and Resilience in Distributed Consensus for Cyberphysical Systems}

\author{Michal Yemini, Angelia Nedi\'c, Andrea J. Goldsmith, Stephanie Gil
\thanks{M.\ Yemini and A.\ J.\ Goldsmith are with the Department of Electrical and Computer Engineering, Princeton University, Princeton,
NJ, 08544 USA, e-mails: \textit {myemini@princeton.edu; goldsmith@princeton.edu}. A.\ Nedi\'c is with the School of Electrical, Computer and Energy Engineering, Arizona State University, Tempe, AZ 85281 USA, email: \textit{Angelia.Nedich@asu.edu}.
S.\ Gil is with the Computer Science Department at the School of Engineering and Applied Sciences, Harvard University,
Cambridge, MA, 02139 USA, email: \textit{sgil@seas.harvard.edu}.}% <-this % stops a space
\thanks{The authors gratefully acknowledge partial funding support from NSF CAREER grant \#2114733, the Alfred P. Sloan Fellowship,  ONR Grant N000141512527 and the Air Force Office of Scientific Research (AFOSR) Grant FA 8750-20-2-0504.}% <-this % stops a space
}

\maketitle

\begin{abstract}
This work considers the problem of resilient consensus where stochastic values of trust between agents are available.  Specifically, we derive a unified mathematical framework to characterize convergence, deviation of the consensus from the true consensus value, and expected convergence rate, when there exists additional information of trust between agents.  We show that under certain conditions on the stochastic trust values and consensus protocol: 1) almost sure convergence to a common limit value is possible even when malicious agents constitute more than half of the network connectivity, 2)  the deviation of the converged limit, from the case where there is no attack, i.e., the true consensus value, can be bounded with probability that approaches 1 exponentially, and 3) correct classification of malicious and legitimate  agents can be attained in finite time almost surely. Further, the expected convergence rate decays exponentially as a function of the quality of the trust observations between agents.
\end{abstract}


% \fi
\IEEEpeerreviewmaketitle

\section{Introduction}
\label{sec:intro}
In this paper we address the problem of multi-agent distributed consensus in the presence of malicious agents. The consensus problem is of core importance to many algorithms and coordinated behaviors in multi-agent systems.  It is well-known however, that these algorithms are vulnerable to malicious activity and/or non-cooperating agents~\cite{origByz,nancy,paxos,lynchPaxos,sundaram} and that several of the existing performance guarantees for the nominal case fail in the absence of cooperation~\cite{bulloUnreliable,security_swarmThreat}. 

Many works have investigated the possibility of attaining resilient consensus, or consensus in the face of non-cooperating or malicious agents. However, achieving resilience often relies on conservative assumptions of connectivity, and bounds on the maximum number of malicious actors.  A classical result for consensus in networks with non-cooperating agents (faulty or malicious) holds from~\cite{origByz,bulloUnreliable}:\\

\noindent Well behaving agents can always agree upon a parameter if and only if the number of malicious agents is \emph{less than 1/2 of the network connectivity}.\footnote{The connectivity of a graph is the maximum number of disjoint paths between any two vertices of the graph.}\\

\noindent This result has long been held as a fundamental requirement for achieving consensus in unreliable networks~\cite{origByz,bulloUnreliable,wmsr}.

For purely cyber systems, where \emph{data} is used to validate information, these limitations are difficult to circumvent.  Cyberphysical systems (CPS) however, are becoming prevalent in many fields, from robotics to power systems to traffic systems and beyond~\cite{CPS_Design}. These systems provide new opportunities to additionally use \emph{physical channels} of information for data validation.  This physical channel data can be exploited to understand ``trust'', or the trustworthiness of data, among the agents in the system. Examples include using observers to compare against system models~\cite{BulloCyberphysSecurity_GeometricPrinciples,pappasDesignActuationSensing}, using camera or other physical channel data~\cite{Candell_PhysicsBasedDetection}, or using transmitted communication channels by comparing against a known carrier signal to find whether a message has been manipulated~\cite{pappasSecretChannelCodes,Sinopoli_Watermark}.  Wireless channels are of particular interest since they often exist in multi-agent systems as a medium for information exchange and can also be used to extract useful information 
such as, for example, tracking motion~\cite{fadel,arrayTrack}, localizing robot agents~\cite{spotFi, RF-Compass}, as well as providing security to wireless systems~\cite{SecureArray}. 

Our recent work in~\cite{AURO} derives the characterization of a stochastic value of trust, $\alpha_{ij}\in(0,1)$, that approaches 1 for trusted transmission and 0 for an untrusted transmission between any two agents $i$ and $j$ in the case of a Sybil Attack.   Importantly, these recent results show that information capturing the trustworthiness of communicating agents can be obtained from analyzing the communication signals themselves - alluding to the fact that the probability of trust between agents can be characterized using information that already exists in multi-agent networks.

However, there remains a need for a unified mathematical framework to guide how a probabilistic understanding of trust between agents can be used to arrive at important performance guarantees for resilient coordination and consensus.

We are motivated by these recent works to derive consensus results that take trust into account. \emph{This paper provides a unified framework that characterizes the value of trust when it comes to achieving resilient consensus in adversarial settings.}
\begin{figure}
  \centering
      \includegraphics[width=0.45\textwidth]{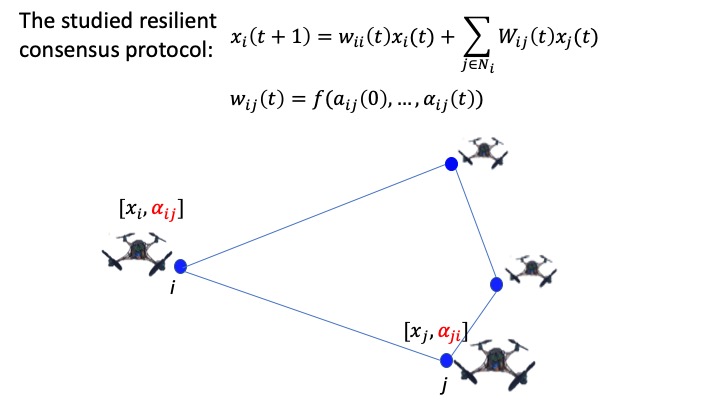}
  \caption{A schematic of a multi-robot system with a consensus protocol modified to use agent transmitted values $x_i$ and $x_j$ as well as observations of inter-agent trust $\alpha_{ij}$. }
  %\label{fig:trust}
  \label{fig:problem}
  \vspace{-0.2in}
\end{figure}
We generalize the concept of trust between agents by considering the availability of stochastic observations $\alpha_{ij}$ that can be thought of as a probability of trust between any two agents $i$ and $j$ in a system (see Fig.~\ref{fig:problem}).  We show that using this additional information can lead to strong results for resilience.  Namely, under certain conditions on $\alpha_{ij}$ that we characterize, we establish the following results:\\
1) \textbf{Convergence:} We show that convergence to a common limit value is possible almost surely, even when the number of malicious agents is larger than less than 1/2 of the network connectivity.\\
2) \textbf{Deviation from nominal average consensus:} Whereas without additional information of trust between agents the presence of malicious agents comprising more than  1/2 of the network connectivity could take the consensus value to any limit, we show that by using $\alpha_{ij}$'s, the influence of malicious agents can be \emph{bounded} and that this bound can be characterized.\\
3) \textbf{Convergence rate:} We show that convergence can be attained in finite time with arbitrarily high probability, and that the expected convergence rate decays exponentially with the quality of the observations $\alpha_{ij}$ and other problem parameters that we characterize.

A roadmap of our paper is as follows: We begin by providing some background and context for the resilient consensus problem in Section~\ref{sec:relatedWork}. In Section~\ref{sec:ProblemDescription} we present our consensus model, our model of trust between agents and the $\alpha_{ij}$ observations, and we characterize the influence of malicious agents on multi-agent consensus. In Section~\ref{sec:Results} we provide our theoretical results and analysis of convergence, deviation, and the convergence rate for consensus with malicious agents. Finally, in Section~\ref{sec:sims} we present simulation results,  followed by Section~\ref{sec:discussion} that provides a discussion and future work directions, and a conclusion in Section~\ref{sec:conclusion}.

\vspace{-0.15in}
\subsection{Related Work}
\label{sec:relatedWork}

The problem of coordination in multi-agent systems and consensus has a long history~\cite{murray,lynchText,pappas,spanos}.  In particular, the consensus problem has enjoyed the development of many fundamental results providing insights into conditions for convergence~\cite{nedicConstrainedConsensus,nedicGeometricConvergence,TsitsiklisSwitchingConvergence,olfati-saber,8340193}, convergence rate~\cite{boyd,nedicConvRate,TsitiklisConvergenceRate}, and general analysis and implementation \cite{additional_reviewer_reference1,additional_reviewer_reference3}.  Similarly, in the current paper we hope to start building a more general mathematical foundation for \emph{resilient consensus} using the concept of \emph{inter-agent trust.}

Resilient consensus refers to consensus in the case where agents are non-cooperative either due to faulty or malicious behavior. This is an important problem due to the prevalence of the non-cooperative case.  For multi-robot systems for example, security is of utmost importance as these robot systems enter the world as autonomous vehicles, delivery drones, or partners in security and defense~\cite{maliciousModel,swarmThreat,surveySecurity}. For cyberphysical systems the issue of security is a highly investigated topic due to the increasing role that cyberphysical systems play in critical infrastructure such as power grids, traffic management systems, and beyond~\cite{sastryBook,scada,Yang13,ControlSysMagazineResilienceIssue,sastry_saurabh_DOS,sastryCPS_survivability}. As a result, there has been increased attention on exposing the vulnerabilities of consensus in multi-agent systems~\cite{origByz,nancy,paxos,lynchPaxos} and on developing new theory to support resilient consensus for these systems~\cite{bulloUnreliable,sundaram,wmsr,prorok,BulloCyberphysSecurity_GeometricPrinciples}.  A common theme in these works is to use transmitted data to identify and thwart attacks, and as a result they are often accompanied by important constraints on the minimum connectivity needed in the graph or requirements on the maximum number of malicious agents that can be tolerated.  Take for example the classical result from~\cite{origByz} which states that well behaving agents can always agree upon a parameter if and only if the number of malicious agents comprises less than 1/2 of the network connectivity. Additionally, the work \cite{wmsr} presents conditions for an asymptotic convergence  to an \textit{arbitrary} point in the convex hull of the initial values of the legitimate agents.

 Cyberphysical systems, including multi-robot systems as a prominent example, offer more physical channels of information that can be used for data validation and for quantifying ``trust'' between agents (robots, sensors, autonomous vehicles, and more). Camera data, communicated wireless signals, lidar, and radar, are a few examples of physical channels of information that often exist in cyberphysical systems and can be used to validate data transmitted between agents.  Indeed, many works have exploited \emph{the physics} of cyberphysical systems to validate transmitted data which presents great promise for the security of these systems~\cite{Candell_PhysicsBasedDetection,dataVeracityCPS_entropyMethods,pappasSecretChannelCodes,pappasDesignActuationSensing,modelingDependability_CPS,Sinopoli_Watermark}.

Advances in sensing using wireless signals is especially promising in this domain as many multi-agent systems transmit messages over a wireless medium~\cite{goldsmith,davidtse}.  The wireless community has demonstrated that these wireless signals themselves contain important information that can be fruitfully exploited for security.  For example it has been shown that wireless signals can be analyzed to provide tracking and localization of agents~\cite{lteye,RF-Compass,fadel,ubicarse,arrayTrack,Wifi_collaborative}, authentication~\cite{SecureArray,kai,ting}, and thwarting of spoofing attacks~\cite{mobicom_adsb}.

\textbf{Deriving Trust Values using Wireless Signals:}
Here we provide some intuition on one particular example of quantifying trust between agents by exploiting physical channels of information; specifically, using wireless channels. Our previous work developed a method for measuring \emph{directional signal profiles} using channel state information (CSI) from the wireless messages over each link $(i,j)$ in the network~\cite{IJRR,ubicarse}. These profiles measure signal strength arriving from every direction in the 3D plane. Directional signal profiles display two important properties: 1) transmissions originating from the same physical agent have very similar profiles and 2) energy can be measured coming from the direct-line path between physical agents. The paper~\cite{AURO} quantifies these properties, providing an analysis that shows both analytically and experimentally, that a single scalar value $\alpha_{ij}\in (-0.5,0.5)$ (shifted by -0.5 from~\cite{AURO}) can be computed for each signal profile that quantifies the likelihood that the transmission is coming from the same physical (spoofed) node, or a unique (legitimate) node; a property critical for thwarting Sybil Attacks. Intuitively, the $\alpha_{ij}$ was shown experimentally and theoretically to be close to -0.5 if one of the agents $j$ is a spoofed node and close to 0.5 if both agents are legitimate nodes in the network~\cite{AURO}. This is captured quantitatively by the bounds on the expectation of the $\alpha_{ij}$. The current paper takes inspiration from this previous work which shows the existence of such $\alpha_{ij}$ variables for stochastically characterizing trust between agents, but goes beyond to define a unified mathematical framework for using trust to obtain certain performance guarantees for consensus in multi-agent systems. As such, Definition~\ref{def:alpha} specifies the mildest characteristics of such as $\alpha_{ij}$ term in order to achieve strong performance guarantees. Thus the results contained herein would be compatible with stochastic variables as derived in~\cite{AURO} or any other physical channels that satisfy Definition~\ref{def:alpha}. We note however, that the focus of this paper is not on the derivation of stochastic trust values. Rather, we focus on the development of the mathematical foundations for why and how such values can be key to achieving certain performance guarantees that would be difficult or impossible to obtain otherwise.

\textbf{Contributions in the context of previous work:} This paper builds most closely off of our previous work in~\cite{AURO,GilLCSS} which shows that a stochastic variable $\alpha_{ij}(t)$ can be extracted from wireless signals between agents to verify the validity of an agent's identity (in a probabilistic sense) during a Sybil Attack. Further,~\cite{GilLCSS,ICRA2019} shows the potential of using these values to arrive at resilient consensus. In the current work we assume the availability of such probabilistic observations of trust between agents, the $\alpha_{ij}$, and focus on the development of the mathematical machinery that can utilize these trust values to arrive at strong results for resilient consensus. Key differentiating results include: 1)~We
provide convergence rate results for the case where the number of malicious agents is larger than 1/2 of the network connectivity, showing that convergence in \emph{finite time} is possible with arbitrarily high probability. To the best of our knowledge, this is the first time that this result has been formally proven. 2)~We divide the consensus algorithm into two critical stages with broad consequences for the attainable performance guarantees. In the first stage the agents \emph{observe the network}, i.e. they collect values of  $\alpha_{ij}$ to detect malicious agents but do not start running the consensus protocol. In the second stage, agents continue to collect values of  $\alpha_{ij}$ and also adapt their data values using a modified consensus protocol that includes received data only from agents they classify as trustworthy. We show that the impact of observing the network for some time window $T_0$ is highly consequential for deviation and convergence, where convergence rate and deviation amount decreases exponentially with increasing $T_0$. 3)~We use weights that define a switching topology allowing for arbitrarily small deviation in the final consensus value. 4)~The current work allows for the use of non-symmetric weights which can accommodate larger classes of multi-agent systems. Finally, the results in this paper are general to many classes of attacks so long as certain characteristics on the probability of trust between agents, $\alpha_{ij}$, are met.  We establish these key characteristics and define their roles in convergence, deviation, and convergence rate to consensus.

\section{Problem Description}
\label{sec:ProblemDescription}

We consider the problem of consensus for multi-agent systems. The agents in the system communicate over a network, which is  represented by an undirected graph, $\mathbb{G}=(\mathbb{V},\mathbb{E})$,
where $\mathbb{V}=\{1,\hdots,n\}$ denotes the set of node indices for the agents and $\mathbb{E}\subset\mathbb{V}\times\mathbb{V}$ 
denotes the set of undirected edges. We use $\{i,j\}$ to denote 
the edge connecting agents $i$ and $j$. 
The set of \emph{neighbors} of node $i$ is denoted by $\mathcal{N}_i= \{j\in\mathbb{V}\mid \{i,j\}\in\mathbb{E}\}$.
Note that by this definition if $i$ is neighbor of $j$, then $j$ is also a neighbor of $i$.
With the graph $\mbbG$, we will associate time-varying nonnegative weights $w_{ij}(t)$, $\{i,j\}\in\mathbb{E},$ which capture the information exchange process among the agents at a given time $t$. We use $w_{ij}(t)$ to denote the $ij$th entry of the weight matrix $W(t)$. For these matrices, we have $w_{ij}(t)\ge 0$ and  $w_{ji}(t)\ge 0$
only if $\{i,j\}\in\mathbb{E}$, and $w_{ij}(t)=w_{ji}(t)=0$ otherwise. In this way, the positive entries of the matrices $W(t)$ are associated only with the edges in the set $\mathbb{E}$. Additionally, the matrices $W(t)$ need not be symmetric.
We consider the case where a subset of nodes with indices denoted by the set $\MM$, $\MM\subset \mathbb{V}$, are either \emph{malicious} or \emph{malfunctioning} and, thus, do not reliably cooperate in the consensus protocol. The set $\MM$ is assumed to be unknown, with the cardinality $|\MM|$. An agent that is not malicious or malfunctioning is termed \textit{legitimate}, and the set of legitimate agents is denoted $\mathcal{L}$, $\mathcal{L}\subset \mathbb{V}$,
with the cardinality $|\LL|$. 

In what follows, we use $\boldsymbol{0}$ to denote the zero matrix, where the dimension of $\boldsymbol{0}$
is to be understood from the context.
We say that a matrix $A$ is nonnegative if its entries are nonnegative, i.e., $A_{ij}\textcolor{black}{\ge}0$ for all $i,j$, and we write $A\ge\mathbf{0}$ to denote that $A$ is nonnegative. We write $A>\mathbf{0}$ to denote a matrix having all entries positive,  i.e., $A_{ij}>0$ for all $i,j$.
A nonnegative matrix is row-stochastic if the sum
of its entries in every row is equal to 1, and it is row-substochastic if the sum of its entries in every row is equal to or less than 1. We use $\1$ to denote the vector with all entries equal to 1. Given a vector $x$, we use
$x_i$ or $[x]_i$ (when $x$ has a subsrcipt)
to denote its $i$th entry. A similar notation is used for the entries of a matrix.

\vspace{-0.1in}
\subsection{The Model}
\label{sec:ConsensusMode}
We are concerned with the problem of consensus where each \textit{legitimate} agent updates its value according to the following update equation (see \cite{sundaram}) for all $t\ge T_0-1$,
\begin{align}
\label{eq:consensusProtocol}
    x_i(t+1)=w_{ii}(t)x_i(t)+\sum_{j\in\mathcal{N}_i}w_{ij}(t)x_j(t),
\end{align}
where $x_i(t)\in\mathbb{R}$, while
the weights $w_{ij}(t)$ are nonnegative and sum to 1, i.e.,
$w_{ii}(t)>0$, $w_{ij}(t)\ge 0$ 
for $j\in {\cal N}_i$ and $w_{ii}(t)+\sum_{j\in{\cal N}_i}w_{ij}(t)=1$.
The process is initiated at some time $T_0\ge0$ with the agents' initial values $x_i(T_0-1)=x_i(0)$ for all $i\in\mathbb{V}$.

\begin{definition}[Malicious agent]
Agent $i$ is said to be \emph{malicious} if it does not follow the rule in~\eqref{eq:consensusProtocol} for updating its value $x_i(t)$ at some time $t\ge T_0-1$.
\end{definition}

Let $x(t)\in\mathbb{R}^n$ denote the vector of values $x_i(t)\in\mathbb{R}$ for all the agents at time $t$. Without loss of generality, we assume that 
the agent indices are ordered in such a way that the vector $x(t)$ can be separated in two components; the first component representing the legitimate agent values $x_\LL(t)$ and the second component representing the malicious agent values $x_{\MM}(t)$. In this way, the consensus dynamics takes the following form:

\begin{align}
\label{eq:dynamics}
\begin{bmatrix}
    x_\LL(t+1)\\
    x_{\MM}(t+1)\\
\end{bmatrix}
    =
    \begin{bmatrix}
    W_{\LL}(t)       & W_{\MM}(t) \\
    \Theta(t)    & \Omega(t) \\
    \end{bmatrix}
\begin{bmatrix}    
    x_{\LL}(t) \\
    x_{\MM}(t) \\
\end{bmatrix},
\end{align}
where $W_\LL(t)\in\mathbb{R}^{|\LL|\times |\LL|}$ is the matrix multiplying the state component corresponding to the legitimate agents' values, and $W_\MM(t)\in\mathbb{R}^{|\LL|\times |\MM|}$ is the matrix multiplying the state component 
corresponding to the malicious agents' values.
%and $\mathbb{W}(t)=[W_\LL(t)\ W_\MM(t)]$. 
The matrices $\Theta(t)$ and $\Omega(t)$ dictate the dynamics of the malicious agents' values and are assumed to be unknown. For the sake of simplicity, we assume that a malicious agent sends the same data value to all of its neighbors, however, our analysis holds for the case where a malicious agent can send different data values to different neighbors at each iteration.
We consider the case when a parameter $\eta>0$ is known to all legitimate agents and it is known that the values $|x_i(t)|$ should not exceed $\eta$ at any time. 
The assumption $|x_j(t)|\leq \eta$  
for all $t$ and $j\in \MM$, is reasonable, since 
an agent $j$ can be classified as malicious based on its value
(i.e., if $|x_j(t)|> \eta$ at some time $t$). 
Thus, in this paper, our focus is on the update rule in~\eqref{eq:dynamics} where $|x_j(t)|\leq \eta$ for all $j\in\mathcal{L}\cup\mathcal{M}$ and time instants $t$.

In this paper we are interested in the case where each 
transmission from agent $j$ to agent $i$ can be tagged with a random observation $\alpha_{ij}(t)\in[0,1]$ of a random variable
$\alpha_{ij}\in[0,1]$. The random variable $\alpha_{ij}$ represents a \emph{probability of trust} that agent $i$ can give to its neighbor $j\in \mathcal{N}_i$. 

\begin{definition}[$\alpha_{ij}$] 
\label{def:alpha}
For every $i\in\LL$ and $j\in\NN_i$,
the random variable $\alpha_{ij}$ taking values in the interval $[0,1]$ represents the probability that 
agent $j\in\mathcal{N}_i$ is a trustworthy neighbor of agent $i$. We assume the availability of such observations $\alpha_{ij}(t)$ throughout the paper.
\end{definition}

We refer to~\cite{AURO} for an example of such an $\alpha_{ij}$ value. Intuitively, a random realization $\alpha_{ij}(t)$ of $\alpha_{ij}$ 
contains useful trust information if the legitimacy of the transmission can be thresholded.  

In this paper we assume that a value of $\alpha_{ij}(t)>1/2$ indicates a legitimate transmission and $\alpha_{ij}(t)<1/2$ indicates a malicious transmission in a stochastic sense (misclassifications are possible).  Note that $\alpha_{ij}(t)=1/2$ means that the observation is completely ambiguous and contains no useful trust information for the transmission at time $t$.

Our ultimate goal is to understand the convergence, and convergence rate properties of consensus under the influence of a malicious attack if such observations of trust are available.
Toward this goal, our first objective is to construct the appropriate weight matrix $W(t)=[W_\LL(t)\ W_\MM(t)]$ 
for legitimate agents' state dynamics in~\eqref{eq:dynamics} by using some function of the observations $\alpha_{ij}(t),\alpha_{ij}(t-1),\ldots,\alpha_{ij}(0)$,  where $i$ is a legitimate agent so that $i\in\LL$
and $j$ is a neighbor of $i$ so that $j\in\mathcal{N}_i$.
In this way, the legitimate agents' ability to isolate the malicious agents is reflected in the convergence of the weight matrices $W_\MM(t)$ to the zero matrix and the convergence of the weight matrices $W_\LL(t)$ to a row-stochastic matrix, where the convergence is in a probabilistic sense (to be defined shortly).
Necessary characteristics of $\alpha_{ij}$ to allow for the convergence of the linear protocol in~\eqref{eq:dynamics}, possibly including its distribution, bounds on expectation, etc., will be derived in this work.

In line with this model, in addition to broadcasting its value $x_i(t)$ at each time $t$, we assume that a message from each agent $i\in \mathbb{V}$ is tagged with an observation $\alpha_{ij}(t)$ for every neighboring agent $j$, i.e.,  $\{i,j\}\in \mathbb{E}$ (see Figure~\ref{fig:problem}). \emph{\textbf{The purpose of this work is to develop the necessary properties of $\alpha_{ij}(t)$ such that consensus can be maintained in the presence of malicious or malfunctioning agents.}} 

Suppose that, in the interval of time $[0,T_0-1)$, the legitimate agents only use their measurement of the trustworthiness values $\alpha_{ij}$, $j\in \mathcal{N}_i$, $i\in \LL$. We refer to this interval of time as an \emph{observation window} where agents can amass a history of trust values $\alpha_{ij}$ before starting the consensus.  At time $T_0-1$, the agents start the data passing phase with the weights $w_{ij}(T_0-1)=f(\alpha_{ij}(0),\ldots,\alpha_{ij}(T_0-1))\cdot\mathbbm{1}_{\{j\in\mathcal{N}_i\}}$, where $f$ is some function capturing the history of all observations $\alpha_{ij}(0),\ldots,\alpha_{ij}(T_0-1)$ and $\mathbbm{1}_{\{\cdot\}}$ denotes the indicator function.  Our goal is to study the resulting consensus dynamics over time. 
Since we start the consensus algorithm~\eqref{eq:dynamics} from time $t=T_0-1$, we have that $x_{\LL}(t)=x_{\LL}(0)$ for all $t\in[0,T_0-1]$. 

Let
$\prod_{k=p}^rH_k$ denote the backward product 
of the matrices $H_p,\ldots, H_r$, i.e.,
\[\prod_{k=p}^rH_k =\begin{cases}
H_{r}\cdots H_{p+1} H_{p}\quad&\text{ if }r\ge p,\\
I \quad&\text{ if }r< p,
\end{cases}
\]
where $I$ denotes the identity matrix.
 Following \eqref{eq:dynamics} the value of the legitimate agents at any point in time $t$ has two salient influence terms, the values contributed by other legitimate agents $\tilde{x}_{\LL}(T_0,t)$ and the values contributed by other malicious agents $\phi_{\MM}(T_0,t)$ so that for all $t\geq T_0$
\begin{align}
\label{eq:leg-state}
    x_{\mathcal{L}}(T_0,t) = \tilde{x}_{\LL}(T_0,t)+\phi_{\MM}(T_0,t),
\end{align}
where
\begin{align}
\label{eq:leg-influence}
    \tilde{x}_{\LL}(T_0,t)=\left(\prod_{k=T_0-1}^{t-1}W_{\LL}(k)\right)x_{\LL}(0),
\end{align}
and 
\begin{align}
\label{eq:mal-influence}
    \phi_{\MM}(T_0,t) = \sum_{k=T_0-1}^{t-1}\left(\prod_{l=k+1}^{t-1}W_{\mathcal{L}}(l)\right)W_{\MM}(k)x_{\MM}(k).
\end{align}

Relation~\eqref{eq:leg-state} shows explicitly how the states of legitimate agents depend on their initial states at time $T_0$ and on the states of malicious agents. The term  $\phi_{\MM}(T_0,t)$ given in~\eqref{eq:mal-influence} captures the total influence of the malicious agents on the states of the legitimate agents from the initial time $T_0$ to the current time $t$ and will be the primary subject of study.  Specifically, we will study attainable bounds on the malicious agents' influence term $\phi_{\MM}(T_0,t)$.

\begin{remark}[Arbitrary starting time $T_0$ of the consensus.] We emphasize that the consensus updates from Eq.~\eqref{eq:dynamics} begin at an arbitrary time $T_0\geq 0$ and that before this time agent values $x_i(t)$ are not updated for $i\in\LL$.  In other words, for $t\in[0:T_0)$ only observations $\alpha_{ij}(t)$ are collected by each agent $i$ for its neighbors $j\in \NN_i$. The choice of $T_0$ is a user-selected parameter that affects the overall influence of malicious agents in the network and the convergence rate according to a relationship that we later characterize. 
\end{remark}

We now construct a modified weight matrix $W(t)=[W_\LL(t)\,\ W_\MM(t)]$
that governs the state dynamics of legitimate agents (cf.~\eqref{eq:dynamics}) and utilizes the $\alpha$ observations. Intuitively, we wish to assign larger weights to transmissions that are most likely to originate from \emph{legitimate} agents, and smaller or zero weights to transmissions that are most likely to originate from \emph{malicious} agents as dictated by the observations $\alpha_{ij}(t)$ for every neighbor $j$ of $i$. Toward this end,
recalling that $\mathcal{N}_i$ is the complete set of neighbors of legitimate agent $i$, we define the following two quantities:

For the function $f$ that captures the history of observations $\alpha_{ij}(t)$, we choose the sum, and define   $\beta_{ij}(t)$ as follows:
\begin{align}
    \label{eq:betas}
    \beta_{ij}(t)=\sum_{k=0}^{t} \left(\alpha_{ij}(k)-1/2\right)
\hbox{ for $t\ge0,i\in\LL, j\in\mathcal{N}_i$}.\end{align}
Intuitively, following the discussion after Definition~\ref{def:alpha} of $\alpha_{ij}$, the $\beta_{ij}(t)$ will tend towards positive values for legitimate agent transmissions $i\in\LL,\,j\in \NN_i\cap\LL$, and will tend towards negative values for malicious agent transmissions where  $i\in\LL,j\in\NN_i\cap\MM$.

We also define a time dependent \emph{trusted neighborhood} for each agent $i$ as:
\begin{align}%\label{def:N_i_t}
\label{eq:Ni_t}
    \mathcal{N}_i(t) = \{j\in\mathcal{N}_i:\:\beta_{ij}(t)\geq0\}.
\end{align} 
This is the set of neighbors that legitimate agent $i$ classifies as its legitimate neighbors at time $t$.
 Denote
\[n_{w_i}(t) = \max\{\kappa,|\mathcal{N}_i(t)|+1\}\geq1\qquad\hbox{for all }i\in\LL,\]
where
 $\kappa >0$ can be thought of as a parameter limiting the maximum influence that the neighbors of agent $i$ are allowed to have on agent $i$'s update values (i.e. a value of $\kappa$ approaching infinity would preclude neighbors from influencing agent $i$'s updated values).
Using these quantities, we define the weight matrix $W(t)$ by choosing its 
entries $w_{ij}(t)$ as follows:
for every $i\in\LL$, $j\in\NN_i$,
\begin{align}
\label{eq:weights}
    w_{ij}(t) %\nonumber\\
    = 
    \begin{cases}
    \frac{1}{n_{w_i}(t)} & \text{if\,} ,j\in \mathcal{N}_i(t),\\
    0 & \text{if\,} ,j\notin \mathcal{N}_i(t)\cup \{i\},\\
    1-\displaystyle{\sum_{m\in\mathcal{N}_i}w_{im}(t)} & \text{if\,} j=i.
    \end{cases}
\end{align}
 Note that the dependence of the weights $w_{ij}(t)$ on the trust observation history $\beta_{ij}(t)$ comes in through the choice of time-dependent trusted neighborhood $\mathcal{N}_i(t)$ (cf. Eqn.~\ref{eq:Ni_t}). Note also that the matrix $W(t)$ is (rectangular) row-stochastic, by construction. Furthermore, since $\alpha_{ij}(t)$ are random, so are the quantities $\beta_{ij}(t)$ and the set of neighbors of agents $i$,  $\mathcal{N}_i(t)$, that agent $i$ classifies as legitimate at time $t$.
Consequently, some entries of the matrix $W(t)$ are also random,
as seen from~\eqref{eq:weights}. We use the parameter $\kappa$ later on, in our analysis of the dynamics in~\eqref{eq:dynamics}, to obtain an upper bound
on the entries of the matrix $W(t)$.

We define convergence as follows:
\begin{definition}[Convergence of the consensus protocol] 
\label{def:convergence}
We define consensus to be achieved if:\\
1) \textbf{\underline{Limit values:}}
There exists, almost surely, a random variable  $y(T_0)$ such that 
\begin{align}
    &\lim_{t\rightarrow \infty}\tilde{x}_{\LL}(T_0,t)=y(T_0)\boldsymbol{1},
\end{align}
where $\boldsymbol{1}\in\mathbb{R}^{|\LL|\times1}$, and
$y(T_0)$ is in the convex hull of 
the legitimate agents' initial values \textcolor{black}{$x_i(0)$}, $i\in\LL$ \textcolor{black}{and $0$}, and its distribution depends on the starting time $T_0$ of the consensus algorithm. Further, almost surely, a limit exists for the legitimate agent values which is a random variable $z(T_0)$ such that 
\begin{align}
    &\lim_{t\rightarrow \infty} x_{\LL}(T_0,t)=z(T_0)\boldsymbol{1}.
\end{align}
Here, $z(T_0)$ is in the convex hull of the initial legitimate agent values \textcolor{black}{$x_i(0)$}, $i\in\LL$ \textcolor{black}{and the malicious inputs $x_i(t), i\in\MM$, $t\geq T_0-1$ such that malicious agent $i$ is misclassified by a legitimate agent at time $t$}.
Its deviation from the nominal consensus value (the case with no malicious agents) depends on the starting time $T_0$ of the consensus algorithm.  See Def.~\ref{def:convergence}.2 below.\\
2) \textbf{\underline{Deviation:}} There is a probabilistic upper bound on the asymptotic
distance between the legitimate agent values and a suitably defined weighted average of their initial values, i.e.,
\begin{spacing}{0.75}
\begin{align}
&\Pr \left(\max_{i\in\LL}\,\limsup_{t\rightarrow \infty} \left| \left[x_{\LL}(T_0,t)-\boldsymbol{1}v'x_{\LL}(0)\right]_i \right| \leq \Delta(T_0, \delta) \right)\nonumber\\
&\qquad \ge 1 - \delta, \\
&\text{ for some finite } \Delta(T_0,\delta) \leq \Delta_\text{max}(T_0,\delta) , \nonumber 
\end{align}
where $\delta>0$ is a user-defined error probability,
while $v$ is the Perron-Frobenius left-eigenvector of the matrix $\overline{W}_{\LL}$ defined in~\eqref{eq:limitweight} such that $v'\boldsymbol{1}=1$.
\end{spacing}
\end{definition}

We next discuss the assumptions that we use in the sequel. For this, we let $\mbbG_{\LL}=(\LL,\mathbb{E}_{\LL})$ be 
the subgraph of the graph $\mbbG$ that is induced by the legitimate agents, where 
$\mathbb{E}_{\LL}\subset\mathbb{E}$ is given by
$\mathbb{E}_{\LL}=\{\{i,j\}\mid \{i,j\}\in\mathbb{E},\ i\in\LL,j\in\LL\}$.
We use the following assumptions throughout the paper:
\begin{assumption}\label{assumptions}
Suppose that the following hold:\\
1. \emph{[Sufficiently connected graph]}   
    The subgraph $\mbbG_{\LL}$ induced by the legitimate agents is connected.\\
2. \emph{[Homogeneity of trust variables]}
    The expectation of the variables $\alpha_{ij}$ are constant for the case of malicious transmissions and legitimate transmissions, respectively, i.e., for some scalars $c,d$ with $c\ne d$, 
    \[c=E(\alpha_{ij}(t))-1/2\qquad
    \hbox{for all $i\in \LL,\ j\in \NN_i\cap\MM$},\] 
    \[d=E(\alpha_{ij}(t))-1/2\qquad
    \hbox{for all $i\in\LL,\ j\in \NN_i\cap\LL$}. \]
3. \emph{[Independence of trust observations]}  
    The observations $\alpha_{ij}(t)$ are independent for all $t$ and all pairs of agents $i$ and $j$, with $i\in \LL$, $j\in\mathcal{N}_i$.
    Moreover, for any $i\in \LL$ and $j\in\mathcal{N}_i$, the observation sequence
    $\{\alpha_{ij}(t)\}$ is identically distributed.
\end{assumption}

Assumption~\ref{assumptions}.1 captures the general connectivity requirement in consensus networks with a fixed topology. 
The assumption on the homogeneity of trust variables (cf.\ Assumption~\ref{assumptions}.2) is made only for the ease of exposition. The assumption can be generalized to the heterogeneous case as discussed in Section~\ref{sec:discussion}.

\begin{figure}
  \centering
      \includegraphics[width=0.35\textwidth]{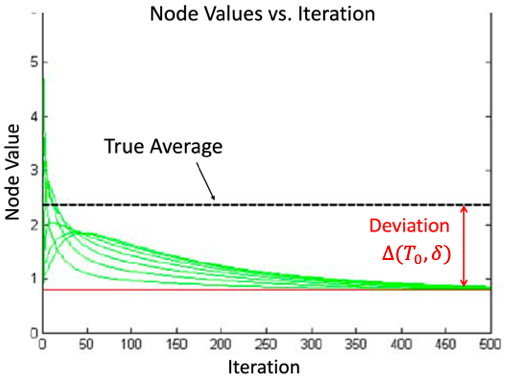}
  \caption{A schematic showing the \emph{deviation} from the nominal average in the case of a malicious agent attack and no resilience in the consensus protocol.\label{fig:deviation}}
  \vspace{-0.2in}
\end{figure}

\subsection{Problem Definition}
\label{sec:problem-defn}

In the following we summarize the three main problems that we aim to address under the assumptions in Assumption~\ref{assumptions}:
\begin{problem}[Convergence of the consensus protocol]
\label{sec:problem1}
We aim to show that in the presence of an adversarial attack, consensus to a common limit value is attainable \emph{even when malicious agents comprise $>1/2$ of the network connectivity} if observations $\alpha_{ij}(t)$ are available. We seek conditions on $\alpha_{ij}(t)$, and a definition of the weights $\W(t)$ that allow convergence.
\end{problem}

\begin{problem}[Bounded deviation for average consensus]
\label{sec:problem2}
For the case of average consensus, we aim to find a bound on deviation from the average consensus value, $\Delta\left(T_0, \delta\right)$ that can be achieved with a probability at least $1-\delta$, where the bound is given as a function of problem parameters (such as the number of legitimate and the number of malicious agents) and the characteristics of $\alpha_{ij}$ (such as bounds on its expected value $E(\alpha_{ij})$).
\end{problem}

\begin{problem}[Characterization of convergence rate $\tau$]
\label{sec:problem3}
We aim to determine a finite convergence time for the correct classification of legitimate and malicious agents almost surely. Additionally, we aim to identify the rate of convergence as a function of the stochastic observations $\alpha_{ij}$ and other problem parameters that we characterize.
\end{problem}

\section{Analysis}
\label{sec:Results}

\subsection{Convergence of Consensus of Legitimate Agents Inputs}

In this section, we analyze the state dynamics among the legitimate agents, that is governed by the weight matrix $W_{\LL}(k)$, as seen from~\eqref{eq:dynamics}. Recalling that 
$W(t)=[W_\LL(t)\, W_\MM(t)]$, from the definition of the matrix $W(t)$ (cf.~\eqref{eq:weights}), we see that 
the matrix $W_\LL(t)$ is a sub-stochastic matrix, since the sum of some rows may be strictly less than one. This section analyzes the convergence of the (backward) product $\prod_{k=0}^{t}W_{\LL}(k)$, as $t\to\infty$. 
We prove the strong ergodicity of this product, and show that the limit matrix $\prod_{k=0}^{\infty}W_{\LL}(k)$ has strictly positive entries, i.e., $\prod_{k=0}^{\infty}W_{\LL}(k)>\boldsymbol{0}$.

Our goal here is to show convergence of the consensus system in Eq.~\eqref{eq:dynamics} using the choice of weights from Eq.~\eqref{eq:weights} that exploit the observations $\alpha_{ij}(t)$. Specifically, we aim to characterize the conditions on $\alpha_{ij}(t)$ necessary to achieve consensus in the presence of malicious agents \emph{even when the number of malicious agents is arbitrarily high}.

We will achieve this in three parts: first we show that for a sufficiently connected network (cf. Assumption~\ref{assumptions}.1) the weight matrix in Eq.~\eqref{eq:weights} reaches its limit described in Eq.~\eqref{eq:limitweight} in finite time $T_f$ almost surely and that this limit is a primitive matrix (cf. Lemma~\ref{lemma:weight_mat_primitive} and Proposition~\ref{proposition:consensus_not_vanishing}). Second we show that this in turn implies that the legitimate agents' influence captured in Eq.~\eqref{eq:leg-influence} approaches a finite limit almost surely (cf. Proposition~\ref{proposition:convergence_legitimate_part}).  By the same token, the malicious agents' influence captured in Eq.~\eqref{eq:mal-influence} also approaches a finite limit almost surely (c.f. Proposition~\ref{proposition:convergence_malicious_part}).  Finally, putting these limits together we show that the values of the legitimate agents, $x_\LL(T_0,t)$, approach a finite limit, i.e. we achieve convergence.

In what follows we use the notion of a primitive matrix,
which is provided next.
\begin{definition}[Primitive matrix]
A nonnegative square matrix 
$A$ is said to be a primitive matrix if there exists 
an integer $k\ge 1$ such that $A^k>\boldsymbol{0}$. 
\end{definition}

Let us define a matrix $\overline{W}_{\LL}$ with the entries given by: for every $i,j\in\LL$,

\begin{align}
\label{eq:limitweight}
    [\overline{W}_\LL]_{ij} = \begin{cases}
    \frac{1}{\max\{|\mathcal{N}_i\cap\LL|+1,\kappa\}} & \text{ if } 
    j\in \mathcal{N}_i\cap \LL,\\
    1-\frac{|\mathcal{N}_i\cap\LL|}{\max\{|\mathcal{N}_i\cap\LL|+1,\kappa\}} & \text{ if } j=i,\\
    0 & \text{otherwise}.
    \end{cases}
\end{align}

For the matrix $\overline{W}_{\LL}$, we have the following result.

\begin{lemma}\label{lemma:weight_mat_primitive}
The matrix $\overline{W}_{\LL}$ is primitive.
\end{lemma}
\begin{IEEEproof}
The result follows directly by the definition of the matrix $\overline{W}_{\LL}$ in~\eqref{eq:limitweight} and the assumption that the graph $\mbbG_\LL$ is connected 
(cf.~Assumption 1.1).
\end{IEEEproof}

Under Assumptions \ref{assumptions}.2--\ref{assumptions}.3, we next provide a result for the random quantities $\beta_{ij}(t)$ defining the weight matrix $W(t)$.
\begin{lemma}
\label{Lemma:concentration_upper}
Consider the random variables $\beta_{ij}(t)$ 
as defined in~\eqref{eq:betas}.
Then, for every $t\geq 0$ and every $i\in \LL$, $j\in\NN_i\cap\LL,$
\begin{align*}
\Pr\left(\beta_{ij}(t)< 0\right)\leq \max\{\exp(-2(t+1)d^2),\mathbbm{1}_{\{d<0\}}\}.
\end{align*}
Additionally, for every $t\geq 0$ and every 
$i\in \LL,j\in \NN_i\cap\MM$,
\begin{align*}
\Pr\left(\beta_{ij}(t)\geq 0\right)\leq \max\{\exp(-2(t+1)c^2),\mathbbm{1}_{\{c>0\}}\}.
\end{align*}
\end{lemma}
\begin{IEEEproof}
By the linearity of the expectation, we have
 $   E(\beta_{ij}(t)) = \sum_{k=0}^{t}
    \left(E(\alpha_{ij}(k))-1/2\right)$.
For every $k\geq 0$, and for all $i,j$ with $i\in\LL$,
$j\in\NN_i\cap\LL$, by
Assumptions~\ref{assumptions}.2--\ref{assumptions}.3,
we have $d=E(\alpha_{ij}(k))-1/2$. 
Additionally, for every $k\geq 0$ and $i\in\LL, j\in\NN_i\cap\MM$, 
by Assumptions~\ref{assumptions}.2--\ref{assumptions}.3,
we have $c=E(\alpha_{ij}(k))-1/2$. 
The result then follows directly from 
the Chernoff-Hoeffding Inequality (see \cite[Theorem 1.1]{Dubhashi2009ConcentrationOM}).
\end{IEEEproof}

From Lemma \ref{Lemma:concentration_upper} one can observe the following. 
\begin{corollary}[Bounds on Expectation of $\alpha_{ij}(t)$]
\label{cor:expectationBoundsAlpha}
For the choice $\beta_{ij}(t) = \sum_{k=0}^{t}[\alpha_{ij}(k)-1/2]$, we must have that
$c=E(\alpha_{ij}(k))-1/2<0$ for all $i\in \LL,\ j\in \MM, k\geq0$ and $d=E(\alpha_{ij}(k))-1/2>0$, for all $i,j\in\LL, k\geq0$ in order to have a decaying misclassification probabilities. 
\end{corollary}

We provide some intuition about this corollary as it has significant implications for the characteristics of $\alpha_{ij}(t)$ necessary for obtaining several important results. The condition that $c$ and $d$ be bounded away from zero in Corollary~\ref{cor:expectationBoundsAlpha} intuitively means that there is information captured by the $\alpha_{ij}(t)$ observations. The case of $c=0$, $d=0$ is equivalent to saying that, in expectation, the observations $\alpha_{ij}(t)$ will be the same regardless of the transmission being legitimate or malicious (the observations $\alpha_{ij}(t)$ contain no useful information). \textit{Hereafter, we assume that $c<0$ and $d>0$}, which means that $\alpha_{ij}(t)$ tends to be less than $1/2$ for malicious transmissions (has expectation value less than $1/2$) and $\alpha_{ij}(t)$ tends to be greater than $1/2$ for legitimate transmissions (has expectation value greater than $1/2$). 

\begin{figure}
  \centering
      \includegraphics[width=0.35\textwidth]{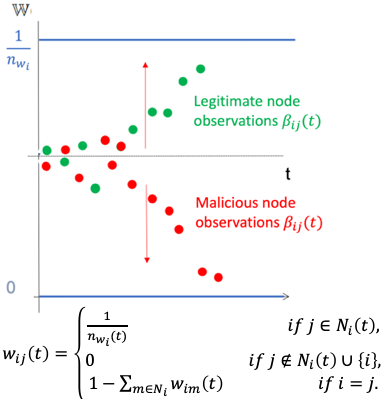}
  \caption{A depiction of the $\beta_{ij}(t)$ values for legitimate and malicious nodes becoming more separated and accurately classifiable with higher probability as a longer history of $\alpha_{ij}(t)$ becomes available (i.e. as $t$ gets larger).}
  \vspace{-0.2in}
\end{figure}

\begin{proposition}
\label{proposition:consensus_not_vanishing} 
There exists a finite time instant $T_f>0$ such that  $W_{\LL}(k)=\overline{W}_{\LL}$ for all $k\ge T_f$ almost surely. Moreover, there exists a stochastic vector $v>\mathbf{0}$ such that 
\[\lim_{k\rightarrow\infty}\overline{W}^k_{\LL}=\1 v'.\]
In particular, 
$\prod_{k=T_0-1}^{\infty}W_{\LL}(k)>\boldsymbol{0}$ for every $T_0\geq0$ almost surely.
\end{proposition}

\begin{IEEEproof}
By Lemma \ref{Lemma:concentration_upper} we have that for every $k\geq0$ and $i\in\LL$, $j\in\NN_i\cap\LL$,
\begin{align}
&\Pr\left(\beta_{ij}(k)<0\right)\leq \exp\left(-2(k+1)d^2\right).
\end{align}
Additionally, for every $k\geq0$ and 
$i\in \LL$, $j\in \NN_i\cap\MM$,
\begin{align}
&\Pr\left(\beta_{ij}(k)\geq 0\right)\leq \exp\left(-2(k+1)c^2\right).
\end{align}
Note that we have
{$\sum_{k=0}^{\infty}\exp\left(-2(k+1)c^2\right)<\infty$} and 
$\sum_{k=0}^{\infty}\exp\left(-2(k+1)d^2\right)<\infty$. Thus, 
by the Borel-Cantelli Lemma, 
the events $\left\{\beta_{ij}(k)\geq 0\:\:\: \forall \: i\in \LL,j\in \NN_i\cap\MM\right\}$ and 
$\left\{\beta_{ij}(k)< 0\:\:\:\forall\: i\in \LL, j\in\NN_i\cap\LL\right\}$ occur only finitely often almost surely. Therefore, 
there exists a (random) finite time
$T_f>0$  such that  $W_{\LL}(k)=\overline{W}_{\LL}$ for all $k\ge T_f$ almost surely.

Now, almost surely, we have
\begin{align*}
        \prod_{k=T_0-1}^{\infty}W_{\LL}(k) 
        &= 
        \prod_{k=\max\{T_f,T_0\}}^{\infty}W_{\LL}(k)\prod_{k=T_0-1}^{\max\{T_f,T_0\}-1}W_{\LL}(k)\cr
        &=\prod_{k=\max\{T_f,T_0\}}^{\infty}
        \overline{W}_{\LL}\prod_{k=T_0-1}^{\max\{T_f,T_0\}-1}W_{\LL}(k)\cr
        &=\lim_{k\to\infty}  \overline{W}_{\LL}^{k-\max\{T_f,T_0\}}
        \prod_{k=T_0-1}^{\max\{T_f,T_0\}-1}W_{\LL}(k).
\end{align*}
By Lemma~\ref{lemma:weight_mat_primitive}, 
the matrix $\overline{W}_\LL$ is primitive.
Therefore,
by the Perron-Frobenius Theorem (see \cite{Perron_Frobenius_ref}), there exists a stochastic column vector 
$v>\mathbf{0}$ such that  $\lim_{k\rightarrow\infty}\overline{W}^k_{\LL}=\1 v'$.

Consequently, almost surely we have
\[\prod_{k=T_0-1}^{\infty}W_{\LL}(k)
= \1v'
\left(\prod_{k=T_0-1}^{\max\{T_f,T_0\}-1}W_{\LL}(k)\right).\]
Since $v>\mathbf{0}$ and the diagonal entries of the matrices $W_{\LL}(k)$ are positive (by construction),
it follows that 
$v'\left(
\prod_{k=T_0-1}^{\max\{T_f,T_0\}-1}W_{\LL}(k)\right)
>\mathbf{0},$
implying that almost surely 
$\prod_{k=T_0-1}^{\infty}W_{\LL}(k)>\mathbf{0}.$
\end{IEEEproof}

Our next results are aimed at proving that the legitimate agents' values reach consensus almost surely.
In particular, using the decomposition $x_\LL(T_0,t) =\tilde x_\LL(T_0,t) + \phi_\MM(T_0,t)$ of the legitimate agent states 
(see~\eqref{eq:leg-state}), we focus on the limiting behavior of 
$\tilde x_\LL(T_0,t)$, as defined 
in~\eqref{eq:leg-influence}. We show that 
$\tilde x_\LL(T_0,t)$ converges almost surely to a consensus value, as $t\to\infty$, as seen in the following proposition.

\begin{proposition}[Convergence of legitimate agents' values]\label{proposition:convergence_legitimate_part}
Consider the dynamics in~\eqref{eq:dynamics}.
Let $x_{\LL}(0)$ be the vector of initial values of legitimate agents, and
define
\begin{align}
    \tilde{x}_{\LL}^{\infty}(T_0)=\left(\prod_{k=T_0-1}^{\infty}W_{\LL}(k)\right)x_{\LL}(0).
\end{align}
Then, we have $[\tilde{x}_\LL^{\infty}(T_0)]_i=
[\tilde{x}_\LL^{\infty}(T_0)]_j$ for every $i,j\in\LL$,
almost surely. 
\end{proposition}

\begin{IEEEproof}
By Proposition~\ref{proposition:consensus_not_vanishing}, there exists a finite time instant $T_f$ such that  $W_{\LL}(k)=\overline{W}_{\LL}$ for all 
$k\ge T_f$ almost surely, and
$\lim_{k\rightarrow\infty}\overline{W}^k_{\LL}=\1 v'$
for a stochastic vector $v>\mathbf{0}$.
Hence, it follows that almost surely
\begin{align}
\tilde{x}_{\LL}^{\infty}(T_0)
    &=\left(\prod_{k=T_0-1}^{\infty}W_{\LL}(k)\right)x_{\LL}(0) \nonumber\\
    &= 
    \left(\prod_{k=\max\{T_f,T_0\}}^{\infty}W_{\LL}(k)\prod_{k=T_0-1}^{\max\{T_f,T_0\}-1}W_{\LL}(k)\right)x_{\LL}(0)\nonumber\\
    &= 
    \left(\prod_{k=\max\{T_f,T_0\}}^{\infty}\overline{W}_{\LL}\prod_{k=T_0-1}^{\max\{T_f,T_0\}-1}W_{\LL}(k)\right)x_{\LL}(0)\nonumber\\
    &= \boldsymbol{1}v'\left(\prod_{k=T_0-1}^{\max\{T_f,T_0\}-1}W_{\LL}(k)\right)x_{\LL}(0).
\end{align}

By letting
$\tilde{x}_{\LL}^{T_f}(T_0)=\left(\prod_{k=T_0-1}^{\max\{T_f,T_0\}-1}W_{\LL}(k)\right)x_{\LL}(0)$, we obtain that
 $   \tilde{x}_\LL^{\infty}(T_0) = \boldsymbol{1}v'\tilde{x}_{\LL}^{T_f}(T_0)$,
almost surely. 
\end{IEEEproof}

We next consider the limiting effect that malicious agents' states have on the legitimate agents' states. Specifically,
in the state decomposition relation 
$x_\LL(T_0,t) =\tilde x_\LL(T_0,t) + \phi_\MM(T_0,t)$ for the legitimate agents' states 
(see~\eqref{eq:leg-state}), we focus on the limiting behavior of 
$\phi_\MM(T_0,t)$ defined 
in~\eqref{eq:mal-influence}, as $t\to\infty$.
We show that $\phi_\MM(T_0,t)$ converges to a consensus
almost surely given our choice of weight matrices, as seen in the following proposition.
\begin{proposition}\label{proposition:convergence_malicious_part}
Consider the vector $\phi_{\MM}^{\infty}(T_0)$ defined by
\begin{align}
    \phi_{\MM}^{\infty}(T_0) = \sum_{k=T_0-1}^{\infty}\left(\prod_{l=k+1}^{\infty}W_{\mathcal{L}}(l)\right)W_{\MM}(k)x_{\MM}(k).
\end{align}
Then, the entries of this vector are all the same almost surely, i.e., $[\phi_\MM^{\infty}(T_0)]_i=[\phi_\MM^{\infty}(T_0)]_j$ almost surely for every $i,j\in\LL$.
\end{proposition}
\begin{IEEEproof}
By Proposition~\ref{proposition:consensus_not_vanishing}, there exists a finite time instant $T_f>0$ such that
almost surely,
\begin{equation}\label{eq-limit}
W_{\LL}(k)=\overline{W}_{\LL}
\qquad\hbox{for all $k\ge T_f$}.
\end{equation}
This implies that $W_{\MM}(k)=\boldsymbol{0}$ for all 
$k\ge T_f$ almost surely.
Therefore, 
\begin{align}
\phi_\MM^{\infty}(T_0)    &=\sum_{k=T_0-1}^{\infty}\left(\prod_{l=k+1}^{\infty}W_{\mathcal{L}}(l)\right)W_{\MM}(k)x_{\MM}(k) \nonumber\\
    &= \sum_{k=T_0-1}^{T_f-1}\left(\prod_{l=k+1}^{\infty}W_{\mathcal{L}}(l)\right)W_{\MM}(k)x_{\MM}(k).
    \end{align}
Now, for the product 
$\prod_{l=k+1}^{\infty}W_{\mathcal{L}}(l)$,
we have that almost surely
\begin{align*}
\prod_{l=k+1}^{\infty}W_{\mathcal{L}}(l)
&=
\prod_{l=\max\{T_f,k+1\}}^{\infty}W_{\mathcal{L}}(l)
\prod_{l=k+1}^{\max\{T_f,k+1\}-1}W_{\mathcal{L}}(l)\cr
&=\lim_{t\to\infty}\overline{W}_\LL^{t-\max\{T_f,k+1\}}\prod_{l=k+1}^{\max\{T_f,k+1\}-1}W_{\mathcal{L}}(l),
\end{align*}
where the second relation follows from~\eqref{eq-limit}.

By Proposition~\ref{proposition:consensus_not_vanishing},
we also have for a stochastic vector $v>\mathbf{0}$,
such that $\lim_{k\rightarrow\infty}\overline{W}^k_{\LL}=\1 v'$.
Hence, almost surely,
    \begin{align*}
        &\phi_\MM^{\infty}(T_0)
        =\sum_{k=T_0-1}^{T_f-1}\left(\prod_{l=k+1}^{\infty}W_{\mathcal{L}}(l)\right)W_{\MM}(k)x_{\MM}(k)\nonumber\cr
        &\hspace{0.3cm} = \boldsymbol{1}v'
        \left[\sum_{k=T_0-1}^{T_f-1}\left(\prod_{l=k+1}^{\max\{T_f,k+1\}-1}W_{\mathcal{L}}(l)\right)W_{\MM}(k)x_{\MM}(k)\right],
    \end{align*}
    for some finite $T_f$, implying that all the coordinates of the vector $\phi_\MM^{\infty}(T_0)$
    are identical almost surely.
\end{IEEEproof}

As a direct consequence of 
Propositions~\ref{proposition:convergence_legitimate_part} 
and~\ref{proposition:convergence_malicious_part}, we have the following result.
\begin{corollary}
\label{cor:limit}
There exists, almost surely, a random variable  $y(T_0)$ such that 
\begin{align}
    &\lim_{t\rightarrow \infty}\tilde{x}_{\LL}(T_0,t)=y(T_0)\boldsymbol{1},
\end{align}
where $y(T_0)$ is in the convex hull of 
the legitimate agents' initial values \textcolor{black}{$x_i(0)$}, $i\in\LL$ \textcolor{black}{and $0$}, and its distribution depends on the 
starting time $T_0$ of the consensus algorithm in~\eqref{eq:dynamics}. 
Furthermore, almost surely, there exists a random variable $z(T_0)$ such that
\begin{align}
    &\lim_{t\rightarrow \infty} x_{\LL}(T_0,t)=z(T_0)\boldsymbol{1},
\end{align}
where $z(T_0)$ is in the convex hull of the initial legitimate agent values \textcolor{black}{$x_i(0)$}, $i\in\LL$ \textcolor{black}{and the malicious inputs $x_i(t), i\in\MM$, $t\geq T_0-1$ such that malicious agent $i$ is misclassified by a legitimate agent at time $t$}, and its deviation from the nominal consensus value (the case with no malicious agents) depends on the starting time $T_0$ of the consensus algorithm. The analysis of the deviation is given in Theorem \ref{thm:main2}.
%where $z(T_0)$ is in the convex hull of the initial values $x_i(T_0), i=1,\ldots,n$, of the legitimate and malicious agents, and its distribution depends on the  starting time $T_0$ of the consensus algorithm.
\end{corollary}
\begin{IEEEproof}
The results follow directly from Proposition \ref{proposition:convergence_legitimate_part} and Proposition \ref{proposition:convergence_malicious_part},
and the consensus algorithm in~\eqref{eq:dynamics}.
\end{IEEEproof}
We note that, in the view of the fact that the legitimate agents start consensus protocol at time $T_0$ \textcolor{black}{and due to the substochasticity of the matrices $W_{\LL}(t)$}, 
we have that the random variable $y(T_0)$ is in 
the convex hull of the initial values $x_i(0),i\in\LL,$ of legitimate agents \textcolor{black}{and $0$}.
The random variable $z(T_0)$ is in the convex hull 
of $x_i(0),\ i\in\LL,$ and the \textcolor{black}{malicious inputs $x_i(t), i\in\MM$, $t\geq T_0-1$ such that malicious agent $i$ is misclassified by a legitimate agent at time $t$},
since the malicious agents may have deviated from their initial values at time\textcolor{black}{s $t\leq T_f$}.

We are now ready to state our main convergence result.

\begin{theorem}[Convergence of Resilient Consensus]
\label{thm:main1}
Consider the consensus system described by the dynamics in Eq.~\eqref{eq:dynamics} with an arbitrary start time of $T_0\geq 0$ and weights as described by Eq.~\eqref{eq:weights}.  Given that 1) $\G$ is a sufficiently connected network and that 2) bounds on the expected values of the observations $\alpha_{ij}(k)$ (see Defn.~\ref{def:alpha}) satisfy $c=E(\alpha_{ij}(k)-1/2<0$ for all $i\in\LL, j\in\NN_i\cap\MM,k\geq 0$ and $d=E(\alpha_{ij}(k)-1/2>0$ for all 
$i\in\LL, j\in\NN_i\cap \LL, k\geq 0$, 
the consensus algorithm converges almost surely, i.e.,
\begin{align*}
    &\lim_{t\rightarrow \infty} x_{\LL}(T_0,t)=z(T_0)\boldsymbol{1},
\end{align*}
almost surely, independently of the number of malicious agents in the network. The value $z(T_0)$ is in the convex hull of the initial values  $x_i(0),\ i\in\LL,$ and the \textcolor{black}{malicious inputs $x_i(t), i\in\MM$, $t\geq T_0-1$ such that malicious agent $i$ is misclassified by a legitimate agent at time $t$.}
\end{theorem}

\begin{IEEEproof}
The proof follows from Propositions~\ref{proposition:convergence_legitimate_part},~\ref{proposition:convergence_malicious_part}, and Corollary~\ref{cor:limit}.  Specifically, we have that the legitimate agent values satisfy $x_{\mathcal{L}}(T_0,t) = \tilde{x}_{\LL}(T_0,t)+\phi_{\MM}(T_0,t)$, where the first term $\tilde{x}_{\LL}(T_0,t)$ converges by Proposition~\ref{proposition:convergence_legitimate_part} and the second term $\phi_{\MM}(T_0,t)$ converges by Proposition~\ref{proposition:convergence_malicious_part}. 
\end{IEEEproof}

Note that this result asserts convergence in the presence of an \emph{arbitrary number of malicious nodes} and characterizes \emph{necessary conditions on the observations $\alpha_{ij}(t)$} to attain this convergence.  However, the distance of the limit value $z(T_0)$ to the nominal average consensus value in the ideal case (when the malicious agents are known at the start time) has not been discussed.  The characterization of this distance is the subject of the following section.

\vspace{-0.1in}
\subsection{Deviation from Nominal Consensus Value}
In this section we characterize the deviation from the nominal consensus value under our consensus model in the presence of a malicious attack. The nominal consensus value is the consensus value over the graph $\mbbG_\LL$ of legitimate agents in the case of no malicious agents in the network. More specifically, by Proposition~\ref{proposition:consensus_not_vanishing} 
we know that there exists a finite time instant $T_f>0$ such that  $W_{\LL}(k)=\overline{W}_{\LL}$ for all $k\ge T_f$ almost surely, where $\overline{W}_{\LL}$ is given in~\eqref{eq:limitweight}. The nominal consensus value 
corresponds to consensus that would have been reached 
among the legitimate agents if the dynamics in~\eqref{eq:dynamics} were governed by $[\overline{W}_{\LL}\ \mathbf{0}]$ 
instead of $[W_\LL(k)\ W_\MM(k)]$. Hence,
in view of Proposition~\ref{proposition:consensus_not_vanishing}, {\it the nominal consensus value over the graph $\mbbG_\LL$ of legitimate agents is
$\1 v'x_\LL(0)$}.

In what follows we investigate the deviation of  
the consensus value 
$\lim_{t\to\infty} x_\LL(T_0,t)$ reached by the legitimate agents (as predicted by Theorem~\ref{thm:main1})
from the nominal value $\1 v'x_\LL(0)$.
The user-defined start time $T_0$, 
for running the consensus dynamics in Eq.~\eqref{eq:dynamics}, plays a major role in determining the amount of the resulting deviation due to the malicious agents' inputs to the system.  Intuitively, the longer $\alpha_{ij}(t)$ values are observed for all pairs $\{i,j\}\in\mathbb{E}$ before starting the consensus, the probability of making a misclassification error (classifying a legitimate agent as malicious and vice versa) decays exponentially.

Our approach in this section is to show that for larger values of $T_0$, the probability of our weight matrix not having reached its limit value decays exponentially (see Lemma~\ref{lemma:not _equal_T_0_upper_bound}). The implication of this manifests itself as a smaller overall deviation. Specifically, deviation of the convergence limit (cf. Theorem~\ref{thm:main1}) from the nominal average consensus value can be bounded by the deviation introduced by legitimate agents being misclassified as malicious and by malicious agents being misclassified as legitimate.  Propositions~\ref{proposition:upper_legitimate_dev_prob} and~\ref{proposition:upper_malicious_dev_prob} proven in this section provide definitive bounds for these deviations with probability $p(T_0)$ which is characterized as a function of $T_0$.

We begin by showing that the probability of the weights from Eq.~\eqref{eq:weights} not taking their ideal form $\overline{W}_{\LL}$ decays exponentially with $T_0$.
\begin{lemma}\label{lemma:not _equal_T_0_upper_bound}
For every $T_0\geq1,$ 
\begin{align*}
&\Pr\left(\exists k\geq T_0-1\: :\: W_{\LL}(k)\neq \overline{W}_{\LL}\right)\nonumber\\
&\leq
|\LL|^2\cdot\frac{\exp(-2T_0d^2)}{1-\exp(-2d^2)}+|\LL||\MM|\cdot\frac{\exp(-2T_0c^2)}{1-\exp(-2c^2)}.
\end{align*}
\end{lemma}
\begin{IEEEproof}
For every $T_0\geq1$, we have 
\begin{align}
\label{eq-oneo}
&\Pr\left(\exists k\geq T_0-1\: :\: W_{\LL}(k)\neq \overline{W}_{\LL}\right)\nonumber\\
&= Pr\left(\bigcup_{k\geq T_0-1} \{W_{\LL}(k)\neq \overline{W}_{\LL}\}\right)\nonumber\\
&\leq\sum_{k=T_0-1}^{\infty}\Pr\left(W_{\LL}(k)\neq \overline{W}_{\LL}\right),
\end{align}
where the inequality follows from applying the union bound to the event $\bigcup_{k\geq T_0-1} \{W_{\LL}(k)\neq \overline{W}_{\LL}\}$.
Further, for every $k$ we have that 
\begin{align}
\label{eq-twoo}
    \hspace{-0.15cm}\{W_{\LL}(k)\neq \overline{W}_{\LL}\}=\bigcup_{\substack{i\in\LL,\\ j\in\NN_i\cap\LL}}\{\beta_{ij}(k)<0\} \bigcup_{\substack{i\in\LL,\\ j\in\NN_i\cap\MM}}\{\beta_{ij}(k)\geq 0\}.
\end{align}
Therefore, by the union bound 
\begin{align*}
    &\Pr\left(W_{\LL}(k)\neq \overline{W}_{\LL}\right)\cr
    &\leq 
    \sum_{i\in\LL,j\in\NN_i\cap\LL}\Pr\left(\beta_{ij}(k)<0\right) +\sum_{i\in\LL,j\in\NN_i\cap\MM}\Pr\left(\beta_{ij}(k)\geq 0\right)\cr
    &\leq |\LL|^2\cdot\exp(-2(k+1)d^2)+|\LL||\MM|
    \cdot\exp(-2(k+1)c^2),
\end{align*}
where the last inequality follows from  
Lemma~\ref{Lemma:concentration_upper}.
The result then follows by combining relations \eqref{eq-oneo} and~\eqref{eq-twoo}.
\end{IEEEproof}

Next, we consider the deviation contributed by misclassifying legitimate agents as malicious
at some times and, thus, discarding their values. In the sequel, we use the bound of Lemma~\ref{lemma:not _equal_T_0_upper_bound} and the bound we present next in Lemma~\ref{lemma:lemma_upper_diff_sub_matrices} to ultimately
bound the deviation, in Proposition~\ref{proposition:upper_legitimate_dev_prob}, 
that is caused by its incorrect classification of the legitimate agents.

%In what follows, we make use of the following lemma.
\begin{lemma}\label{lemma:lemma_upper_diff_sub_matrices}
Let $\tilde{X},X\in\mathbb{R}_+^{m\times m}$ be two row-substochastic matrices, and let $\gamma>0$ be such that $\tilde{X}_{ii}\geq\gamma>0$ and $X_{ii}\geq \gamma$ for all $i$. Then $\left[\left|\tilde{X}-X\right|\boldsymbol{1}\right]_i\leq 2(1-\gamma)$,
for all $i$, where $|X|$ denote the matrix with entries $|x_{ij}|$. 
\end{lemma}
\begin{IEEEproof}
Let coordinate index $i$ be arbitrary, $1\le i\le m$.
First, we prove that $\left|\tilde{x}_{ii}-x_{ii}\right|+2\gamma\leq \tilde{x}_{ii}+ x_{ii}$.
By the triangle inequality, we obtain
\begin{align*}
    \left|\tilde{x}_{ii}-x_{ii}\right|+2\gamma &= \left|(\tilde{x}_{ii}-\gamma)-(x_{ii}-\gamma)\right|+2\gamma\nonumber\\
    &=\left|\tilde{x}_{ii}-\gamma\right|+\left|x_{ii}-\gamma\right|+2\gamma.
\end{align*}
Since $\tilde{x}_{ii}-\gamma\geq0$ and $x_{ii}-\gamma\geq0$ we have that 
\begin{align*}
    \left|\tilde{x}_{ii}-\gamma\right|+\left|x_{ii}-\gamma\right|+2\gamma = 
    \tilde{x}_{ii} + x_{ii}.
\end{align*}
Therefore, it follows that
\begin{align*}
     \left[\left|\tilde{X}-X\right|\boldsymbol{1}\right]_i+2\gamma 
     &=  
     \sum_{j=1}^m\left|\tilde{x}_{ij}-x_{ij}\right|+2\gamma\nonumber\\
     &=\tilde{x}_{ii} + x_{ii}
     +\sum_{j=1,\, j\ne i}^{m}
     \left|\tilde x_{ij}-x_{ij} \right|\cr
     &\le \tilde{x}_{ii} + x_{ii}
     +\sum_{j=1,\, j\ne i}^{m}
     \left(\tilde{x}_{ij} + x_{ij}
     \right),
     %&\leq \sum_{j\in %|\LL|}\left(\tilde{x}_{ij} + x_{ij}
     %\right).
\end{align*}
where the last inequality is obtained 
by using the triangle inequality and the fact that the matrices $X$ and $\tilde X$ have nonnegative entries. 
Since both $\tilde{X}$ and $X$ are stochastic matrices, we further have that
$\sum_{j=1}^{m}\left[\tilde{x}_{ij}+x_{ij}\right]\leq 2$.
From the preceding two relations, it follows that
$\left[\left|\tilde{X}-X\right|\boldsymbol{1}\right]_i\leq 2(1-\gamma)$.
\end{IEEEproof}

Denote by $\varphi_i(T_0,t)$ the deviation 
suffered by legitimate agent $i$ that is caused by its incorrect classification of the legitimate agents, i.e., for all $i\in\LL$,
\begin{equation}
\label{eq:devleg}
\varphi_i(T_0,t)=\left|\left[ \tilde{x}_{\LL}(T_0,t)- \left(\prod_{k=T_0}^{t-1}\overline{W}_{\LL}\right)x_{\LL}(0)\right]_i\right|.
\end{equation}
Then we have the following probabilistic bound on the deviation.
\begin{proposition}\label{proposition:upper_legitimate_dev_prob}
Given an error tolerance $\delta>0$,
for the deviations $\varphi_i(T_0,t)$, $i\in\LL$, as defined in~\eqref{eq:devleg},
we have
\begin{align*}
    \Pr\left(\max_{i\in\LL}\,\limsup_{t\to\infty}\varphi_{i}(T_0,t)> \tilde{g}_{\LL}(T_0,\delta)\right)<\delta,
\end{align*}
where 
\begin{align}
\label{eq:boundgleg}
\tilde{g}_{\LL}(T_0,\delta)&=\frac{2\eta |\LL|^2}{\delta}\cdot\frac{\exp(-2T_0d^2)}{1-\exp(-2d^2)}\nonumber\\
&\qquad+\frac{2\eta |\LL| \,|\MM|}{\delta}\cdot\frac{\exp(-2T_0c^2)}{1-\exp(-2c^2)},
\end{align}
and $\eta\geq\sup_{i\in\LL\cup\MM,t\in\mathbb{N}}|x_i(t)|$ is a finite upper bound on the absolute input value.
\end{proposition}
\begin{IEEEproof}
Let $T_f(T_0,t)$ be a random variable equal to $0$ if $W_{\LL}(k)=\overline{W}_{\LL}$ for all $k\in [T_0-1,t-1]$, and equal to $\sup\{l+1:W_{\LL}(l+T_0-1)\textcolor{black}{\neq}\overline{W}_{\LL},\:l\in[0,t-T_0]\}$, otherwise.
In view of the evolution of 
$\tilde{x}_{\LL}(T_0,t)$ as given in~\eqref{eq:leg-influence}, we have
\[
\varphi_i(T_0,t)
= \left|\left[
\left(\left(\prod_{k=T_0-1}^{t-1}W_{\LL}(k)\right)
-\overline{W}_{\LL}^{t-T_0}\right)
x_{\LL}(0)\right]_i\right|.
\]

Define
\begin{align*}
&\hspace{-0.2cm}\Delta(W_\LL,T_f)=\left(\prod_{k=T_0-1}^{T_0+T_f(T_0,t)-2}W_{\LL}(k)\right)-\left(\prod_{k=T_0-1}^{T_0+T_f(T_0,t)-2}\overline{W}_{\LL}\right).
\end{align*}
Then,
\begin{align}
    \varphi_i(T_0,t)&=\left|\left[\left(\prod_{k=T_0+T_f(T_0,t)-1}^{t-1}\overline{W}_{\LL}\right)\Delta(W_{\LL},T_f)x_{\LL}(0)\right]_{i}\right|\nonumber\\
&\stackrel{(a)}{\leq} 
\max_{i\in\LL} \left|\left[\Delta(W_{\LL},T_f)x_{\LL}(0)
\right]_{i}\right|,
\end{align}
where $(a)$ follows since $\overline{W}_{\LL}$ is a row-stochastic matrix.
Further, since $|x_\LL(0)|\le\eta$, it follows that
\[\max_{i\in\LL} \left|\left[\Delta(W_{\LL},T_f)x_{\LL}(0)\right]_{i}\right|\le \eta \max_{i\in\LL} \left[\left|\Delta(W_{\LL},T_f)\right|
\boldsymbol{1}\right]_{i},\]
where for a matrix $A=[a_{ij}]$, we use $|A|$ to denote the matrix with entries $|a_{ij}|$.
Therefore,
\[\varphi_i(T_0,t)\le \eta \max_{i\in\LL} \left[\left|\Delta(W_{\LL},T_f)\right|
\boldsymbol{1}\right]_{i}.\]

Let $n_w=\max\{|\LL|+|\MM|,\kappa\}$
and note that $\overline{W}_{\LL}$ is a stochastic matrix with $[\overline{W}_{\LL}]_{ii}\geq\frac{1}{n_w}$ (see~\eqref{eq:limitweight}). Additionally, for every $k$, the matrix $W_{\LL}(k)$ is substochastic with $[W_{\LL}(k)]_{ii}
\geq\frac{1}{n_w}$  for every $i\in\LL$. Therefore, by Lemma \ref{lemma:lemma_upper_diff_sub_matrices} we have 
\begin{align*}
 \left[\left|\Delta(W_{\LL},T_f)\right|\1\right]_{i}
 \leq 2\left[1-\left(\frac{1}{n_w}\right)^{T_f(T_0,t)}\right].
\end{align*}
Define
\[\overline{\varphi}(T_0,t)= \left[1-\left(\frac{1}{n_{w}}\right)^{T_f(T_0,t)}\right]2\eta.\]
Then, we have
$\varphi_i(T_0,t)\leq \overline{\varphi}(T_0,t)$ for all $i\in\LL$ and $t\geq T_0$, thus implying that
\[\max_{i\in\LL}\,\limsup_{t\rightarrow\infty}\varphi_i(T_0,t)\leq \lim_{t\rightarrow\infty}\overline{\varphi}(T_0,t).\]
By Proposition~\ref{proposition:consensus_not_vanishing},
for $t\ge T_f$, we have $T_f(T_0,t)=t$ almost surely,
implying that $\overline{\varphi}(T_0,t)$ almost surely converges as $t\to\infty$. 
By Markov's inequality it follows that
\begin{align*}
\Pr\left(\lim_{t\rightarrow\infty}\overline{\varphi}(T_0,t)> \tilde{g}_{\LL}(T_0,\delta)\right)
\leq \frac{E\left(\lim_{t\rightarrow\infty}\overline{\varphi}(T_0,t)\right)}{\tilde{g}_{\LL}(T_0,\delta)}.
\end{align*}
By definition, for every $t\geq T_0$ we have that  $0\leq T_f(T_0,t)\leq T_f(T_0,t+1)$, therefore, 
$0\leq\overline{\varphi}(T_0,t)\leq\overline{\varphi}(T_0,t+1)\leq 2$ for every $t\geq T_0$. Hence, by the Monotone Convergence Theorem (see \cite{durrett_2010}), it follows that
\begin{align*}
E\left(\lim_{t\rightarrow\infty}\overline{\varphi}(T_0,t)
\right) 
    &= \lim_{t\rightarrow\infty} E\left(\overline{\varphi}(T_0,t)\right)\nonumber\\
    &=2\eta\left[1-\lim_{t\rightarrow\infty} E\left(\left(\frac{1}{n_{w}}\right)^{T_f(T_0,t)}\right)\right]\nonumber\\
    &\leq 2\eta\left[1-\lim_{t\rightarrow\infty} \Pr\left(T_f(T_0,t)=0\right)\right]\cr
    &= 2\eta\lim_{t\rightarrow\infty}\left[1- \Pr\left(T_f(T_0,t)=0\right)\right].
\end{align*}    
Since $1- \Pr\left(T_f(T_0,t)=0\right)= \Pr\left(T_f(T_0,t)>0\right)$,
it further follows that
\begin{align*}
E\left(\lim_{t\rightarrow\infty}\overline{\varphi}(T_0,t)
\right) 
    &= 2\eta\lim_{t\rightarrow\infty} \Pr\left(T_f(T_0,t)>0\right)\cr
    &\leq 
    2\eta\Pr\left(\exists k\geq T_0: W_{\LL}(k)\neq \overline{W}_{\LL}\right)\cr
    &\leq\delta\cdot\tilde{g}_{\LL}(T_0,\delta),
\end{align*} 
where the last inequality follows from 
Lemma~\ref{lemma:not _equal_T_0_upper_bound}
and the definition of $\tilde{g}_{\LL}(T_0,\delta)$
in~\eqref{eq:boundgleg}.
\end{IEEEproof}

Now we consider the deviation contributed by malicious agents that are misclassified as legitimate.
We denote by $\phi_{i}(T_0,t)$, $i\in\LL$, the worst case effect on legitimate agent $i$
due to incorrect malicious agents' classification (labeling an untrustworthy agent as trustworthy) \textcolor{black}{committed by any legitimate agent at some time $t\geq T_0-1$.}  Specifically, let $\eta=\max_{i\in\LL\cup\MM}|[x_\LL(0)]_i|$
and define 
\begin{align}
\label{eq:devmal}
    \phi_{i}(T_0,t) = 
    \eta\sum_{k=T_0-1}^{t-1}\sum_{j\in \textcolor{black}{\MM}}
    \left[\left(\prod_{l=k+1}^{t-1}W_{\LL}(l)\right)W_{\MM}(k)\right]_{ij},
\end{align}
for all $i\in\LL$.
We note that these quantities are nonnegative, since the matrices $W_{\LL}(k)$ and $W_{\MM}(k)$ are nonnegative for all $k$.
Looking at the malicious agents' influence vector $\phi_\MM(T_0,t)$ defined in~\eqref{eq:mal-influence},
we see that
\begin{align*}
    |\left[\phi_{\MM}(T_0,t)\right]_i|\leq \max_{j\in\LL}\phi_j(T_0,t), \:\:\forall\: i\in\LL.
\end{align*}
Our next result provides a probabilistic bound
on the maximal influence $\max_{j\in\LL}\phi_j(T_0,t)$
of malicious agents inputs on the legitimate agents' values.
To do so, given $\delta>0$, we define
\begin{align}
\label{eq:boundgmal}
\tilde{g}_\MM(T_0,\delta)
=\frac{\eta |\LL||\MM|}{\delta\cdot\kappa}\cdot\frac{\exp(-2T_0c^2)}{1-\exp(-2c^2)},
\end{align}
where $\eta\geq\sup_{i\in\LL\cup\MM,t\in\mathbb{N}}|x_i(t)|$.

\begin{proposition}\label{proposition:upper_malicious_dev_prob}
Given an error tolerance $\delta>0$,
for the deviations $\phi_i(T_0,t)$, $i\in\LL$, as defined in~\eqref{eq:devmal},
we have
\begin{align}
    \Pr\left(\max_{i\in\LL}\,
    \limsup_{t\rightarrow\infty}\phi_{i}(T_0,t)> \tilde{g}_{\MM}(T_0,\delta)\right)<\delta,
\end{align}
where $\tilde{g}_{\MM}(T_0,\delta)$ is given 
by~\eqref{eq:boundgmal}.
\end{proposition}

\begin{IEEEproof}
\textcolor{black}{By the substochasticity of $W_{\LL}(k)$ and the nonnegativity of $W_{\MM}(k)$ for all $k\geq0$, we have for every $i\in\LL$ and $t\geq T_0\geq1$,
\begin{flalign}\label{eq:phi_i_upper_bound}
     \phi_{i}(T_0,t)&\leq
     \eta\sum_{k=T_0-1}^{t-1}\sum_{i\in\LL}\sum_{j\in\MM}\left[W_{\MM}(k)\right]_{ij}\leq  \overline{\phi}_{i}(T_0,t),
\end{flalign}
where
\begin{flalign}
    \overline{\phi}_{i}(T_0,t) \triangleq  \eta\sum_{k=T_0-1}^{t-1}\sum_{\ell\in\LL}\sum_{j\in\MM} \frac{1}{\kappa}\cdot \mathbbm{1}_{\{\beta_{\ell j}(k)\geq 0\}}.
\end{flalign}
Since the definition of $\overline{\phi}_{i}(T_0,t)$ is identical for all $i\in\LL$, we denote $\overline{\phi}(T_0,t)\equiv \overline{\phi}_{i}(T_0,t)$. }

\textcolor{black}{Observe that by \eqref{eq:phi_i_upper_bound} we have that
\begin{flalign}
\max_{i\in\LL}\limsup_{t\rightarrow\infty}\phi_{i}(T_0,t)
&\leq \limsup_{t\rightarrow\infty}\max_{i\in\LL}\phi_{i}(T_0,t)\nonumber\\
&\leq \limsup_{t\rightarrow\infty}\overline{\phi}(T_0,t).
\end{flalign}
Consequently,
\begin{align*}
&\Pr\left(\max_{i\in\LL}\,\limsup_{t\rightarrow\infty}
    \phi_{i}(T_0,t)> \tilde{g}_{\MM}(T_0,\delta)\right)\\
    &\leq \Pr\left(\limsup_{t\rightarrow\infty}\overline{\phi}(T_0,t)> \tilde{g}_{\MM}(T_0,\delta)\right).
\end{align*}
We can now apply Markov's inequality on the right-hand-side to deduce that
\begin{align*}
\Pr\left(\limsup_{t\rightarrow\infty}\overline{\phi}(T_0,t)> \tilde{g}_{\MM}(T_0,\delta)\right)\leq \frac{E\left[\limsup_{t\rightarrow\infty}\overline{\phi}(T_0,t)\right]}{\tilde{g}_{\MM}(T_0,\delta)}.
\end{align*}
It is easy to verify by definition that $\overline{\phi}(T_0,t)$ is a monotonically non-decreasing function of $t$, for every $T_0$.
Additionally, for any realization of trust observations $\{\beta_{ij}(k)\}_{i\in\LL,j\in\MM,k\geq T_0-1}$ we have the pointwise convergence
\[\lim_{t\rightarrow\infty}\overline{\phi}(T_0,t)=\frac{\eta}{\kappa}\sum_{k=T_0-1}^{\infty}\sum_{i\in\LL}\sum_{j\in\MM} \mathbbm{1}_{\{\beta_{ij}(k)\geq 0\}},\]
note that this limit may be $\infty$. %\textcolor{red}{[(MY) I added this part since we have to show point-wise convergence in all cases.]}.
Thus, by the Monotone Convergence Theorem \cite{durrett_2010} we can replace the order of expectation and $\limsup_{t\rightarrow\infty}$ as follows:  
  \begin{align}
  \label{eq:five}
      E\left(\limsup_{t\to\infty}
      \overline{\phi}(T_0,t)\right)
      &=E\left(\lim_{t\to\infty}\overline{\phi}(T_0,t)\right)\cr
  &=\lim_{t\to\infty}E\left(\overline{\phi}(T_0,t)\right).
  \end{align}
Now, 
\begin{flalign}
E\left(\overline{\phi}(T_0,t)\right) &= \frac{\eta}{\kappa}\sum_{k=T_0-1}^{t-1}\sum_{\ell\in\LL}\sum_{j\in\MM} E\left(\mathbbm{1}_{\{\beta_{\ell j}(k)\geq 0\}}\right)\nonumber\\
&\leq \frac{\eta}{\kappa}\sum_{k=T_0-1}^{t-1}\sum_{\ell\in\LL}\sum_{j\in\MM} E\left(\mathbbm{1}_{\{\beta_{\ell j}(k)\geq 0\}}\right)\nonumber\\
&= \frac{\eta}{\kappa}\sum_{k=T_0-1}^{t-1}\sum_{\ell\in\LL}\sum_{j\in\MM} \Pr\left(\beta_{\ell j}(k)\geq0\right)\nonumber\\
&\leq \frac{\eta}{\kappa}\sum_{k=T_0-1}^{t-1}\sum_{\ell\in\LL}\sum_{j\in\MM} \exp(-2(k+1)c^2)\nonumber\\
&\leq \frac{\eta|\LL||\MM|}{\kappa}\sum_{k=T_0-1}^{\infty} \exp(-2(k+1)c^2)\nonumber\\
&\leq \frac{\eta|\LL||\MM|}{\kappa}\cdot\frac{\exp(-2T_0c^2)}{1-\exp(-2c^2)}.
\end{flalign}
Thus, 
\begin{align}
    \Pr\left(\max_{i\in\LL}\,
    \limsup_{t\rightarrow\infty}\phi_{i}(T_0,t)> \tilde{g}_{\MM}(T_0,\delta)\right)<\delta,
\end{align}
where 
\begin{align}
\label{eq:boundgmal}
\tilde{g}_\MM(T_0,\delta)
=\frac{\eta |\LL||\MM|}{\delta \kappa}\cdot\frac{\exp(-2T_0c^2)}{1-\exp(-2c^2)},
\end{align}
as Proposition 5 states.}

\end{IEEEproof}

We are now ready to prove our main deviation result which states that the maximum deviation of the converged consensus value following the dynamics in Eq.~\eqref{eq:dynamics} with the weights described in Eq.~\eqref{eq:weights} can be bounded by a finite value $\Delta_\text{max}(T_0, \delta)$ that we characterize.
\begin{theorem}[Deviation from Nominal Average Consensus]
\label{thm:main2}
Consider the consensus system described by the dynamics in Eq.~\eqref{eq:dynamics} with an arbitrary start time of $t=T_0\geq 0$ and weights as described by Eq.~\eqref{eq:weights}. Given that 1) $\G$ is a sufficiently connected network and that 2) bounds on the expected values of the observations $\alpha_{ij}(k)$ (see Defn.~\ref{def:alpha}) satisfy $c=E(\alpha_{ij}(k)-1/2<0$ for all 
$i\in\LL, j\in\NN_i\cap \MM, k\geq 0$, and $d=E(\alpha_{ij}(k)-1/2>0$ 
for all $i\in\LL, j\in\NN_i\cap \LL, k\geq 0$,
for a given error probability $\delta>0$, we have the following result
\begin{align*}
&\Pr \left(\max_{i\in\LL}\,\limsup_{t\rightarrow \infty} \left|\left[ x_{\LL}(T_0,t)-\boldsymbol{1}v'x_{\LL}(0)\right]_i \right| \leq \Delta_\text{max}(T_0, \delta) \right)\nonumber\\
&\qquad \ge 1 - \delta, 
\end{align*}
where
\begin{align}\label{eq:def:g:sum_L_M}
    \Delta_\text{max}(T_0,\delta) = 2\left[\tilde{g}_{\LL}(T_0,\delta)+\tilde{g}_{\MM}(T_0,\delta)\right],
\end{align}
with $\tilde{g}_{\LL}(T_0,\delta)$ and $\tilde{g}_\MM(T_0,\delta)$ respectively given by
\begin{align*}
&\tilde{g}_{\LL}(T_0,\delta)\nonumber\\
&=\frac{\eta |\LL|^2}{\delta}\cdot\frac{\exp(-2T_0d^2)}{1-\exp(-2d^2)}+\frac{\eta |\LL| |\MM|}{\delta}\cdot\frac{\exp(-2T_0c^2)}{1-\exp(-2c^2)},
\end{align*}
and
\begin{align*}
&\tilde{g}_{\MM}(T_0,\delta)=\frac{\eta |\LL| |\MM|}{\delta\cdot \kappa}\cdot\frac{\exp(-2T_0c^2)}{1-\exp(-2c^2)}.
\end{align*}

\end{theorem}
\begin{IEEEproof}
First, note that by the triangle inequality
\begin{align*}
    &\left| \left[x_{\LL}(T_0,t)-\boldsymbol{1}v'x_{\LL}(0)\right]_i \right| \nonumber\\
    &\leq \left|\left[x_{\LL}(T_0,t)-\left(\prod_{k=T_0}^{t-1}\overline{W}_{\LL}\right)x_{\LL}(0)\right]_i\right|\nonumber\\
    &\qquad+\left|\left[\left(\prod_{k=T_0}^{t-1}\overline{W}_{\LL}\right)x_{\LL}(0)-\boldsymbol{1}v'x_{\LL}(0)\right]_i \right|.
\end{align*}
Since $\overline{W}_{\LL}$ is a primitive stochastic matrix, by the Perron-Frobenius Theorem it follows that
\[\lim_{t\to\infty}\left|\left[\left(\prod_{k=T_0}^{t-1}\overline{W}_{\LL}\right)x_{\LL}(0)-\boldsymbol{1}v'x_{\LL}(0)\right]_i \right|=0.\]
Thus, 
\begin{align*}
    &\Pr \left(\max_{i\in\LL}\,\limsup_{t\rightarrow \infty} \left|\left[ x_{\LL}(T_0,t)-\boldsymbol{1}v'x_{\LL}(0)\right]_i \right| \geq 
    \Delta_\text{max}(T_0,\delta)
    %g(\mathcal{P},T_0, \delta) 
    \right)\nonumber\\
    &=\Pr \Bigg(\max_{i\in\LL}\,\limsup_{t\rightarrow \infty} \left| \left[x_{\LL}(T_0,t)- \left(\prod_{k=T_0}^{t-1}\overline{W}_{\LL}\right)x_{\LL}(0)\right]_i\right| \nonumber\\
    &\hspace{6cm}\geq \Delta_\text{max}(T_0,\delta)
    %g(\mathcal{P},T_0, \delta) 
    \Bigg).
\end{align*}
Now, since $x_{\mathcal{L}}(T_0,t) = \tilde{x}_{\LL}(T_0,t)+\phi_{\MM}(T_0,t)$, by the triangle inequality we have
\begin{align*}
    &\left|\left[ x_{\LL}(T_0,t)- \left(\prod_{k=T_0}^{t-1}\overline{W}_{\LL}\right)x_{\LL}(0)\right]_i\right|\nonumber\\
    &\leq \left|\left[ \tilde{x}_{\LL}(T_0,t)- \left(\prod_{k=T_0}^{t-1}\overline{W}_{\LL}\right)x_{\LL}(0)\right]_i\right|
   + \left| \left[\phi_{\MM}(T_0,t)\right]_i\right|.
\end{align*}
It follows by the definition of $\Delta_\text{max}(T_0,\delta)$
in~\eqref{eq:def:g:sum_L_M} and the union bound that
\begin{flalign*}
&\Pr \left(\max_{i\in\LL}\,\limsup_{t\rightarrow \infty} \left |\left[ x_{\LL}(T_0,t)-\boldsymbol{1}v'x_{\LL}(0)\right]_i \right| \geq \Delta_\text{max}(T_0,\delta)
\right)\nonumber\\
&\leq\Pr \left(\max_{i\in\LL}\,\limsup_{t\rightarrow \infty} \left|\left[\tilde{x}_{\LL}(T_0,t)- \Bigg(\prod_{k=T_0}^{t-1}\overline{W}_{\LL}\right)x_{\LL}(0)\right]_i\right| \nonumber\\
&\hspace{5cm}\geq 
2\tilde g_{\LL}(T_0, \delta) \Bigg)\nonumber\\
&\quad+\Pr \left(\max_{i\in\LL}\limsup_{t\rightarrow \infty} \left|\left[ \phi_{\MM}(T_0,t)\right]_i\right|\geq 
2\tilde g_{\MM}(T_0, \delta) \right).
\end{flalign*}
Since $2\tilde g_\LL(T_0, \delta)=
\tilde g_\LL(T_0, \delta/2)$ and 
$2\tilde g_\MM(T_0, \delta)=
\tilde g_\MM(T_0, \delta/2)$, 
it follows that
\begin{flalign*}
&\Pr \left(\max_{i\in\LL}\,\limsup_{t\rightarrow \infty} \left |\left[ x_{\LL}(T_0,t)-\boldsymbol{1}v'x_{\LL}(0)\right]_i \right| \geq \Delta_\text{max}(T_0,\delta)\right)\nonumber\\
&\leq\Pr \left(\max_{i\in\LL}\,\limsup_{t\rightarrow \infty} \left|\left[\tilde{x}_{\LL}(T_0,t)- \Bigg(\prod_{k=T_0}^{t-1}\overline{W}_{\LL}\right)x_{\LL}(0)\right]_i\right| \nonumber\\
&\hspace{5cm}\geq 
\tilde g_{\LL}(T_0, \delta/2) \Bigg)\nonumber\\
&\quad+\Pr \left(\max_{i\in\LL}\limsup_{t\rightarrow \infty} \left|\left[ \phi_{\MM}(T_0,t)\right]_i\right|\geq 
\tilde g_{\MM}(T_0, \delta/2) \right).
\end{flalign*}
Finally, by 
Proposition~\ref{proposition:upper_legitimate_dev_prob} we have
\begin{align*}
    \frac{\delta}{2}\geq&\Pr \left(\max_{i\in\LL}\,\limsup_{t\rightarrow \infty} \left|\left[ \tilde{x}_{\LL}(T_0,t)- \Bigg(\prod_{k=T_0}^{t-1}\overline{W}_{\LL}\right)x_{\LL}(0)\right]_i\right| \nonumber\\
&\hspace{5cm}\geq 
\tilde g_{\LL}(T_0, \delta/2) \Bigg),
\end{align*}  
while, similarly, by 
Proposition~\ref{proposition:upper_malicious_dev_prob}
we have
\begin{align*}
    \frac{\delta}{2}\geq\Pr \left(\max_{i\in\LL}\,\limsup_{t\rightarrow \infty} \left|\left[ \phi_{\MM}(T_0,t)\right]_i\right|\geq 
    \tilde g_{\MM}(T_0, \delta/2) \right).
\end{align*}
\end{IEEEproof}

Thus, we have solved Problem~\ref{sec:problem2} by showing that beyond achieving convergence of consensus in the face of a malicious attack, the deviation of the converged value from the nominal average consensus value (in the case of no malicious agents) can be characterized as a function of user-defined parameters $T_0$ and failure probability $\delta$. In the next section we characterize convergence rate of our consensus protocol.

\vspace{-0.1in}
\subsection{Convergence Rate of Resilient Consensus}

In this section we discuss convergence rate for the consensus protocol in Eq.~\eqref{eq:dynamics} using the weights as defined in Eq.~\eqref{eq:weights}. Let
\[\|z\|_v=\sqrt{\sum_{i=1}^{n_\LL}v_iz_i^2}.\]
We start by presenting a useful lemma showing a contraction property for a consensus step using a primitive matrix that we will later employ. 

\begin{lemma}\label{Lemma_reversible_weight_matrix_limit}
Assume that  $j\in \mathcal{N}_i$ if and only if $i\in \mathcal{N}_j$ for every $i,j\in\LL$ (i.e., bidirectional links between legitimate agents).
Let $\rho_2$ denote the second largest eigenvalue modulus\footnote{The second largest absolute value of the eigenvalues.} of  $\overline{W}_{\LL}$ and let $v>\mathbf{0}$
be the stochastic Perron vector satisfying $v'\overline{W}_{\LL}=v'$. Then,
\[\|\overline{W}_{\LL}^tx(0)-\boldsymbol{1}v'x(0)\|_v\leq\rho_2^t\cdot\|x(0)-\boldsymbol{1}v'x(0)\|_v.\]
\end{lemma}
\begin{IEEEproof}
It can be seen that the Perron vector $v$ has entries $v_i=\frac{\max\{|\mathcal{N}_i\cap\LL|+1,\kappa\}}{\sum_{j\in\LL}\max\{|\mathcal{N}_j\cap\LL|+1,\kappa\}}$ for all $i$, implying that $(\overline{W}_{\LL},v)$ is reversible Markov chain. Thus, the result follows from the convergence rate of reversible Markov chains (see \cite{markov_convergence_rate}).
\end{IEEEproof}

We are now ready to prove our convergence rate results. 
Namely, we show that the rate of convergence of legitimate agents to their limit value is governed by the second largest eigenvalue modulus of the ideal weight matrix $\overline{W}_{\LL}$ with a probability that approaches 1 exponentially, as $T_0$ increases.

\begin{theorem}[Convergence Rate of Resilient Consensus]\label{thm:main3}
  Consider the consensus system described by the dynamics in Eq.~\eqref{eq:dynamics} with an arbitrary start time of $T_0\geq 0$ and weights as described by Eq.~\eqref{eq:weights}.  Assume that 1) $\G$ is a sufficiently connected network and that 2) bounds on the expected values of the observations $\alpha_{ij}$ (see Defn.~\ref{def:alpha}) satisfy $c=E(\alpha_{ij}(k))-1/2<0$ 
for all $i\in\LL, j\in\NN_i\cap\MM$, for all $k\geq 0,$ and $d=E(\alpha_{ij}(k))-1/2>0$ for all $i\in\LL, j\in\NN_i\cap\LL$, for all $k\geq 0$.
Let $\rho_2$ denote the second largest eigenvalue modulus of  $\overline{W}_{\LL}$ and let $v>\mathbf{0}$
be the stochastic Perron vector satisfying $v'\overline{W}_{\LL}=v'$. Additionally, let
$\eta\geq\sup_{i\in\LL\cup\MM,t\in\mathbb{N}}|x_i(t)|$ and assume that  $j\in \mathcal{N}_i$ if and only if $i\in \mathcal{N}_j$ for every $i,j\in\LL$.
For every $T_0>0$, if we start the consensus protocol from time $T_0$, then 
for every $t\geq T_0$ and $m=T_0-1,\ldots,t-1$, we have
\begin{align}
    \left\|x_{\LL}(T_0,t)-z(T_0)\boldsymbol{1}\right\|_v
        \leq 2(m-T_0+2)\rho_2^{t-m}\eta, 
\end{align}
with a probability greater than
\begin{flalign}\label{eq:upper_bound_no_error}
1-|\LL|^2\cdot\frac{\exp(-2(m+1)d^2)}{1-\exp(-2d^2)}-|\LL||\MM|\cdot\frac{\exp(-2(m+1)c^2)}{1-\exp(-2c^2)}.
\end{flalign}
\end{theorem}
\begin{IEEEproof}
Recall that $\lim_{t\rightarrow\infty}x_{\LL}(T_0,t)=z(T_0)\boldsymbol{1}$ almost surely.
Since $x_{\LL}(T_0,t)=\tilde{x}_{\LL}(T_0,t)+\phi_{\MM}(T_0,t)$, it follows by the triangle inequality that
\begin{flalign*}
  & \left\|x_{\LL}(T_0,t)-z(T_0)\boldsymbol{1}\right\|_v\nonumber\\
& =  \left\|x_{\LL}(T_0,t)-\lim_{\tau\rightarrow\infty}x_{\LL}(T_0,\tau)\right\|_v\nonumber\\
&\leq \left\|\tilde{x}_{\LL}(T_0,t)-\lim_{\tau\rightarrow\infty}\tilde{x}_{\LL}(T_0,\tau)\right\|_v\nonumber\\
&\quad+\left\|\phi_{\MM}(T_0,t)-\lim_{\tau\rightarrow\infty}\phi_{\MM}(T_0,\tau)\right\|_v.
\end{flalign*}
Now, assume that for all  $k\geq m$, $W_{\LL}(k)=\overline{W}_{\LL}$ and $W_{\MM}(k)=\boldsymbol{0}$. Then, by the Perron-Frobenius Theorem and the definition of $\eta$,
\begin{flalign*}
&\left\|\tilde{x}_{\LL}(T_0,t)-\lim_{\tau\rightarrow\infty}\tilde{x}_{\LL}(T_0,\tau)\right\|_v\nonumber\\
&=\left\|
\left(\overline{W}_{\LL}^{t-m}
-\boldsymbol{1}v\right)
\left(\prod_{k=T_0-1}^{m-1}W_{\LL}(k)\right)x_{\LL}(0)\right\|_v\nonumber\\
&\leq \rho_2^{t-m}2\eta.
\end{flalign*}
Additionally,
\begin{flalign*}
&\left\|\phi_{\MM}(T_0,t)-\lim_{\tau\rightarrow\infty}\phi_{\MM}(T_0,\tau)\right\|_v\nonumber\\
&=\left\|\sum_{k=T_0-1}^{m-1}\left(\prod_{l=k+1}^{t-1}W_{\mathcal{L}}(l)\right)W_{\MM}(k)x_{\MM}(k)\right.\nonumber\\
&\quad\left.-\lim_{\tau\rightarrow\infty}\sum_{k=T_0-1}^{m-1}\left(\prod_{l=k+1}^{\tau-1}W_{\mathcal{L}}(l)\right)W_{\MM}(k)x_{\MM}(k)\right\|_v\nonumber\\
&\hspace{-0.5cm}=\left\|\sum_{k=T_0-1}^{m-1}\left(\overline{W}_{\LL}^{t-m}
-\boldsymbol{1}v\right)
\left(\prod_{k=T_0-1}^{m-1}W_{\LL}(k)\right)W_{\MM}(k)x_{\MM}(k)\right\|_v\nonumber\\
&\leq (m-T_0+1)\rho_2^{t-m}2\eta.
\end{flalign*}
Finally, following the proof of Lemma \ref{lemma:not _equal_T_0_upper_bound},  we can lower bound the  probability of the event that  $W_{\LL}(k)=\overline{W}_{\LL}$ and $W_{\MM}(k)=\boldsymbol{0}$ for all  $k\geq m$ by \eqref{eq:upper_bound_no_error}.
\end{IEEEproof}
As an immediate consequence of Theorem~\ref{thm:main3},
we have the following result.
\begin{corollary}\label{corollry_con_rate}
Under the conditions of Theorem~\ref{thm:main3},
for every $T_0\ge 1$ and $t\ge T_0$
 we have
\begin{align*}
       & E\left(\left\| x_{\LL}(T_0,t)-z(T_0)\boldsymbol{1}\right\|_v
       \right)\nonumber\\
       &\leq \min_{m\in\{T_0-1,\ldots,t-1\}}\bigg\{2(m-T_0+2)\rho_2^{t-m} \eta+2\eta \cdot\cr 
        & \hspace{0.2cm}\left(\frac{|\LL|^2\exp(-2(m+1)d^2)}{1-\exp(-2d^2)}+\frac{|\LL||\MM|\exp(-2(m+1)c^2)}{1-\exp(-2c^2)}\right)\bigg\}. 
\end{align*}
\end{corollary} 
Note that choosing $m=\frac{t+T_0}{2}$ yields the bound
\begin{align*}
       & E\left(\left\| x_{\LL}(T_0,t)-z(T_0)\boldsymbol{1}\right\|_v
       \right)\nonumber\\
       &\leq 2\left(\frac{t-T_0}{2}+2\right)\rho_2^{\frac{t-T_0}{2}} \eta+ \bigg(\frac{|\LL|^2\exp(-(t+T_0+2)d^2)}{1-\exp(-2d^2)}\nonumber\\
       &\hspace{0.5cm}+\frac{|\LL||\MM|\exp(-(t+T_0+2)c^2)}{1-\exp(-2c^2)}\bigg)2\eta\nonumber\\
       &=O\left(|\LL|\cdot\max\left\{|\LL|,|\MM|\right\}te^{-\gamma t}\right), 
\end{align*}
where $\gamma>0$.

\vspace{-0.1in}
\subsection{Tightening the probabilistic bounds}
Up until this point we have placed a mild information assumption on the $\alpha_{ij}$ (this was a consequence of Lemma~\ref{Lemma:concentration_upper}); namely that the expected values of the $\alpha_{ij}$ are strictly bounded away from 1/2 as required in Corollary~\ref{cor:expectationBoundsAlpha}.  In words this means that there is some information contained in the observations $\alpha_{ij}$ such that in expectation these values are closer to 0 for transmissions from legitimate agents and closer to 1 for transmissions from malicious agents. An expected value of $\alpha_{ij}$ equal to 1/2 would be the case of no information.  Thus this is the mildest assumption possible for the $\alpha_{ij}$ and is the key to the decaying probabilities from Lemma~\ref{Lemma:concentration_upper} that underpin the majority of our presented results in this paper. However, if it is possible to obtain more information on the $\alpha_{ij}$, for example knowing a bound on its \emph{variance}, then the probabilistic bounds can be made tighter and the resulting performance guarantees such as bounds on the deviation can be made stronger.  This observation is substantiated with analysis in the current section. 

Lemma \ref{Lemma:concentration_upper} lays the foundation for deriving the convergence of consensus of legitimate agents  presented in Theorem \ref{thm:main1}, the deviation probabilistic bound presented in Theorem \ref{thm:main2} and and the probabilistic upper bound on the convergence rate presented in Theorem \ref{thm:main3}.  Lemma \ref{Lemma:concentration_upper} considers a large family of probability measures for the random variables $\alpha_{ij}(t)$ where only the expectation of the variables are fixed to known values.  While the bounds derived in Lemma \ref{Lemma:concentration_upper} decay exponentially over time, they are somewhat loose when higher moments of the  random variables $\alpha_{ij}(t)$ are known. In the special case when the random variables $\alpha_{ij}(t)$  belong to a family of probability  measures of known variance values, we can use the  concentration inequalities that consider both the first and second moments, instead of the  Chernoff-Hoeffding concentration inequality used in the proof of Lemma \ref{Lemma:concentration_upper}. Prominent examples of such inequalities are the  Berstein inequality, 
Bennet's inequality and the improved Bennet's inequality, providing more refined results than the Chernoff-Hoeffding inequality, are included here for completeness.

\begin{theorem}[Berstein inequality]
	Assume that $x_1,\ldots,x_n$ are independent random variables and $E(x_i)=0$, $E(x_i^2)=\sigma_i^2$ and $|x_i|\leq M$ almost surely. Then,
	\begin{align*}
	\Pr\left(\sum_{i=1}^nx_i\geq b \right)\leq \exp\left(\frac{\frac{1}{2}b^2}{\sum_{i=1}^n\sigma_i^2+\frac{1}{3}M b}\right).
	\end{align*}
\end{theorem}
\begin{theorem}[Bennet's inequality \cite{doi:10.1080/01621459.1962.10482149}]
	Assume that $x_1,...,x_n$ are independent random variables and $E(x_i)=0$, $E(x_i^2)=\sigma_i^2$ and $|x_i|\leq M$ almost surely. Then, for any $0\leq \textcolor{black}{b}<nM$,
	\begin{align*}
	\Pr\left(\sum_{i=1}^nx_i\geq b\right)\leq \exp\left(-\frac{n\sigma^2}{M^2}h\left(
	\frac{b M}{n\sigma^2}\right)\right),
	\end{align*}
	where $h(x)=(1+x)\ln(1+x)-x$ and $n\sigma^2=\sum_{i=1}^n\sigma_i^2$.
\end{theorem}
\begin{theorem}[Improved Bennet's inequality \cite{doi:10.1080/03610926.2017.1367818}]\label{theorem:improved_bennett}\label{thm:improved_bennet_inq}
	Assume that $x_1\ldots,x_n$ are independent random variables and $E(x_i)=0$, $E(x_i^2)=\sigma_i^2$ and $|x_i|\leq M$ almost surely.  Additionally, let 
	\begin{flalign*}
	A = \frac{M^2}{\sigma^2}+\frac{nM}{b}-1 \text{ and }\: B = \frac{nM}{b} -1,
	\end{flalign*}
	and $\Lambda=A-W(Be^A)$, where $W(\cdot)$ is the Lambert $W$ function.
	Let $\sigma^2=\frac{1}{n}\sum_{i=1}^n\sigma_i^2$, then for any $0\leq b<nM$,
	\begin{align*}
	&\Pr\left(\sum_{i=1}^nx_i\geq b\right) \nonumber\\
	&\hspace{1.5cm}\leq\exp\left[-\frac{\Lambda b}{M}+n\ln\left(1+\frac{\sigma^2}{M^2}\left(e^{\Lambda}-1-\Lambda\right)\right)\right].
	\end{align*}
\end{theorem}
As stated in \cite{Jabera2018}, Bennet's inequality yields strictly tighter approximation than Bernstein's inequality. Bennet's inequality, in turn,  is tightened by the improved Bennet's inequality. For this reason, we present in Lemma \ref{Lemma:concentration_upper_variance}  refined concentration inequalities that improve those presented in Lemma \ref{Lemma:concentration_upper} using the improved Bennet's inequality. Recall Assumption~\ref{assumptions}.2 and Assumption~\ref{assumptions}.3, and additionally assume that \[\text{var}(\alpha_{ij}(t)-1/2)=\sigma_d^2\qquad\hbox{for all $t\geq 0$, $i\in \LL$, $j\in\NN_i\cap\LL$},\] and that 
\[\text{var}(\alpha_{ij}(t)-1/2)
=\sigma_c^2\qquad\hbox{for all $t\geq 0$, $i\in \LL,j\in \NN_i\cap\MM$}.\] 

Applying these tighter concentration inequalities to the history of observations $\beta_{ij}(t)=\sum_{k=0}^{t} \left(\alpha_{ij}(k)-1/2\right)$, we can obtain a faster exponential rate result for the decay of misclassification probabilities (misclassifying a legitimate agent as malicious or a malicious agent as legitimate) than that of  Lemma~\ref{Lemma:concentration_upper}, as seen in the following.

\begin{lemma}\label{Lemma:concentration_upper_variance}
Consider the random variables $\beta_{ij}(t)$ defined in Equation~\eqref{eq:betas} for every legitimate node $i\in\LL$ and every neighbor $j\in\mathcal{N}_i$. 
Then, for every $t\geq 0$ and every $i\in \LL$, $j\in\NN_i\cap\LL,$
\begin{align*}
&\Pr\left(\beta_{ij}(t)< 0\right)\leq \max\Bigg\{\mathbbm{1}_{\{d<0\}},\\ &%\hspace{1.4cm}
\exp\left[-d(t+1)\Lambda_d+
(t+1)\ln\left(1+\sigma_d^2\left(e^{\Lambda_d}-1-\Lambda_d\right)\right)\right]\Bigg\},
\end{align*}
where $d = E(\alpha_{ij}(t)-1/2)$
for $i\in \LL$, $j\in\NN_i\cap\LL$, and
\[\Lambda_d = A_d-W(B_de^{A_d}),\]
with $W(\cdot)$ being the Lambert $W$ function, and 
%scalars $A_d$ and $B_d$ given by
\[A_d = \frac{1}{\sigma_d^2}+\frac{1}{d}-1, \quad B_d = \frac{1}{d}-1.\]
Moreover, for every $t\geq 0$ and every 
$i\in \LL,j\in \NN_i\cap\MM$,
\begin{align*}
&\Pr\left(\beta_{ij}(t)\geq 0\right)\leq\max\Bigg\{\mathbbm{1}_{\{c>0\}},\\ &%\hspace{1.4cm}
\exp\left[c(t+1)\Lambda_c
+(t+1)\ln\left(1+\sigma_c^2\left(e^{\Lambda_c}-1-\Lambda_c\right)\right)\right]\Bigg\},\end{align*}
where $c = E(\alpha_{ij}(t)-1/2)$ for $i\in \LL,j\in \NN_i\cap\MM$, and  $\Lambda_c = A_c-W(B_ce^{A_c})$,
with scalars $A_c$ and $B_c$ given by
\[A_c=\frac{1}{\sigma_c^2}-\frac{1}{c}-1, \quad B_c = -\frac{1}{c}-1.\]
\end{lemma}
\vspace{-0.15in}
\begin{IEEEproof}
Let $i\in \LL$, $j\in\NN_i\cap\LL$, 
and note that the given probabilistic bound holds trivially when $d<0$, so assume that $d>0$ and
define 
\[\gamma_{ij}(t) =d+1/2 -\alpha_{ij}(t).\]
Since $d = E(\alpha_{ij}(t)-1/2)$, 
for the random variable $\gamma_{ij}(t)$ we have $E(\gamma_{ij}(t))=0$ and $\text{var}(\gamma_{ij}(t)) = \sigma_d^2$. Since $\alpha_{ij}(t)\in[0,1]$, it follows that $d\le 1/2$ implying that $\gamma_{ij}(t)\le 1$. Also, since $d>0$ and $\alpha_{ij}(t)\in[0,1]$
it follows that $\gamma_{ij}>-1/2$. 
Thus, $|\gamma_{ij}|\le1$. Using the definition of $\beta_{ij}(t)$ in~\eqref{eq:betas}, we have
\begin{align*}
    \Pr\left(\beta_{ij}(t)< 0\right) 
    &= \Pr\left(\sum_{k=0}^{t}(\alpha_{ij}(k)-1/2)< 0\right)\nonumber\\
    &= \Pr\left(\sum_{k=0}^{t}
    (\alpha_{ij}(k)-1/2-d)< -(t+1)d\right)\nonumber\\
    &= \Pr\left(-\sum_{k=0}^{t}\gamma_{ij}(t)< -(t+1)d\right)\nonumber\\
    &= \Pr\left(\sum_{k=0}^{t}\gamma_{ij}(t)>(t+1)d\right).
\end{align*}
We now invoke Theorem~\ref{thm:improved_bennet_inq} for the sequence of random variables $\gamma_{ij}(0),\ldots,\gamma_{ij}(t)$,
with $M=1$, $n=t+1$ and $b=(t+1)d$, which yields the desired relation.
Note that the condition $b <n M$ in Theorem~\ref{thm:improved_bennet_inq} is satisfied in our case here, since $b=(t+1) d$ with $d\le 1/2$ and $nM=t+1$.

Consider now the case $i\in \LL$, $j\in\NN_i\cap\MM$.
If $c>0$, then the result holds trivially.
Thus, assume that $c<0$,
and consider
\[\gamma_{ij}(t) =\alpha_{ij}(t)-1/2-c.\]
Since $c = E(\alpha_{ij}(t)-1/2)$,
we have that $E(\gamma_{ij}(t))=0$ and $\text{var}(\gamma_{ij}(t)) = \sigma_c^2$. Additionally, since $\alpha_{ij}(t)\in[0,1]$, we can see (similar to the preceding case) 
that $|\gamma_{ij}(t)|\leq 1$. 
Hence, 
\begin{align*}
    \Pr\left(\beta_{ij}(t)\ge 0\right) 
    &= \Pr\left(\sum_{k=0}^{t}(\alpha_{ij}(k)-1/2)\ge  0\right)\nonumber\\
    &= \Pr\left(\sum_{k=0}^{t}
    \left(\alpha_{ij}(k)-1/2-c\right)\ge -(t+1)c\right)\nonumber\\
    &= \Pr\left(\sum_{k=0}^{t}\gamma_{ij}(t)\ge -(t+1)c\right).
\end{align*}
Recall that $c<0$, so we can invoke Theorem~\ref{thm:improved_bennet_inq} for the sequence of random variables $\gamma_{ij}(0),\ldots,\gamma_{ij}(t)$, with $\alpha=-(t+1)c>0$, $M=1$, and $n=t+1$.  
Since $c = E(\alpha_{ij}(t)-1/2)$ and $\alpha_{ij}(t)\in[0,1]$, it follows that $-1/2\le c<0$.
Thus, in our case here, the condition $b <n M$ in Theorem~\ref{thm:improved_bennet_inq} reduces to $-(t+1)c< t+1,$
which holds since $-1/2\le c<0$.
\end{IEEEproof}

Therefore, in the case that additional information on the observations $\alpha_{ij}$ is available, in particular that a bound on the variance of the $\alpha_{ij}$ is known, then the results of Lemma~\ref{Lemma:concentration_upper} can be replaced with Lemma~\ref{Lemma:concentration_upper_variance}, resulting in tighter deviation and convergence rate results when applied to Theorems~\ref{thm:main2} and~\ref{thm:main3}. This shows the generality of our framework and the potential to improve our derived performance guarantees as more information on the $\alpha_{ij}$ observations becomes available.
Similarly, 
Lemma~\ref{Lemma:concentration_upper_variance} can be extended to the case where $\text{var}(\beta_{ij}(t))$ varies over time and is link dependent.

\section{Numerical Results}\label{sec:sims}
\begin{figure}
    \centering
    \includegraphics[scale=1,trim={3cm 4cm 2cm 1.5cm},clip]{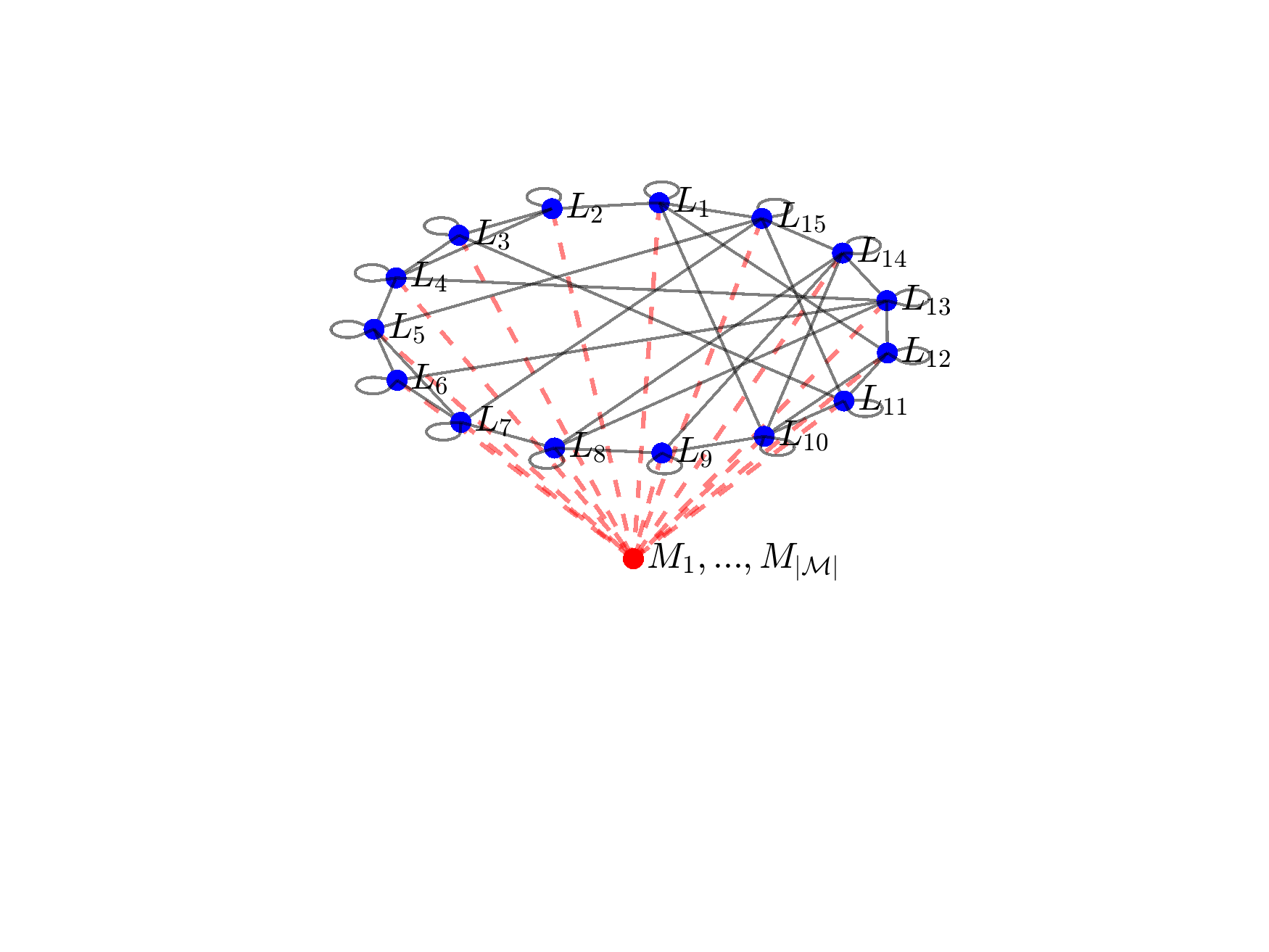}
    \caption{The undirected subgraph  $\G$ of legitimate and malicious agents. Two agents are neighbors if they are connected by an edge. Legitimate agents are depicted by blue nodes and malicious agents are depicted by a red node. Edges between legitimate agents are depicted by black solid lines. Edges between legitimate and malicious agents are depicted by red dashed lines.}
    \label{fig:legitimate_connectivity_graph}
\end{figure}

\begin{figure}
    \centering
    \includegraphics[scale=0.53]{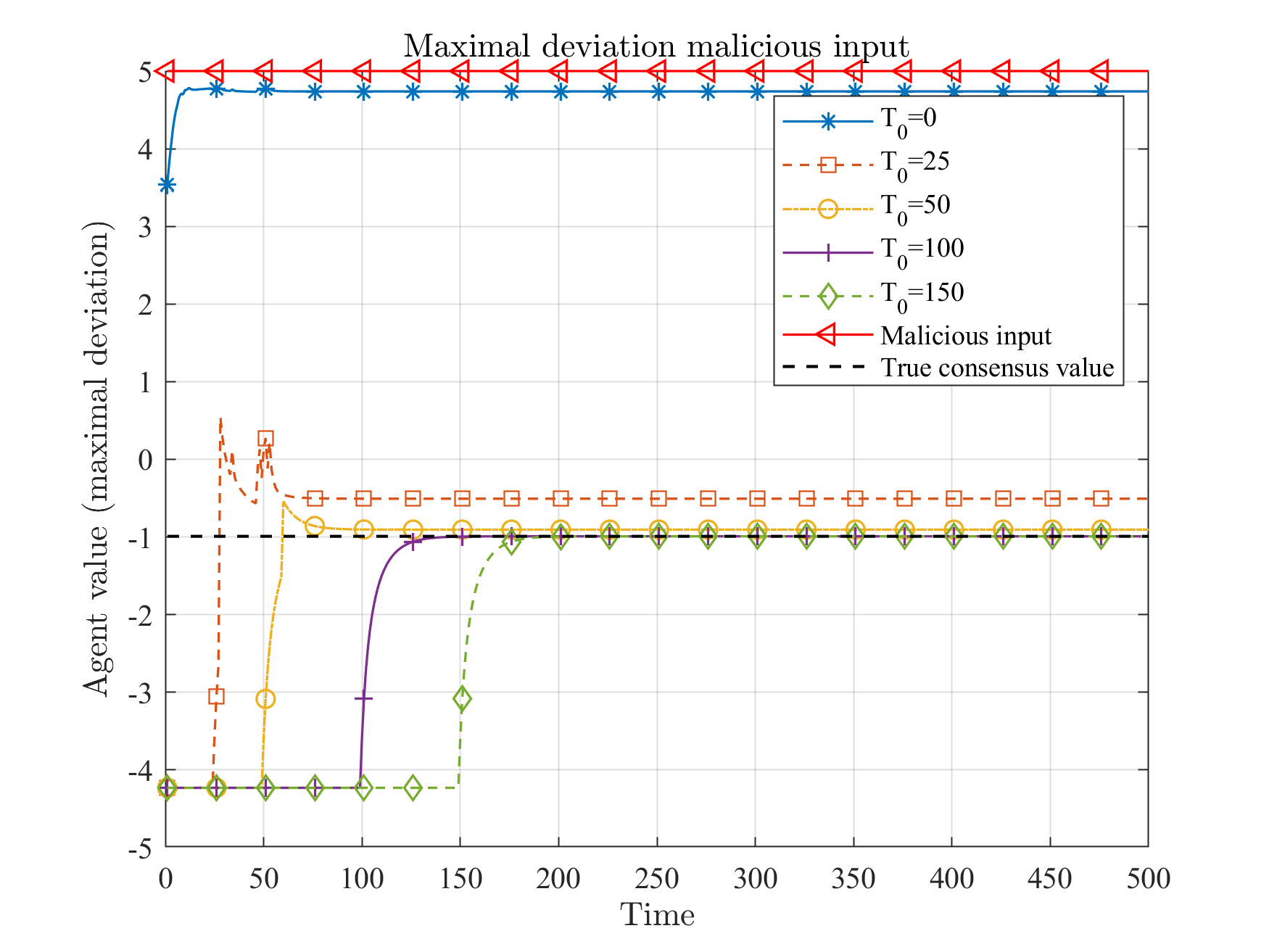}
    \caption{This figure depicts the agent value with the maximal deviation from true consensus value for a single system realization for a system with 15 malicious agents,  $\ell=0.4$ and maximal deviation malicious inputs.    }
    \label{fig:agent_val_realization_max_dev}
    \vspace{-0.2cm}
\end{figure}

\begin{figure}
    \centering
    \includegraphics[scale=0.53]{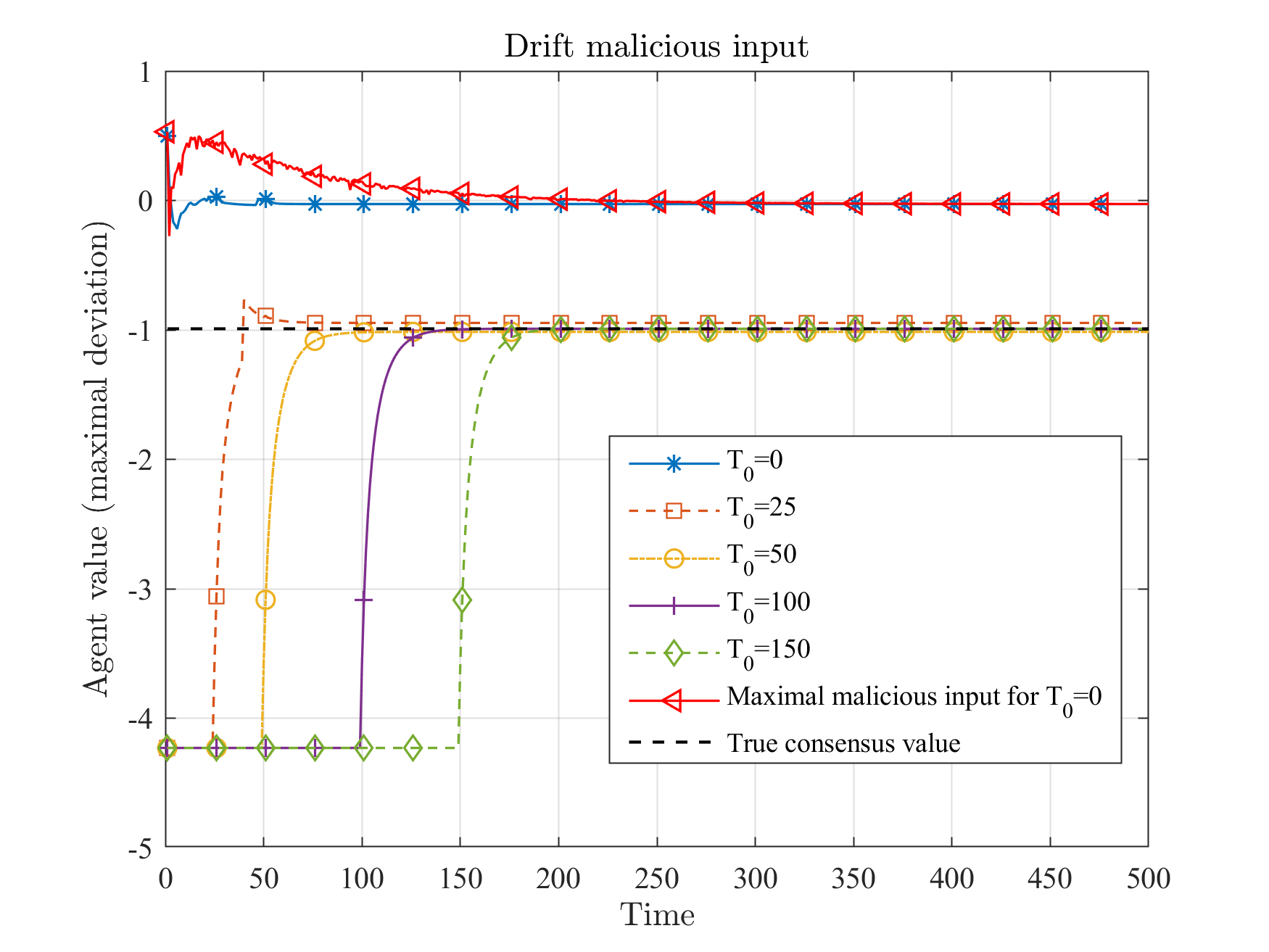}
        \caption{This figure depicts the agent value with the maximal deviation from true consensus value for a single system realization for a system with 15 malicious agents,  $\ell=0.4$ and drift malicious inputs. }
    \label{fig:agent_val_realization_drift}
    \vspace{-0.2in}
\end{figure}

\begin{figure*}
\def\twidth{0.36}
 \subfloat[][(R1,C1):\:\: 5 malicious agents $\ell= 0.2$]{\includegraphics[width=\twidth\textwidth]{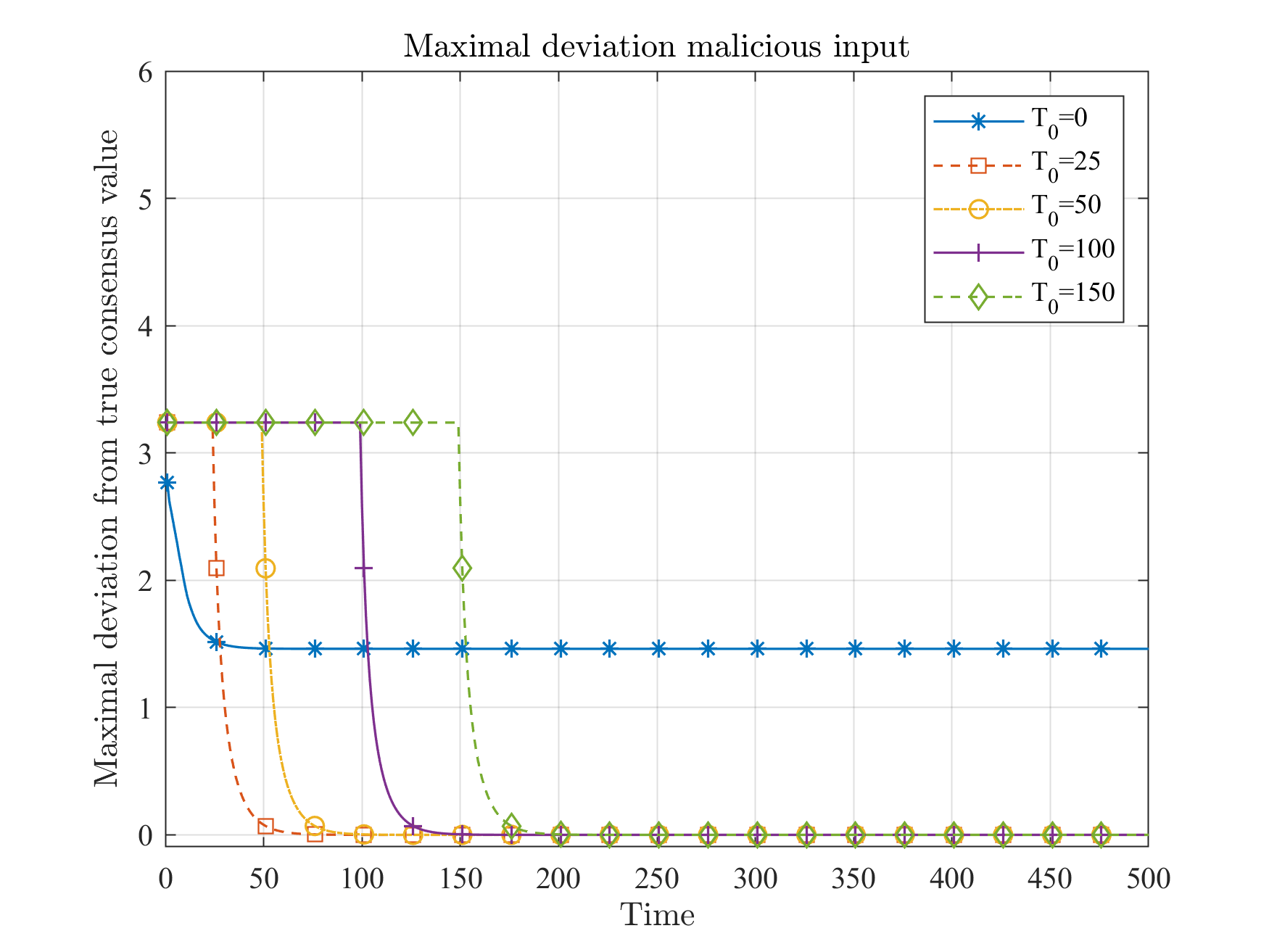}\label{fig:R1_C1}}\hspace{-0.5cm}
 \subfloat[][(R1,C2):\:\: 5 malicious agents $\ell= 0.4$]{\includegraphics[width=\twidth\textwidth]{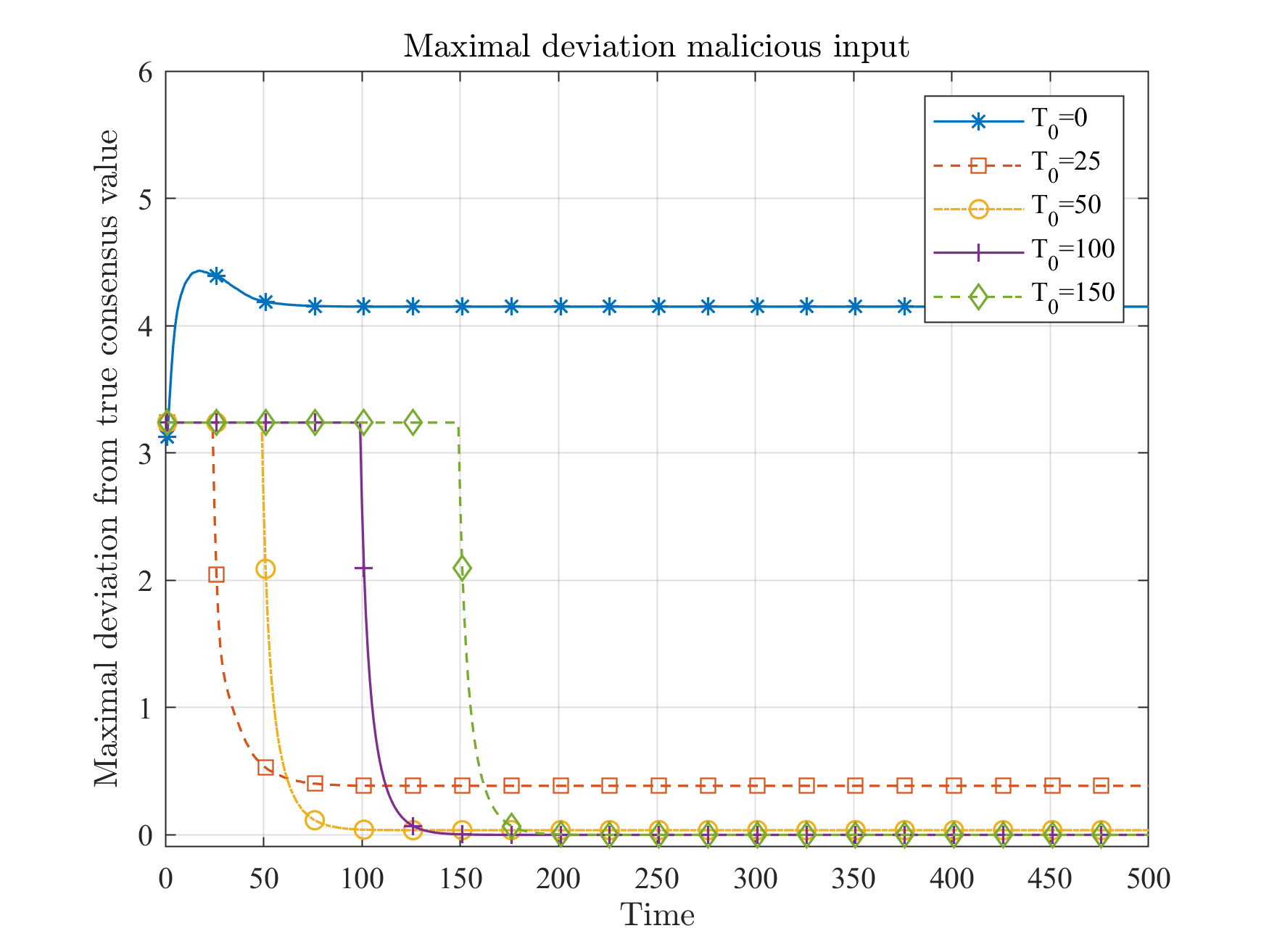}\label{fig:R1_C2}}\hspace{-0.5cm}
 \subfloat[][(R1,C3):\:\: 5 malicious agents $\ell= 0.6$]{\includegraphics[width=\twidth\textwidth]{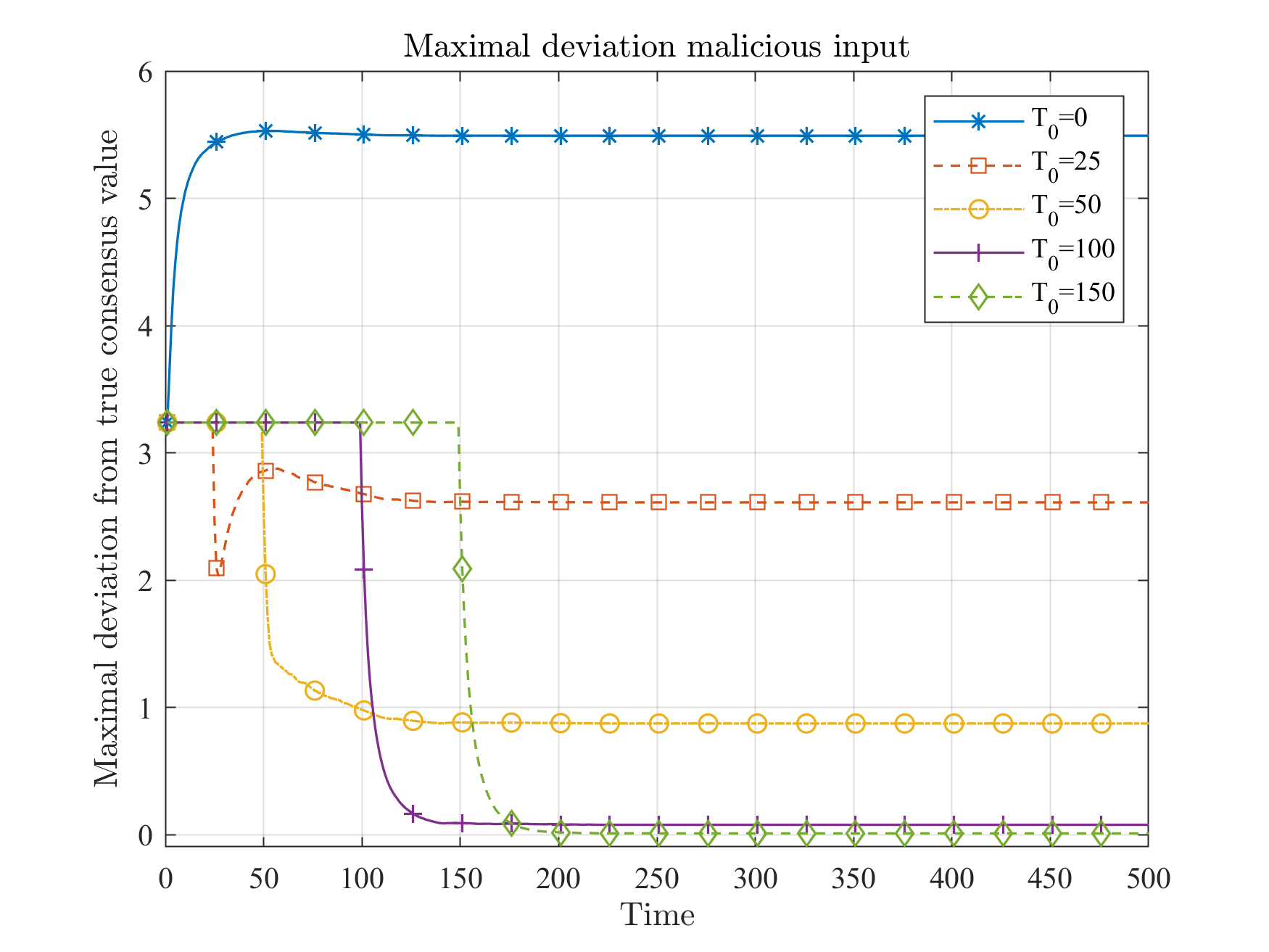}\label{fig:R1_C3}}\\
 \subfloat[][(R2,C1):\:\:15 malicious agents $\ell= 0.2$]{\includegraphics[width=\twidth\textwidth]{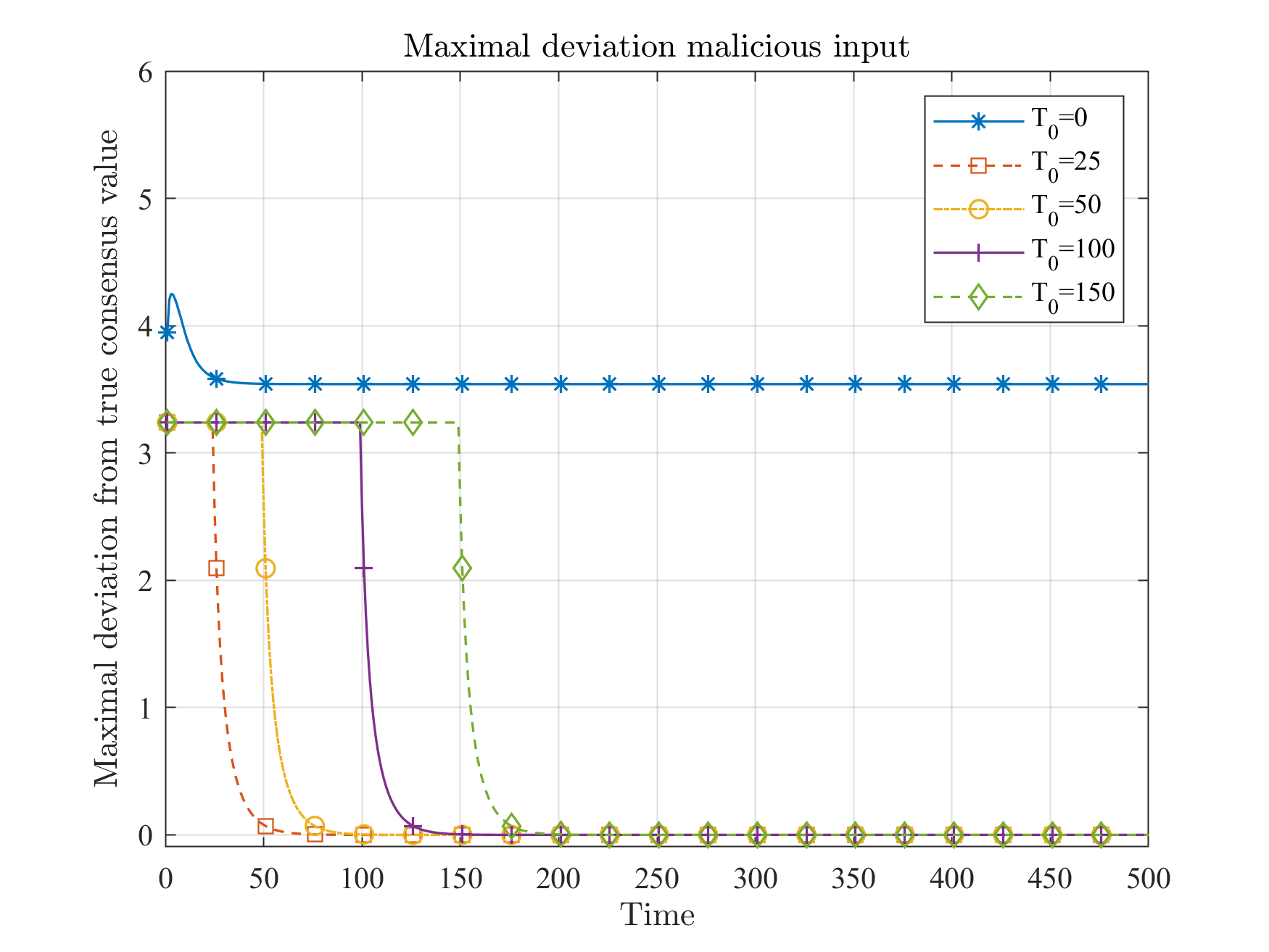}\label{fig:R2_C1}}\hspace{-0.5cm}
 \subfloat[][(R2,C2):\:\:15 malicious agents $\ell= 0.4$]{\includegraphics[width=\twidth\textwidth]{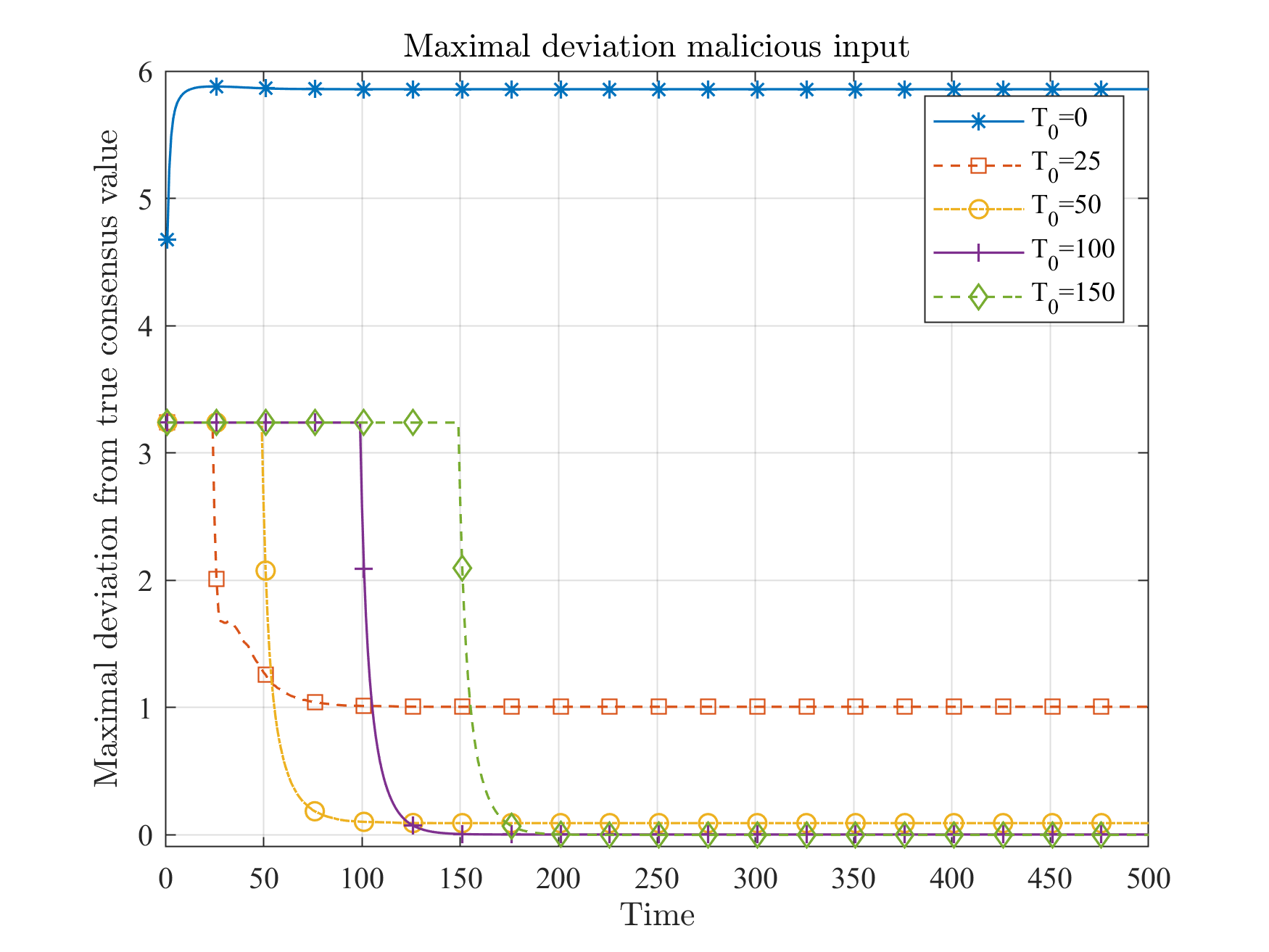}\label{fig:R2_C2}}\hspace{-0.5cm}
 \subfloat[][(R2,C3):\:\:15 malicious agents $\ell= 0.6$]{\includegraphics[width=\twidth\textwidth]{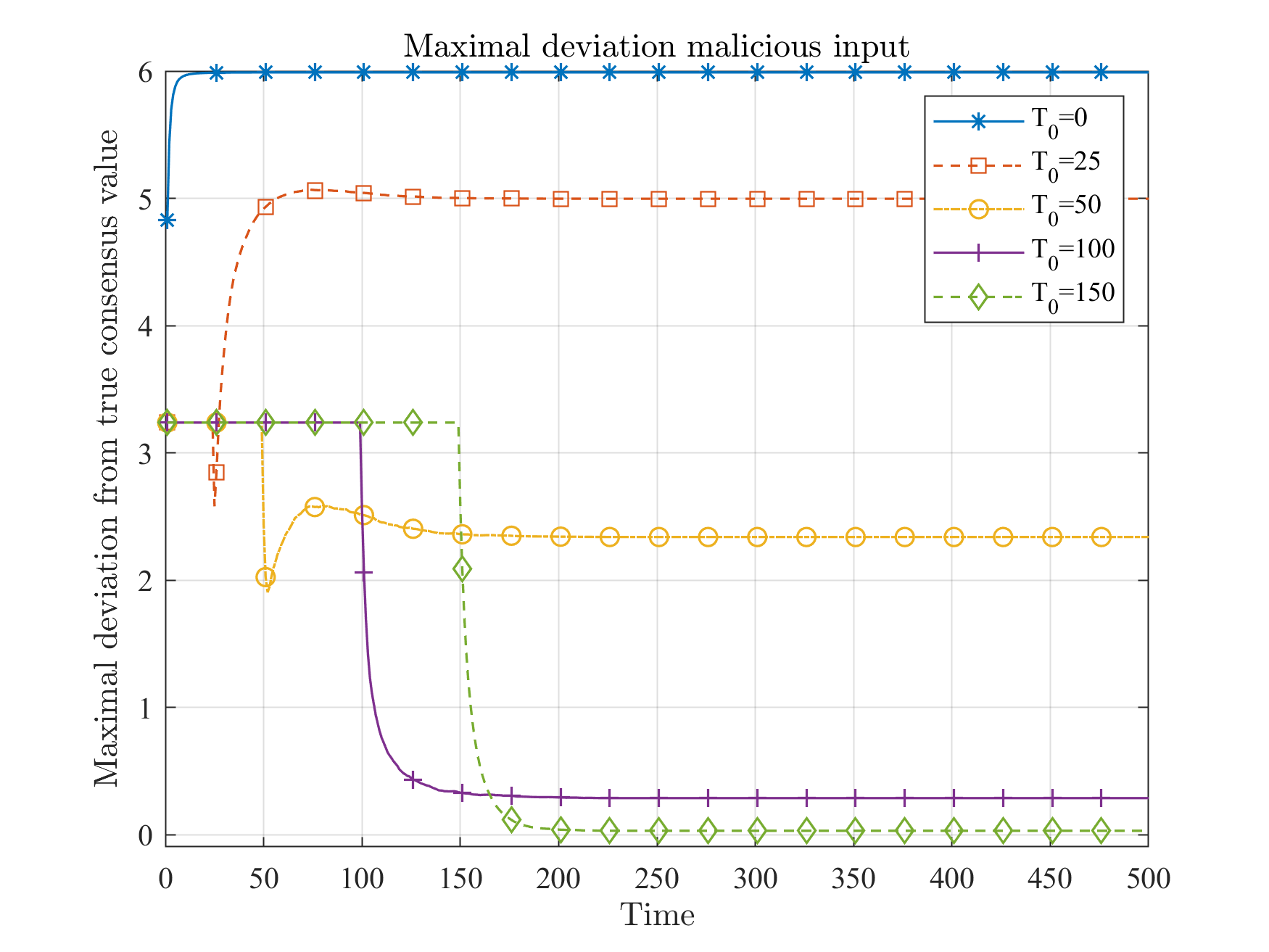}\label{fig:R2_C3}}\\ 
 \subfloat[][(R3,C1):\:\:30 malicious agents $\ell= 0.2$]{\includegraphics[width=\twidth\textwidth]{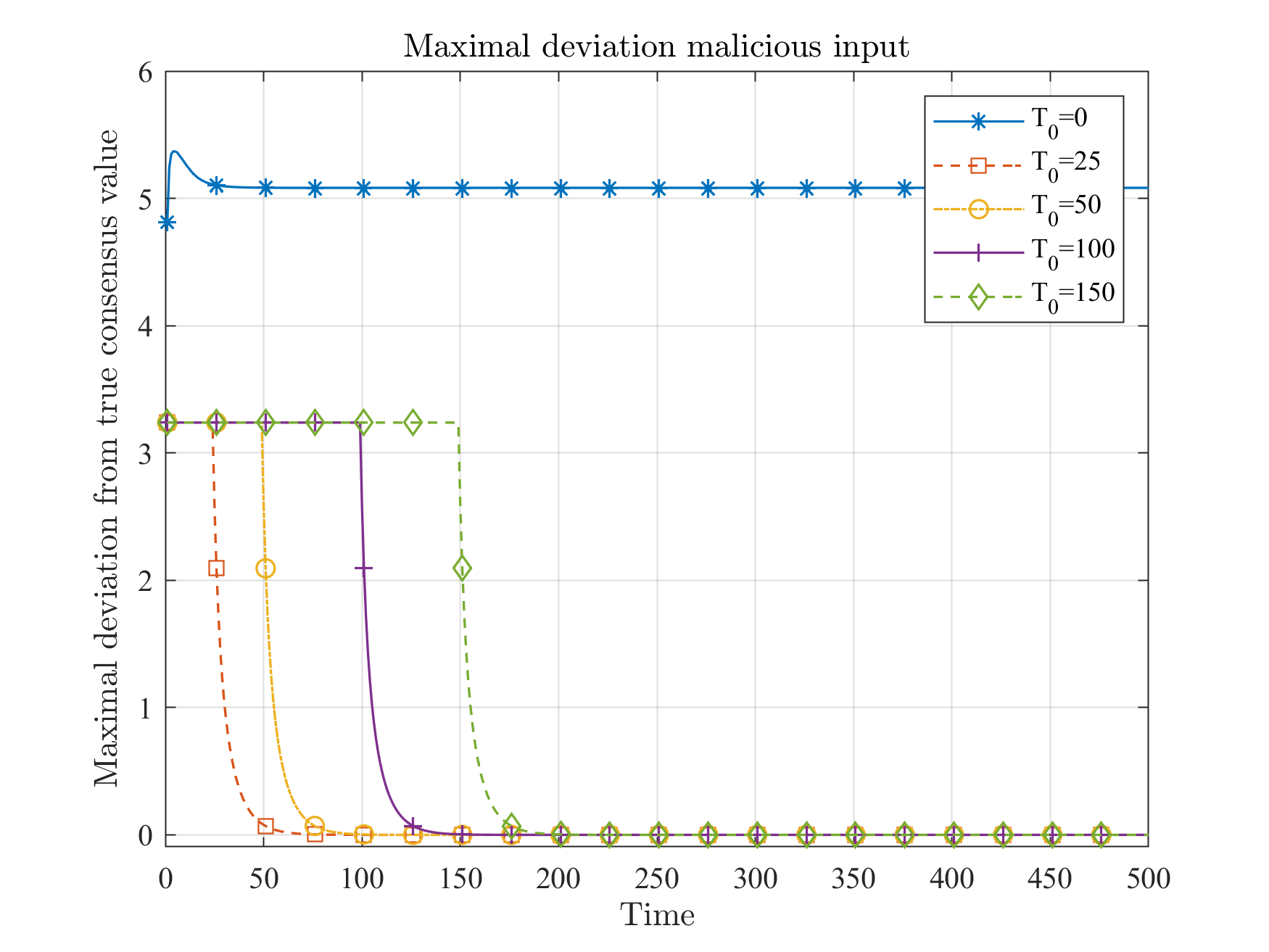}\label{fig:R3_C1}}\hspace{-0.5cm}
 \subfloat[][(R3,C2):\:\:30 malicious agents $\ell= 0.4$]{\includegraphics[width=\twidth\textwidth]{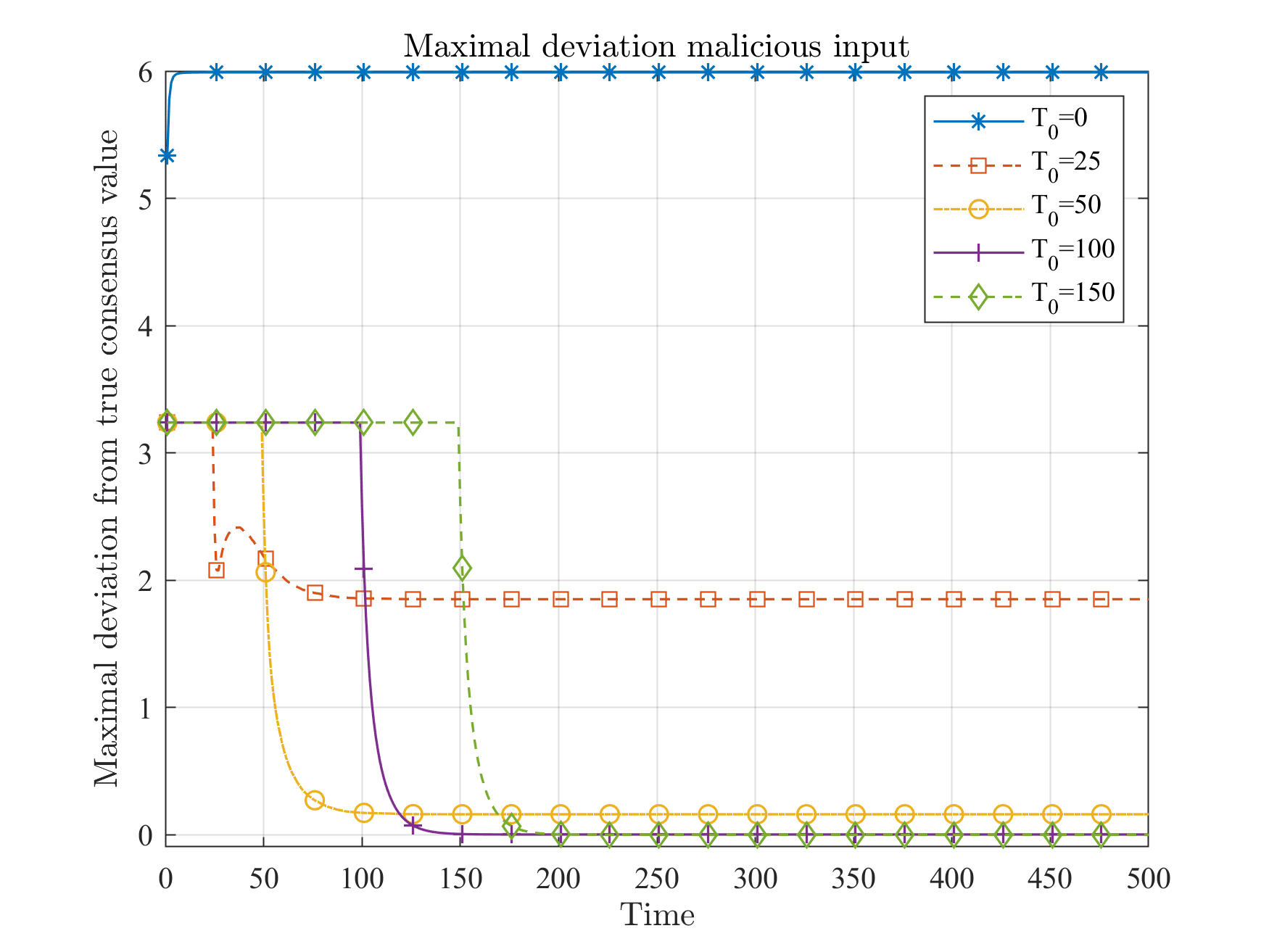}\label{fig:R3_C2}}\hspace{-0.5cm}
 \subfloat[][(R3,C3):\:\:30 malicious agents $\ell= 0.6$]{\includegraphics[width=\twidth\textwidth]{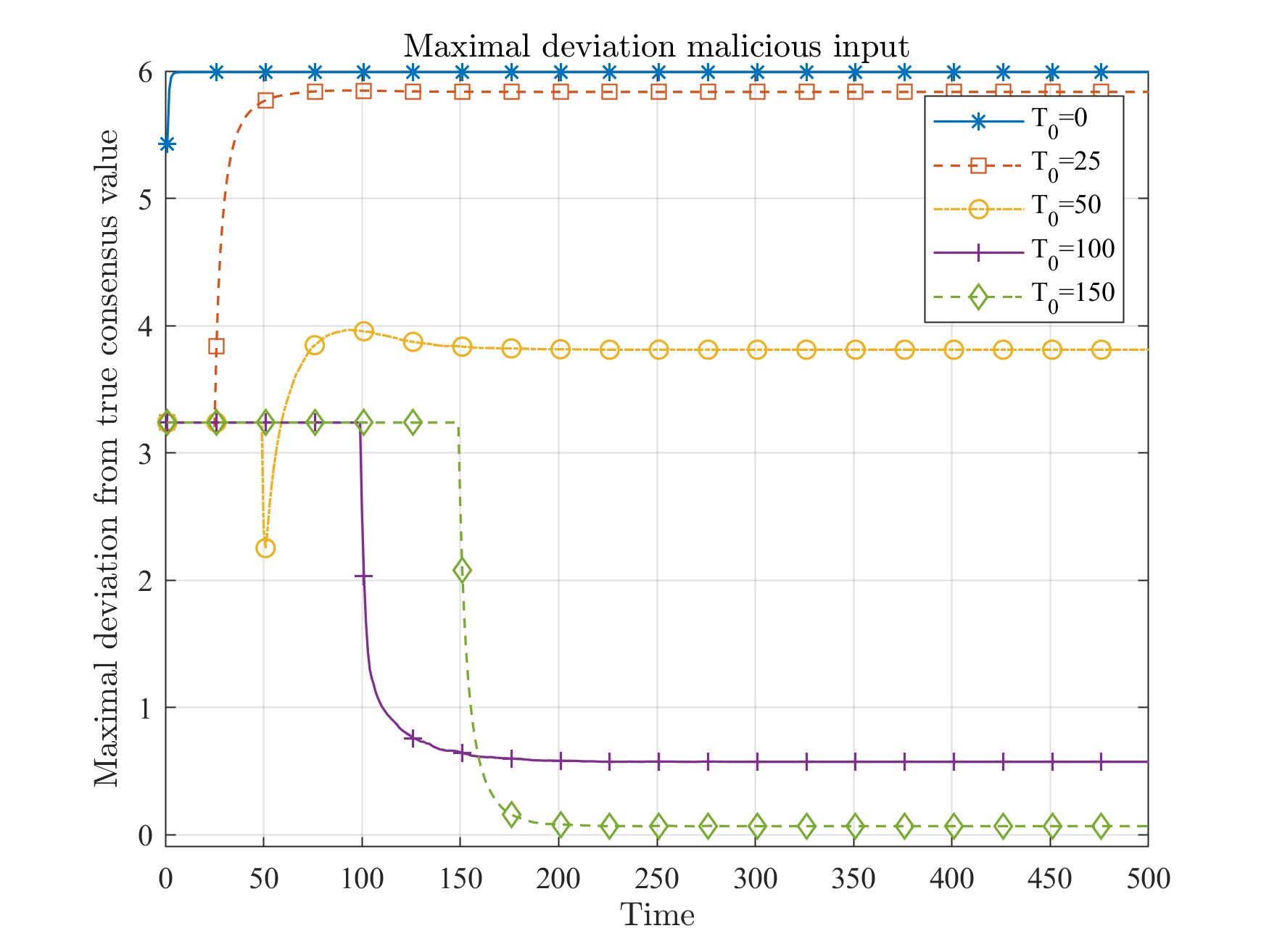}\label{fig:R3_C3}}
 \caption{Maximal deviation input}\label{max_deviation_fig}
  \vspace{-0.1in}
\end{figure*}

In this section we evaluate the performance of the proposed scheme in this work using Monte Carlo simulations. We aim at investigating the effects, on several performance metrics, of the following system parameters: the number $|\MM|$ of malicious agents, the starting time $T_0$ of the data passing stage and the variance of the random variables $\alpha_{ij}$. Performance metrics that we study include deviation from the true consensus value and the convergence rate. We evaluate a system setup with  $15$ legitimate agents, and with $\kappa=10$. We consider the following values for the starting time $T_0$: 0, 25, 50, 100, 150.
The adjacency matrix of the subgraph $\G_{\LL}$ of links among legitimate agents is 
\setcounter{MaxMatrixCols}{20}
\begin{align}\label{adj_mat_num_rslt}
\small{
    \begin{pmatrix}
1 & 1 & 0 & 0 & 0 & 0 & 0 & 0 & 0 & 1 & 0 & 1 & 0 & 0 & 1 \\
1 & 1 & 1 & 1 & 0 & 0 & 0 & 0 & 0 & 0 & 0 & 0 & 0 & 0 & 0 \\
0 & 1 & 1 & 1 & 0 & 0 & 0 & 0 & 0 & 0 & 1 & 0 & 0 & 0 & 0 \\
0 & 1 & 1 & 1 & 1 & 0 & 0 & 0 & 0 & 0 & 0 & 0 & 1 & 0 & 0 \\
0 & 0 & 0 & 1 & 1 & 1 & 1 & 0 & 0 & 0 & 0 & 0 & 0 & 0 & 1 \\
0 & 0 & 0 & 0 & 1 & 1 & 1 & 0 & 0 & 0 & 0 & 0 & 1 & 0 & 0 \\
0 & 0 & 0 & 0 & 1 & 1 & 1 & 1 & 0 & 0 & 0 & 0 & 0 & 0 & 1 \\
0 & 0 & 0 & 0 & 0 & 0 & 1 & 1 & 1 & 0 & 0 & 0 & 1 & 1 & 0 \\
0 & 0 & 0 & 0 & 0 & 0 & 0 & 1 & 1 & 1 & 0 & 0 & 0 & 1 & 0 \\
1 & 0 & 0 & 0 & 0 & 0 & 0 & 0 & 1 & 1 & 1 & 1 & 0 & 1 & 0 \\
0 & 0 & 1 & 0 & 0 & 0 & 0 & 0 & 0 & 1 & 1 & 1 & 0 & 0 & 1 \\
1 & 0 & 0 & 0 & 0 & 0 & 0 & 0 & 0 & 1 & 1 & 1 & 1 & 0 & 0 \\
0 & 0 & 0 & 1 & 0 & 1 & 0 & 1 & 0 & 0 & 0 & 1 & 1 & 1 & 0 \\
0 & 0 & 0 & 0 & 0 & 0 & 0 & 1 & 1 & 1 & 0 & 0 & 1 & 1 & 1 \\
1 & 0 & 0 & 0 & 1 & 0 & 1 & 0 & 0 & 0 & 1 & 0 & 0 & 1 & 1
\end{pmatrix}},
\end{align}
and we assume that every malicious agent is connected to all the legitimate agents in the system which is the worst-case scenario.
We depict the topology corresponding to the connectivity graph $\G$ in Figure \ref{fig:legitimate_connectivity_graph}; since every malicious agent is connected to all legitimate agents we depict the malicious agents in Figure \ref{fig:legitimate_connectivity_graph} by a single node.
The vector of initial values of legitimate agents is
\begin{align*}
x_{\LL}(0)=&(-2.59, -2.44, -4.23, -1.45 , -1.46 ,0.871,-0.51,\nonumber\\
&  -3.19  ,-0.59  ,-3.31  ,-2.25  ,1.31  ,1.87,  1.34 
,1.76 )'
\end{align*}
where $(\cdot)'$ denotes the transpose operator, and we use $\eta=5$ as the upper bound for the agents' initial values.
Additionally, for every legitimate agent $i$, we have $E(\alpha_{ij}) = 0.55$ if $j\in\mathcal{N}_i\cap \LL$, 
and $E(\alpha_{ij}) = 0.45$ if $j\in\mathcal{N}_i\cap \MM$. The random variable $\alpha_{ij}$ is uniformly distributed on the interval $\left[E(\alpha_{ij})-\frac{\ell}{2},E(\alpha_{ij})+\frac{\ell}{2}\right]$. We consider the following three values for $\ell$: $0.2,0.4,0.6$, for which the standard deviation values are $0.0577, 0.1155, 0.1732$, respectively. Note that the larger the $\ell$, the larger the variance of $\alpha_{ij}$ about its mean value. 
%For the sake of evaluating the worst effect possible by the inputs of malicious agents, we assume that every malicious agent is connected to all the legitimate agents. 
We consider the following input scenarios presented in Figures \ref{fig:agent_val_realization_max_dev} and \ref{fig:agent_val_realization_drift}. Both attacks assume that the malicious agents are allowed to communicate and coordinate attacks.
\paragraph{Maximal deviation malicious input}
To measure the maximal possible deviation from true consensus value we consider the case were the malicious agents inputs are the furthest from the true consensus value. That is, if the true consensus value is positive than the malicious input is $-\eta$, and if the true consensus value is negative than the malicious input is $\eta$.
 In our scenario, all the malicious agents choose the maximum possible input $\eta=5$.  Figures \ref{fig:agent_val_realization_max_dev} and \ref{max_deviation_fig}  measure the effect of the malicious agent inputs that lead to the maximum deviation from the true consensus value. Figure~\ref{fig:agent_val_realization_max_dev} depicts agent values over time when running consensus according to the protocol from Equation~\eqref{eq:consensusProtocol} and the weights defined by Equation~\eqref{eq:weights}.  Specifically, 
 Figure \ref{fig:agent_val_realization_max_dev} depicts the  value of the agent with the maximal deviation from the true consensus value for a single system realization with maximum deviation malicious input. Additionally, Figure \ref{max_deviation_fig} depicts the  absolute value of the maximum deviation from the true consensus value averaged over 500 system realizations. Figures \ref{fig:agent_val_realization_max_dev} and \ref{max_deviation_fig} additionally show the differences in converged values and deviation respectively for different consensus start times $T_0$.
 
 The maximum deviation malicious input leads to an upper bound to all possible malicious attacks, and depicts a worst case scenario in terms of deviation from true consensus value.
  However, though this attack is useful to analyze the worst case effect of malicious inputs it is easy to detect using outlier rejection methods for example~\cite{wmsr}. For this reason, we consider the following additional attack scenario that is harder to detect.

\begin{figure*}
\def\twidth{0.36}
 \subfloat[(R1,C1):\:\:5 malicious agents $\ell= 0.2$]{\includegraphics[width=\twidth\textwidth]{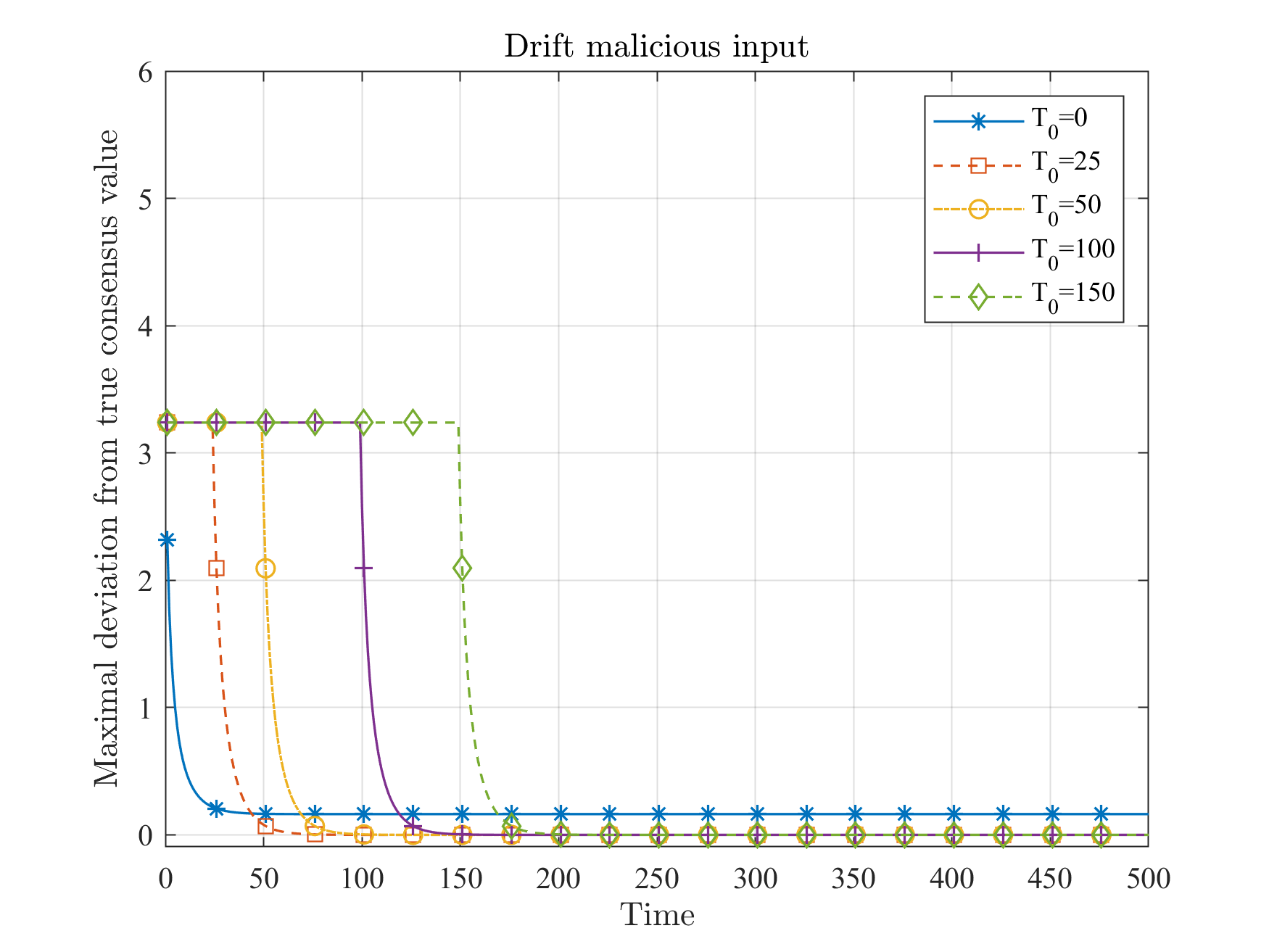}\label{fig:R1_C1_drift}}\hspace{-0.5cm} 
 \subfloat[(R1,C2):\:\:5 malicious agents $\ell= 0.4$]{\includegraphics[width=\twidth\textwidth]{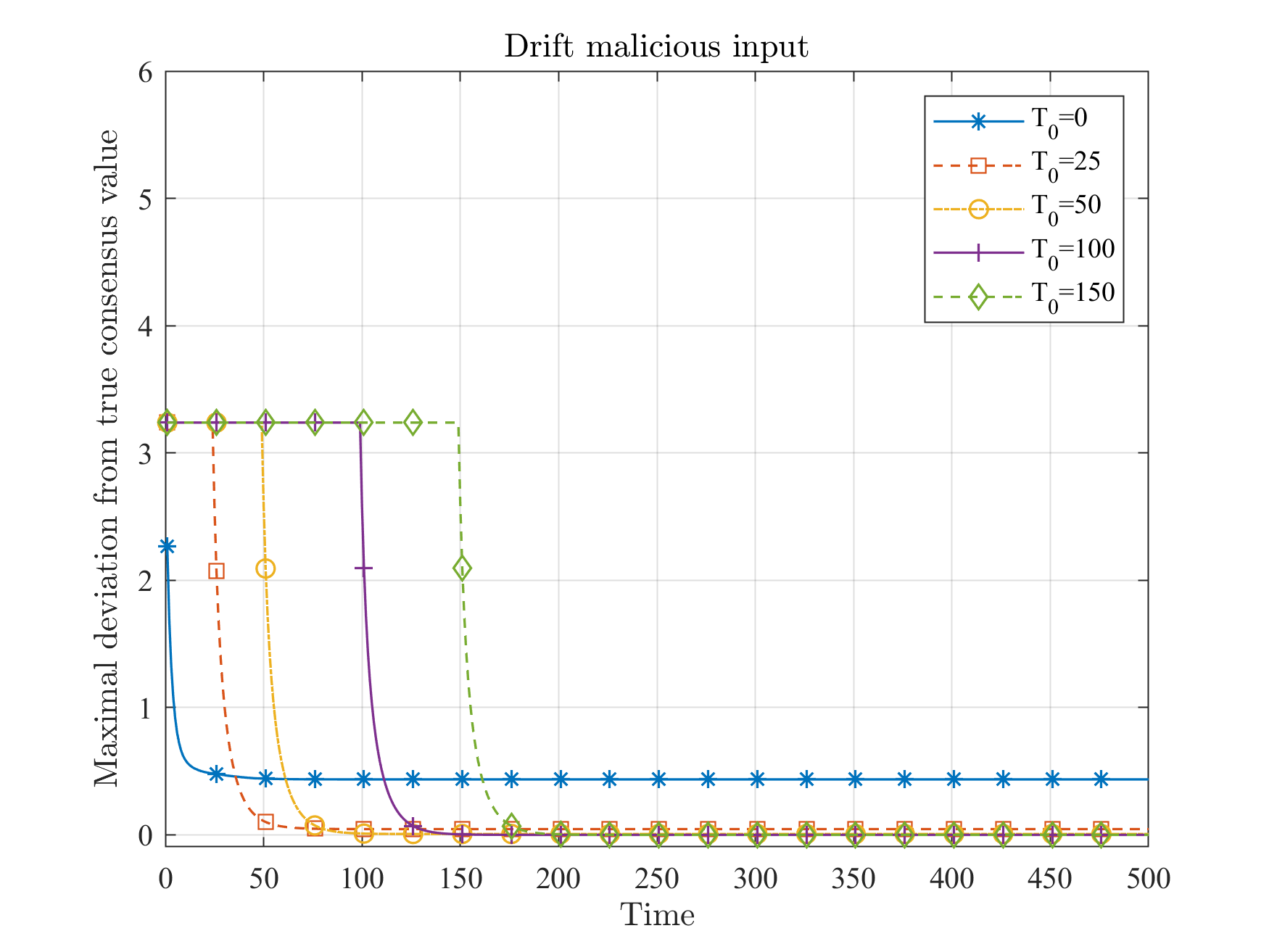}\label{fig:R1_C2_drift}}\hspace{-0.5cm}
 \subfloat[(R1,C3):\:\:5 malicious agents $\ell= 0.6$]{\includegraphics[width=\twidth\textwidth]{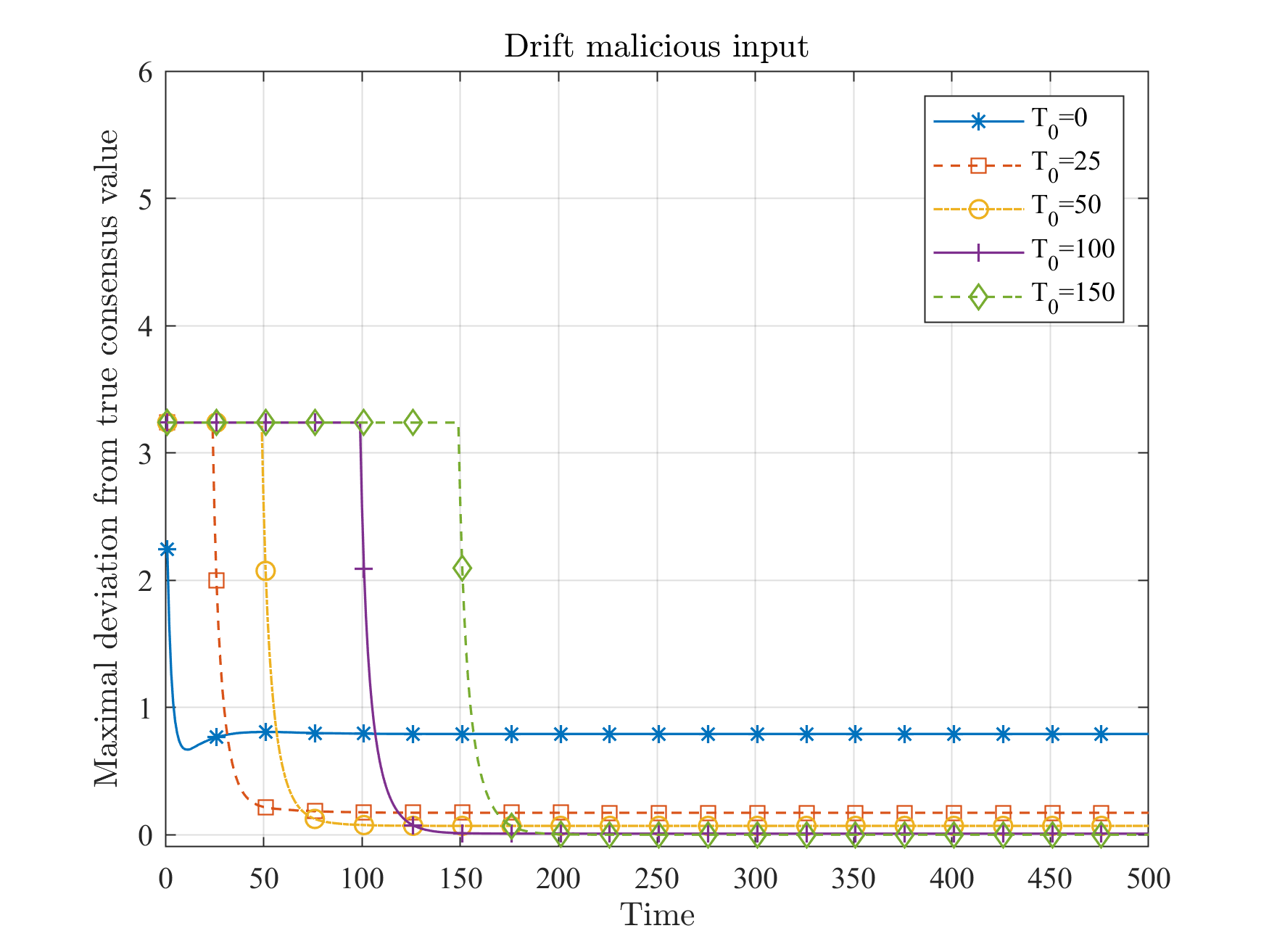}\label{fig:R1_C3_drift}}\\
 \subfloat[(R2,C1):\:\:15 malicious agents $\ell= 0.2$]{\includegraphics[width=\twidth\textwidth]{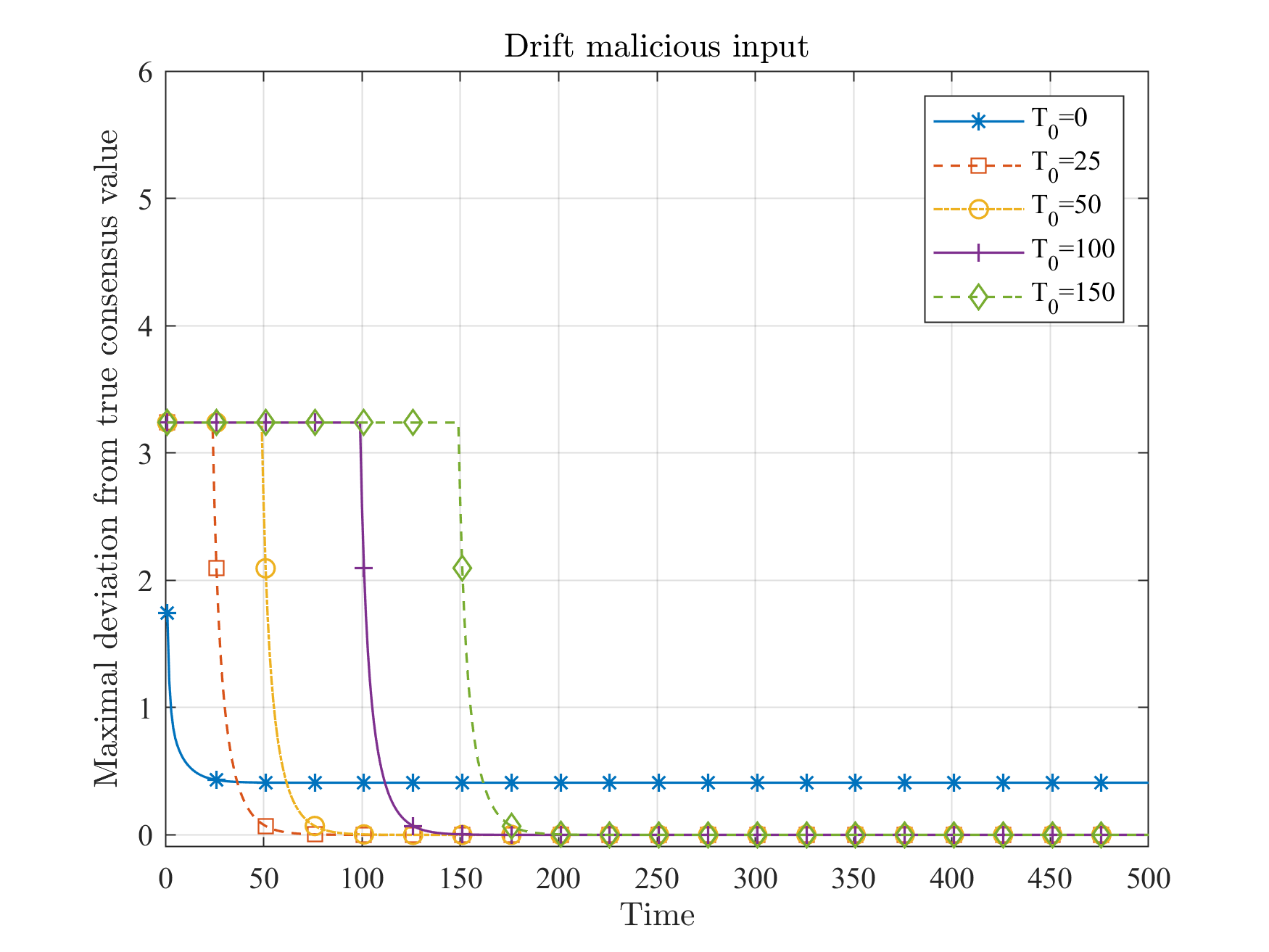}\label{fig:R2_C1_drift}}\hspace{-0.5cm}
 \subfloat[(R2,C2):\:\:15 malicious agents $\ell= 0.4$]{\includegraphics[width=\twidth\textwidth]{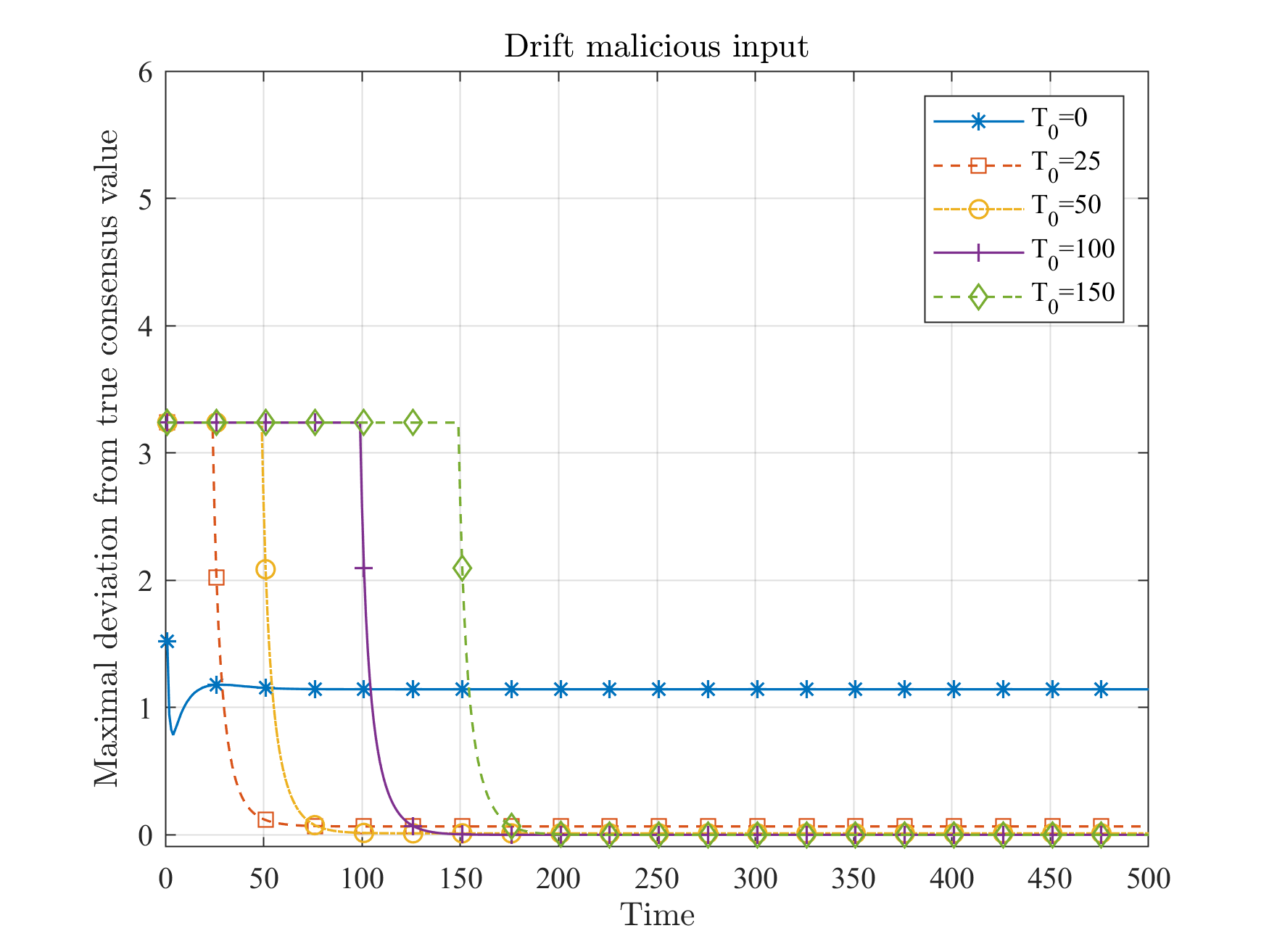}\label{fig:R2_C2_drift}}\hspace{-0.5cm}
 \subfloat[(R2,C3):\:\:15 malicious agents $\ell= 0.6$]{\includegraphics[width=\twidth\textwidth]{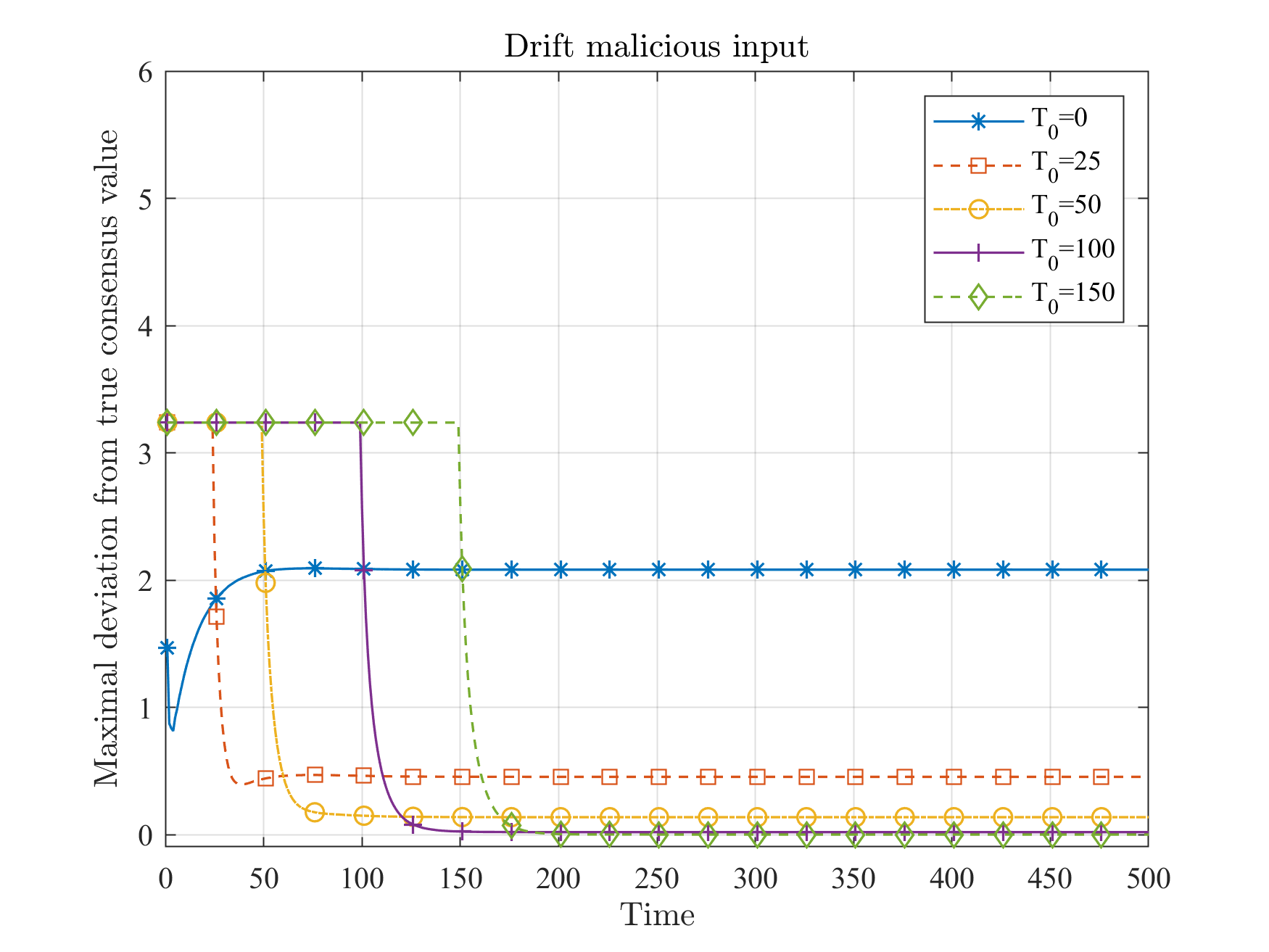}\label{fig:R2_C3_drift}}\\
 \subfloat[(R3,C1):\:\:30 malicious agents $\ell= 0.2$]{\includegraphics[width=\twidth\textwidth]{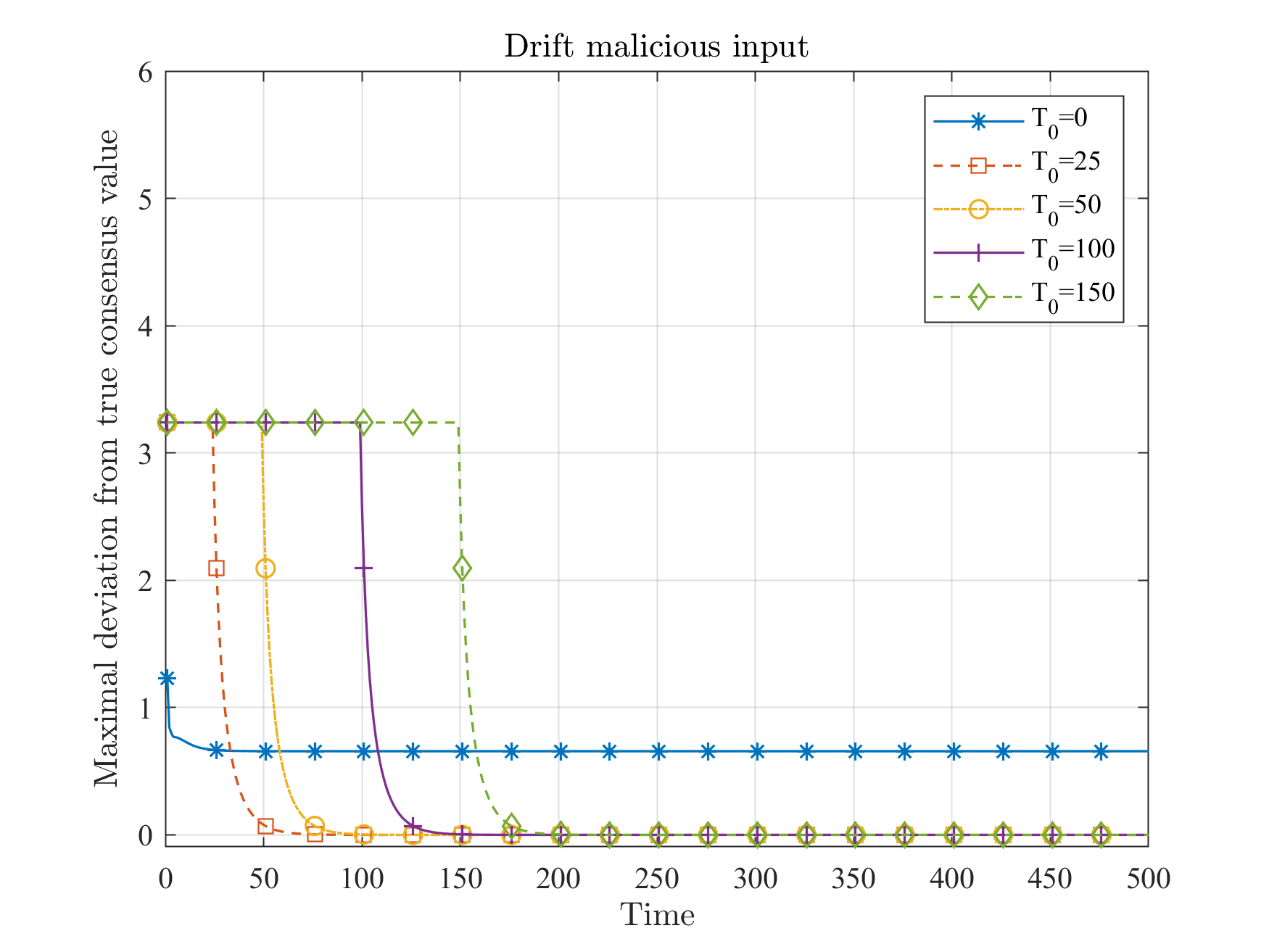}\label{fig:R3_C1_drift}}\hspace{-0.5cm}
 \subfloat[(R3,C2):\:\:30 malicious agents $\ell= 0.4$]{\includegraphics[width=\twidth\textwidth]{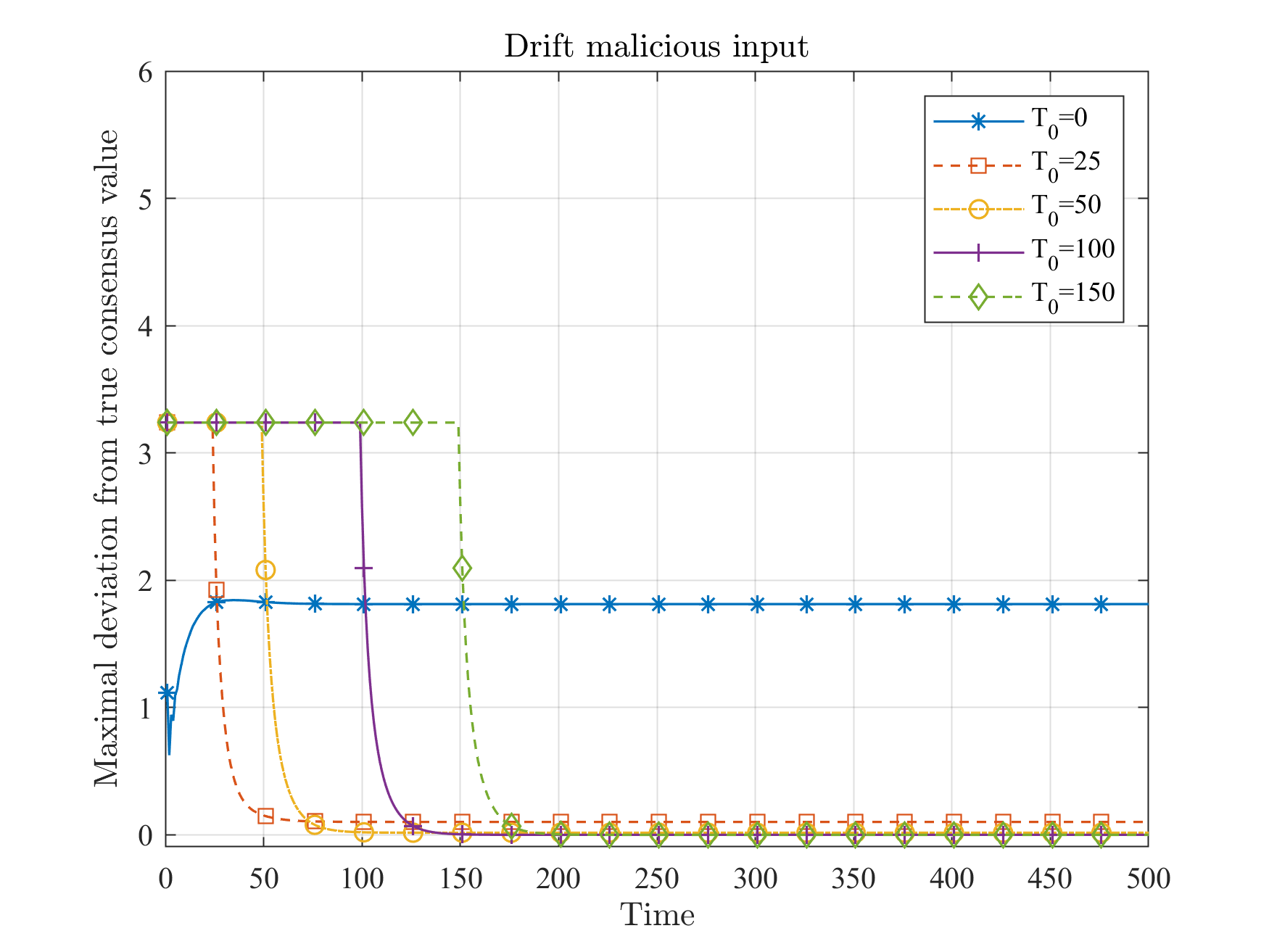}\label{fig:R3_C2_drift}}\hspace{-0.5cm}
 \subfloat[(R3,C3):\:\:30 malicious agents $\ell= 0.6$]{\includegraphics[width=\twidth\textwidth]{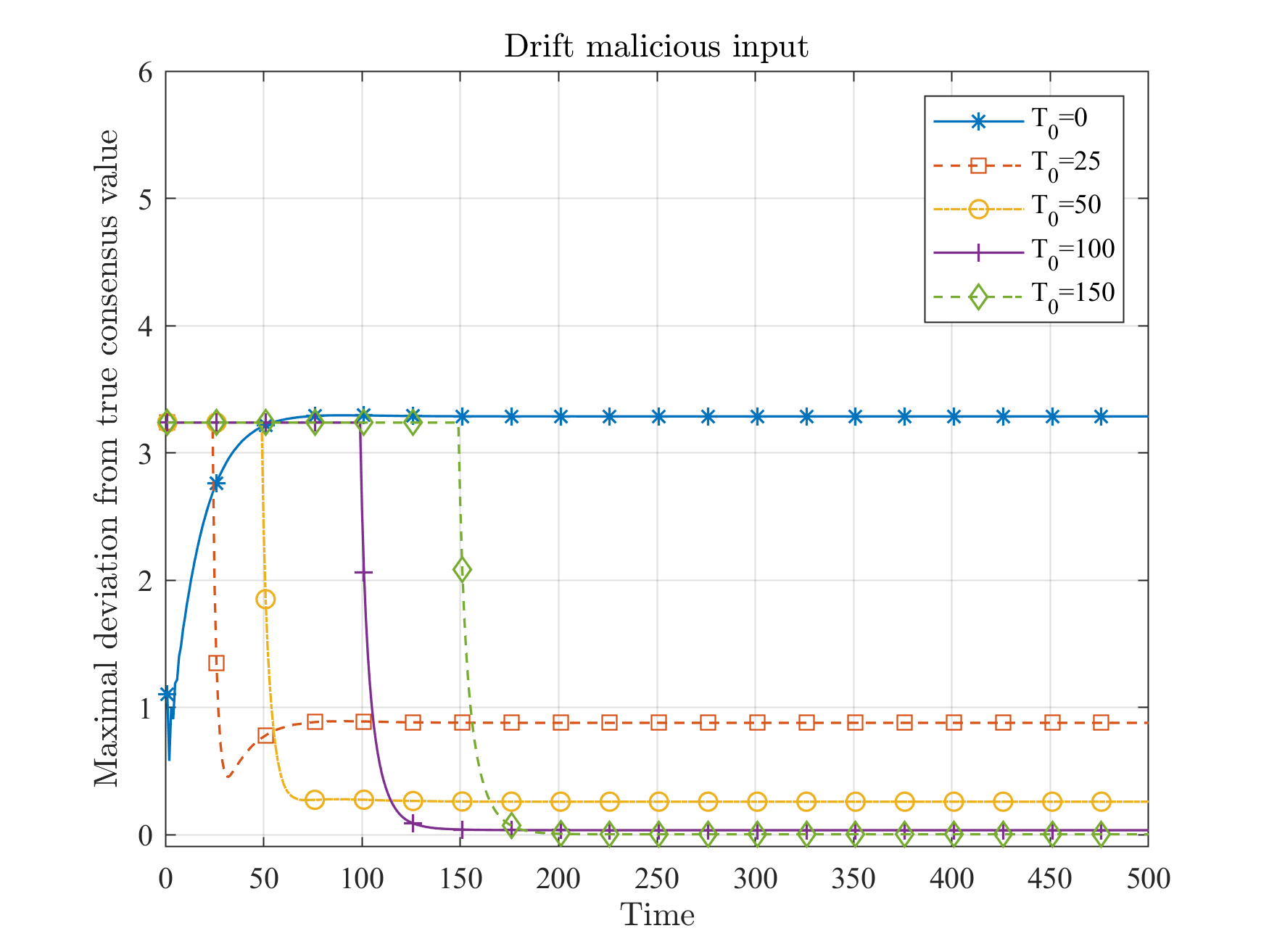}\label{fig:R3_C3_drift}}
 \caption{Drift input}\label{drift_fig}
  \vspace{-0.1in}
\end{figure*}

\paragraph{Drift malicious input}
Figure \ref{fig:agent_val_realization_drift} and \ref{drift_fig}  measure the effect of a less transparent  malicious agents' attack where every malicious agent adds a drift term to the average legitimate agent values. More specifically, let $d_{\LL}(m,t)$ be the average value at time $t$ of the legitimate agents that are neighbors of the malicious agent $m$. Denote $d_{\LL}(t) = \frac{1}{|\mathcal{M}|}\sum_{m=1}^{|\mathcal{M}|}d_{\LL}(m,t)$. Then, at time $t$, the malicious agent $m$ chooses its state value $x_m(t)$ according to the following rule: 
\[x_m(t) =  d_{\LL}(t-1)+0.15\times d_{\MM}(m,t),\]
where $x_m(t=T_0-1)$ is chosen randomly from the set $[-0.15\times\eta,0.15\times\eta]$ and for every $t\geq T_0$
\[d_{\MM}(m,t) = -\text{sign}(d_{\LL}(1))\times \eta\times(0.75)^{0.05(t-T_0)}u(m,t),\]
 $u(m,t)\sim U[0,1]$. Thus malicious agent inputs are \emph{time varying} and strategically close to legitimate agent values to avoid easy detection in this case. Additionally, if $x_m(t)> \eta$ then $x_m(t)$ is generated randomly and uniformly from the interval $[\eta-0.05,\eta]$. Similarly if $x_m(t)< -\eta$ then $x_m(t)$ is generated randomly and uniformly from the interval $[-\eta,-\eta+0.05]$. 
 We note that since in our scenario every malicious agents is a neighbor of all the legitimate agents,  $d_{\LL}(m_1,t)=d_{\LL}(m_2,t)$ for all $m_1,m_2\in\mathcal{M}$.
 Figure \ref{fig:agent_val_realization_drift} depicts the  value of the agent with the maximal deviation from the true consensus value for a single system realization with drift malicious input. 
 Additionally, Figure \ref{drift_fig} depicts the resulting absolute value of the maximum deviation from the true consensus value for a drift malicious input averaged over 500 system realizations.

%In addition to Figures \ref{fig:agent_val_realization_max_dev} and \ref{fig:agent_val_realization_drift} that show an average deviation from true consensus value, we also present in 
Figures \ref{fig:agent_val_realization_max_dev} and \ref{fig:agent_val_realization_drift} depict a single realization of the agent value with the maximal deviation from true consensus value for a system with 15 malicious agents and $\ell=0.4$.  Figure \ref{fig:agent_val_realization_max_dev} depicts the agent value with the maximal deviation from true consensus value for the maximal deviation malicious inputs where Figure \ref{fig:agent_val_realization_drift} depicts the agent value with the maximal deviation from true consensus value for the drift malicious inputs.

 In addition to Figures \ref{fig:agent_val_realization_max_dev} and \ref{fig:agent_val_realization_drift} that depict  the legitimate agent values for a \textit{single system realization} we  depict in Figures \ref{max_deviation_fig} and \ref{drift_fig} the \textit{average deviation} of the legitimate agent values from true consensus value. The average deviation is calculated using 500 system realizations. 
 Figures \ref{max_deviation_fig} and \ref{drift_fig} show that the deviation from true consensus value grows both as $\ell$ grows and as the number $|\MM|$ of malicious agents grows. In each figure we label the plot by (R\#,C\#) where R\# denotes the row number and represents the increase in the number of malicious agents, and C\# denotes the column number and represents the growth of $\ell$. Figures \ref{max_deviation_fig} and \ref{drift_fig} show the significance of introducing the variable $T_0$ that sets the starting time of the data passing between agents. The figures indicate  that increasing $T_0$ can significantly reduce the deviation from true consensus value.  Figures \ref{max_deviation_fig} and \ref{drift_fig} also show that our scheme can combat an attack of a large number of malicious agents, even when $|\MM|=30$, that is, there is twice the number of malicious agents than legitimate agents in the system. \emph{This is a significant improvement over classical results} that depend on the network connectivity \cite{origByz,bulloUnreliable,sundaram}. Using the upper bound $\min_{i\in\LL\cup\MM}|\mathcal{N}_i|$ on the network connectivity we can conclude that classical results will fail to detect all the malicious agents when their number is larger than the number of legitimate agents; in this case the undetected malicious agents can control the data values of the legitimate agents. More specifically, we can see from the choice of adjacency matrix \eqref{adj_mat_num_rslt} that for our setup, the network connectivity is at most $3+|\MM|$ and thus  classical results can tolerate less than $\frac{3+|\MM|}{2}$ malicious agents. Since there are $|\MM|$ malicious agents, the consensus algorithm can fail and be controlled by the malicious agents. Furthermore, for the connectivity graph that we consider in the numerical results, the maximal number of malicious agents leading to a graph connectivity guaranteeing their detection and preventing them from controlling the network is $2$, since we must fulfill the condition that $|\MM|<\frac{3+|\MM|}{2}$. Finally, our numerical results show that there exists a finite time $T_0$ such that the probability to classify all the user in the network correctly for every $t>T_0$ is sufficiently high. We can see from Figures \ref{max_deviation_fig} and \ref{drift_fig} that the minimal value of  $T_0$ that yields a negligible deviation from true consensus value depends on probability of miss-classifying users, i.e., classifying a legitimate agent as malicious and vice-versa. This probability is increased, for example, as the number of malicious agents grows and as the variance of $\alpha_{ij}$ grows. Figures \ref{max_deviation_fig} and \ref{drift_fig} show that as we increase the number of malicious agents and the deviation $\ell$ that governs the variance of $\alpha_{ij}$ we should increase $T_0$ to limit the effect of the malicious agents on the consensus value. Additionally, Figures \ref{max_deviation_fig} and \ref{drift_fig} show 
 that the for $T_0=150$ the deviation from true consensus value for all considered scenarios is negligible. Finally, for $T_0=100$ the deviation from true consensus value is sufficiently small and the effect of malicious agents is limited even for $|\MM|=30$.

Figure~\ref{max_deviation_fig} shows that as  $\ell$ becomes larger (and therefore also the variance), the effect of malicious agents on deviation gets worse.  However, even in this case increasing $T_0$ can recover the nominal consensus value with deviation approaching zero.  The exact relationship between the deviation $\maxdev$ and $T_0$ is given by Theorem~\ref{thm:main2}.

In the case of a time varying malicious agent input as with the \emph{drift malicious input} model described in part (b) of this section, nominal average consensus can also be recovered by making $T_0$ larger, as depicted in Figure~\ref{drift_fig}.  This is true even with variable variance on the $\alpha_{ij}$ values.

\vspace{-0.1in}
\section{Discussion}
\label{sec:discussion}
In this section we discuss the significance of some of our characterizations on the trust observations $\alpha_{ij}$ and some of our assumptions. We also discuss directions for future work.  Firstly, the results of this paper assume a fixed network topology. However,  we believe that it is possible to relax Assumption~\ref{assumptions}.1 to include the case where the graph over legitimate agents is not always connected. In such a case Assumption~\ref{assumptions}.1 can be replaced with the assumption that the topology of the legitimate agents is B-connected.\footnote{See definition in \cite{8340193}.} We believe that the convergence results from \cite[Proposition 1]{Nedic_Ozdaglar_covergence_rate_delay} can be adapted to this case along the same vein as in Theorem \ref{thm:main3}.  Finally we discuss a few interesting observations about the trust values $\alpha_{ij}$. Our results show that the more informed our $\alpha_{ij}$ values are, the tighter the performance guarantees achievable become (cf. Lemma~\ref{Lemma:concentration_upper_variance}). This also opens a future avenue of investigation of finding $\alpha_{ij}$ observations from physical channels in multi-agent systems which can be characterized as fully as possible. The theory in this paper provides the mathematical framework for understanding \emph{which characteristics of} $\alpha_{ij}$ are most critical for attaining performance guarantees that are important for multi-agent coordination. We hope that this can guide continued investigation into trust-based resilience, particularly for cyberphysical systems where additional information on inter-agent trust can be attained.

\section{Conclusion}
\label{sec:conclusion}
This paper presents a unified mathematical theory to treat distributed multi-agent consensus in the presence of malicious agents when stochastic values of trust between agents are available. Under this model, we present three new performance guarantees for consensus systems in the presence of malicious agents:  1) convergence of consensus is possible with probability 1 even when the number of malicious agents is far larger than 1/2 of the network connectivity-- in contrast to classical results in resilience, 2) the deviation of the converged value from average consensus can be bounded and further, can be driven arbitrarily close to zero by using our derived edge weights and by allowing an observation time $T_0$ of trust values over the network, and 3) convergence to the correct classification of agents can be attained in finite time with probability 1 with an exponential rate that we derive. Additionally, we show that our performance guarantees can be made stronger if more information on the trust observations, such as bounds on their variance, is available. Taken together, these results point to the inherent value of quantifying and exploiting trust between agents for consensus. We believe that the mathematical formulations and framework of this paper hold promise for achieving a new generation of resilient multi-agent systems.
%\input{appendices}

%\begin{spacing}{0.85}
%\bibliographystyle{IEEEtran}
%\bibliographystyle{plainnat}
%\vspace{-0.1in}
%\bibliography{references}

\vspace{-0.1in}
\begin{thebibliography}{10}
\providecommand{\url}[1]{#1}
\csname url@rmstyle\endcsname
\providecommand{\newblock}{\relax}
\providecommand{\bibinfo}[2]{#2}
\providecommand\BIBentrySTDinterwordspacing{\spaceskip=0pt\relax}
\providecommand\BIBentryALTinterwordstretchfactor{4}
\providecommand\BIBentryALTinterwordspacing{\spaceskip=\fontdimen2\font plus
\BIBentryALTinterwordstretchfactor\fontdimen3\font minus
  \fontdimen4\font\relax}
\providecommand\BIBforeignlanguage[2]{{%
\expandafter\ifx\csname l@#1\endcsname\relax
\typeout{** WARNING: IEEEtran.bst: No hyphenation pattern has been}%
\typeout{** loaded for the language `#1'. Using the pattern for}%
\typeout{** the default language instead.}%
\else
\language=\csname l@#1\endcsname
\fi
#2}}

\bibitem{origByz}
\BIBentryALTinterwordspacing
D.~Dolev, ``The byzantine generals strike again,'' \emph{Journal of
  Algorithms}, vol.~3, no.~1, pp. 14 -- 30, 1982. [Online]. Available:
  \url{http://www.sciencedirect.com/science/article/pii/0196677482900049}
\BIBentrySTDinterwordspacing

\bibitem{nancy}
M.~J. Fischer, N.~A. Lynch, and M.~S. Paterson, ``Impossibility of distributed
  consensus with one faulty process,'' \emph{J. ACM}, vol.~32, no.~2, p.
  374–382, Apr. 1985.

\bibitem{paxos}
L.~Lamport, ``Paxos made simple,'' \emph{ACM SIGACT News 32}, pp. 18--25, Dec
  2001.

\bibitem{lynchPaxos}
R.~D. Prisco, B.~Lampson, and N.~Lynch, ``Revisiting the paxos algorithm,''
  \emph{Theoretical Computer Science}, vol. 243, no.~1, pp. 35 -- 91, 2000.

\bibitem{sundaram}
S.~{Sundaram} and C.~N. {Hadjicostis}, ``Distributed function calculation via
  linear iterative strategies in the presence of malicious agents,'' \emph{IEEE
  Transactions on Automatic Control}, vol.~56, no.~7, pp. 1495--1508, 2011.

\bibitem{bulloUnreliable}
F.~{Pasqualetti}, A.~{Bicchi}, and F.~{Bullo}, ``Consensus computation in
  unreliable networks: A system theoretic approach,'' \emph{IEEE Transactions
  on Automatic Control}, vol.~57, no.~1, pp. 90--104, Jan 2012.

\bibitem{security_swarmThreat}
F.~Higgins \emph{et~al.}, ``Threats to the swarm: Security considerations for
  swarm robotics,'' \emph{International Journal on Advances in Security}, 2009.

\bibitem{wmsr}
H.~J. {LeBlanc}, H.~{Zhang}, X.~{Koutsoukos}, and S.~{Sundaram}, ``Resilient
  asymptotic consensus in robust networks,'' \emph{IEEE Journal on Selected
  Areas in Communications}, vol.~31, no.~4, pp. 766--781, April 2013.

\bibitem{CPS_Design}
S.~K. {Khaitan} and J.~D. {McCalley}, ``Design techniques and applications of
  cyberphysical systems: A survey,'' \emph{IEEE Systems Journal}, vol.~9,
  no.~2, pp. 350--365, 2015.

\bibitem{BulloCyberphysSecurity_GeometricPrinciples}
F.~{Pasqualetti}, F.~{Dorfler}, and F.~{Bullo}, ``Control-theoretic methods for
  cyberphysical security: Geometric principles for optimal cross-layer
  resilient control systems,'' \emph{IEEE Control Systems Magazine}, vol.~35,
  no.~1, pp. 110--127, 2015.

\bibitem{pappasDesignActuationSensing}
S.~{Pequito}, F.~{Khorrami}, P.~{Krishnamurthy}, and G.~J. {Pappas}, ``Analysis
  and design of actuation–sensing–communication interconnection structures
  toward secured/resilient lti closed-loop systems,'' \emph{IEEE Transactions
  on Control of Network Systems}, vol.~6, no.~2, pp. 667--678, 2019.

\bibitem{Candell_PhysicsBasedDetection}
J.~Giraldo, D.~Urbina, A.~Cardenas, J.~Valente, M.~Faisal, J.~Ruths, N.~O.
  Tippenhauer, H.~Sandberg, and R.~Candell, ``A survey of physics-based attack
  detection in cyber-physical systems,'' \emph{ACM Comput. Surv.}, vol.~51,
  no.~4, Jul. 2018.

\bibitem{pappasSecretChannelCodes}
A.~{Tsiamis}, K.~{Gatsis}, and G.~J. {Pappas}, ``State-secrecy codes for
  networked linear systems,'' \emph{IEEE Transactions on Automatic Control},
  vol.~65, no.~5, pp. 2001--2015, 2020.

\bibitem{Sinopoli_Watermark}
Y.~{Mo}, S.~{Weerakkody}, and B.~{Sinopoli}, ``Physical authentication of
  control systems: Designing watermarked control inputs to detect counterfeit
  sensor outputs,'' \emph{IEEE Control Systems Magazine}, vol.~35, no.~1, pp.
  93--109, 2015.

\bibitem{fadel}
F.~Adib, Z.~Kabelac, D.~Katabi, and R.~C. Miller, ``3d tracking via body radio
  reflections,'' in \emph{11th {USENIX} Symposium on Networked Systems Design
  and Implementation ({NSDI} 14)}.\hskip 1em plus 0.5em minus 0.4em\relax
  Seattle, WA: {USENIX} Association, Apr. 2014, pp. 317--329.

\bibitem{arrayTrack}
J.~Xiong and K.~Jamieson, ``Arraytrack: A fine-grained indoor location
  system,'' in \emph{Presented as part of the 10th {USENIX} Symposium on
  Networked Systems Design and Implementation ({NSDI} 13)}.\hskip 1em plus
  0.5em minus 0.4em\relax Lombard, IL: {USENIX}, 2013, pp. 71--84.

\bibitem{spotFi}
\BIBentryALTinterwordspacing
M.~Kotaru, K.~Joshi, D.~Bharadia, and S.~Katti, ``Spotfi: Decimeter level
  localization using wifi,'' \emph{SIGCOMM Comput. Commun. Rev.}, vol.~45,
  no.~4, p. 269–282, Aug. 2015. [Online]. Available:
  \url{https://doi.org/10.1145/2829988.2787487}
\BIBentrySTDinterwordspacing

\bibitem{RF-Compass}
J.~Wang, F.~Adib, R.~Knepper, D.~Katabi, and D.~Rus, ``Rf-compass: Robot object
  manipulation using rfids,'' ser. MobiCom, 2013.

\bibitem{SecureArray}
J.~Xiong and K.~Jamieson, ``Securearray: Improving wifi security with
  fine-grained physical-layer information,'' in \emph{Proceedings of the 19th
  Annual International Conference on Mobile Computing and Networking}, ser.
  MobiCom, 2013.

\bibitem{AURO}
S.~Gil, S.~Kumar, M.~Mazumder, D.~Katabi, and D.~Rus, ``Guaranteeing
  spoof-resilient multi-robot networks,'' \emph{AuRo}, p. 1383–1400, 2017.

\bibitem{murray}
J.~Fax and R.~Murray, ``Information flow and cooperative control of vehicle
  formations,'' \emph{Automatic Control, IEEE Transactions on}, 2004.

\bibitem{lynchText}
N.~A. Lynch, \emph{Distributed Algorithms}.\hskip 1em plus 0.5em minus
  0.4em\relax San Francisco, CA, USA: Morgan Kaufmann Publishers Inc., 1996.

\bibitem{pappas}
M.~M. {Zavlanos}, M.~B. {Egerstedt}, and G.~J. {Pappas}, ``Graph-theoretic
  connectivity control of mobile robot networks,'' \emph{Proceedings of the
  IEEE}, vol.~99, no.~9, pp. 1525--1540, 2011.

\bibitem{spanos}
D.~P. Spanos and R.~M. Murray, ``Motion planning with wireless network
  constraints,'' in \emph{ACC}, 2005.

\bibitem{nedicConstrainedConsensus}
A.~Nedi\'c, A.~{Ozdaglar}, and P.~A. {Parrilo}, ``Constrained consensus and
  optimization in multi-agent networks,'' \emph{IEEE Transactions on Automatic
  Control}, vol.~55, no.~4, pp. 922--938, 2010.

\bibitem{nedicGeometricConvergence}
A.~Nedi\'c, A.~Olshevsky, and W.~Shi, ``Achieving geometric convergence for
  distributed optimization over time-varying graphs,'' \emph{SIAM Journal on
  Optimization}, vol.~27, no.~4, pp. 2597--2633, 2017.

\bibitem{TsitsiklisSwitchingConvergence}
J.~M. {Hendrickx} and J.~N. {Tsitsiklis}, ``Convergence of type-symmetric and
  cut-balanced consensus seeking systems,'' \emph{IEEE Transactions on
  Automatic Control}, vol.~58, no.~1, pp. 214--218, 2013.

\bibitem{olfati-saber}
R.~Olfati-Saber and R.~M. Murray, ``Consensus problems in networks of agents
  with switching topology and time-delays,'' \emph{IEEE Transactions on
  Automatic Control}, 2004.

\bibitem{8340193}
A.~{Nedi\'c}, A.~{Olshevsky}, and M.~G. {Rabbat}, ``Network topology and
  communication-computation tradeoffs in decentralized optimization,''
  \emph{Proceedings of the IEEE}, vol. 106, no.~5, pp. 953--976, 2018.

\bibitem{boyd}
L.~Xiao and S.~Boyd, ``Fast linear iterations for distributed averaging,''
  \emph{Systems and Control Letters}, 2004.

\bibitem{nedicConvRate}
A.~Nedi\'c, ``Convergence rate of distributed averaging dynamics and
  optimization in networks,'' \emph{Foundations and Trends® in Systems and
  Control}, vol.~2, no.~1, pp. 1--100, 2015.

\bibitem{TsitiklisConvergenceRate}
A.~Olshevsky and J.~N. Tsitsiklis, ``Convergence speed in distributed consensus
  and averaging,'' \emph{SIAM Rev.}, vol.~53, no.~4, p. 747–772, Nov. 2011.

\bibitem{additional_reviewer_reference1}
X.~{Liu} and J.~S. {Baras}, ``Using trust in distributed consensus with
  adversaries in sensor and other networks,'' in \emph{17th International
  Conference on Information Fusion (FUSION)}, 2014, pp. 1--7.

\bibitem{additional_reviewer_reference3}
W.~{Tong}, X.~{Dong}, and J.~{Zheng}, ``Trust-pbft: A peertrust-based practical
  byzantine consensus algorithm,'' in \emph{2019 International Conference on
  Networking and Network Applications (NaNA)}, 2019, pp. 344--349.

\bibitem{maliciousModel}
I.~Sargeant and A.~Tomlinson, ``Modelling malicious entities in a robotic
  swarm,'' in \emph{DASC}, 2013.

\bibitem{swarmThreat}
F.~Higgins, A.~Tomlinson, and K.~M. Martin, ``Threats to the swarm: Security
  considerations for swarm robotics,'' \emph{International Journal on Advances
  in Security}, vol.~2, pp. 288--297, 2009.

\bibitem{surveySecurity}
F.~{Higgins}, A.~{Tomlinson}, and K.~M. {Martin}, ``Survey on security
  challenges for swarm robotics,'' in \emph{Fifth International Conference on
  Autonomic and Autonomous Systems (ICAS) 2009}, 2009.

\bibitem{sastryBook}
A.~A. C\'ardenas, T.~Roosta, G.~Taban, and S.~Sastry, \emph{Cyber Security:
  Basic Defenses and Attack Trends}, 2008, pp. 73--101.

\bibitem{scada}
A.~A. C\'{a}rdenas \emph{et~al.}, ``Attacks against process control systems:
  Risk assessment, detection, and response,'' ser. ASIACCS '11, 2011.

\bibitem{Yang13}
J.~{Yang}, Y.~{Chen}, W.~{Trappe}, and J.~{Cheng}, ``Detection and localization
  of multiple spoofing attackers in wireless networks,'' \emph{IEEE
  Transactions on Parallel and Distributed Systems}, vol.~24, no.~1, pp.
  44--58, 2013.

\bibitem{ControlSysMagazineResilienceIssue}
H.~{Sandberg}, S.~{Amin}, and K.~H. {Johansson}, ``Cyberphysical security in
  networked control systems: An introduction to the issue,'' \emph{IEEE Control
  Systems Magazine}, vol.~35, no.~1, pp. 20--23, 2015.

\bibitem{sastry_saurabh_DOS}
S.~Amin, A.~A. Cárdenas, and S.~S. Sastry, ``Safe and secure networked control
  systems under denial-of-service attacks,'' \emph{HSCC}, pp. 31--45, 2009.

\bibitem{sastryCPS_survivability}
A.~A. {Cardenas}, S.~{Amin}, and S.~{Sastry}, ``Secure control: Towards
  survivable cyber-physical systems,'' in \emph{2008 The 28th International
  Conference on Distributed Computing Systems Workshops}, 2008, pp. 495--500.

\bibitem{prorok}
K.~Saulnier, D.~Saldaña, A.~Prorok, G.~J. Pappas, and V.~Kumar, ``Resilient
  flocking for mobile robot teams,'' \emph{IEEE Robotics and Automation
  Letters}, vol.~2, no.~2, pp. 1039--1046, April 2017.

\bibitem{dataVeracityCPS_entropyMethods}
M.~Krotofil, J.~Larsen, and D.~Gollmann, ``The process matters: Ensuring data
  veracity in cyber-physical systems,'' in \emph{Proceedings of the 10th ACM
  Symposium on Information, Computer and Communications Security}.\hskip 1em
  plus 0.5em minus 0.4em\relax New York, NY, USA: Association for Computing
  Machinery, 2015, p. 133–144.

\bibitem{modelingDependability_CPS}
J.~{Lin}, S.~{Sedigh}, and A.~{Miller}, ``A general framework for quantitative
  modeling of dependability in cyber-physical systems: A proposal for doctoral
  research,'' in \emph{2009 33rd Annual IEEE International Computer Software
  and Applications Conference}, vol.~1, 2009, pp. 668--671.

\bibitem{goldsmith}
A.~Goldsmith, \emph{Wireless Communications}.\hskip 1em plus 0.5em minus
  0.4em\relax New York, NY, USA: Cambridge University Press, 2005.

\bibitem{davidtse}
D.~Tse and P.~Vishwanath, \emph{Fundamentals of Wireless Communications}.\hskip
  1em plus 0.5em minus 0.4em\relax Cambridge University Press, 2005.

\bibitem{lteye}
S.~Kumar, E.~Hamed, D.~Katabi, and L.~Erran~Li, ``{LTE} radio analytics made
  easy and accessible,'' in \emph{SIGCOMM}, 2014.

\bibitem{ubicarse}
S.~Kumar \emph{et~al.}, ``Accurate indoor localization with zero start-up
  cost,'' ser. MobiCom, 2014.

\bibitem{Wifi_collaborative}
N.~Jadhav, W.~Wang, D.~Zhang, O.~Khatib, S.~Kumar, and S.~Gil, ``{WSR: A WiFi}
  sensor for collaborative robotics,'' in \emph{arXiv:2012.04174}, 2020.

\bibitem{kai}
K.~{Zeng}, K.~{Govindan}, and P.~{Mohapatra}, ``Non-cryptographic
  authentication and identification in wireless networks [security and privacy
  in emerging wireless networks],'' \emph{IEEE Wireless Communications},
  vol.~17, no.~5, pp. 56--62, 2010.

\bibitem{ting}
T.~Wang and Y.~Yang, ``Analysis on perfect location spoofing attacks using
  beamforming,'' in \emph{INFOCOM}, 2013.

\bibitem{mobicom_adsb}
D.~Moser, P.~Leu, V.~Lenders, A.~Ranganathan, F.~Ricciato, and S.~Capkun,
  ``Investigation of multi-device location spoofing attacks on air traffic
  control and possible countermeasures,'' in \emph{Proceedings of the 22Nd
  Annual International Conference on Mobile Computing and Networking}, ser.
  MobiCom '16.\hskip 1em plus 0.5em minus 0.4em\relax New York, NY, USA: ACM,
  2016, pp. 375--386.

\bibitem{IJRR}
S.~Gil, S.~Kumar, M.~Mazumder, D.~Katabi, and D.~Rus, ``Adaptive communication
  in multi-robot systems using directionality of signal strength,''
  \emph{IJRR}, 2015.

\bibitem{GilLCSS}
S.~{Gil}, C.~{Baykal}, and D.~{Rus}, ``Resilient multi-agent consensus using
  wi-fi signals,'' \emph{IEEE Control Systems Letters}, vol.~3, no.~1, pp.
  126--131, Jan 2019.

\bibitem{ICRA2019}
T.~{Wheeler}, E.~{Bharathi}, and S.~{Gil}, ``Switching topology for resilient
  consensus using wi-fi signals,'' in \emph{2019 International Conference on
  Robotics and Automation (ICRA)}, 2019, pp. 2018--2024.

\bibitem{Dubhashi2009ConcentrationOM}
D.~P. Dubhashi and A.~Panconesi, ``Concentration of measure for the analysis of
  randomized algorithms.''\hskip 1em plus 0.5em minus 0.4em\relax Cambridge
  University Press, 2009.

\bibitem{Perron_Frobenius_ref}
S.~U. {Pillai}, T.~{Suel}, and {Seunghun Cha}, ``The perron-frobenius theorem:
  some of its applications,'' \emph{IEEE Signal Processing Magazine}, vol.~22,
  no.~2, pp. 62--75, 2005.

\bibitem{durrett_2010}
R.~Durrett, \emph{Probability: Theory and Examples}, 4th~ed.\hskip 1em plus
  0.5em minus 0.4em\relax Cambridge University Press, 2010.

\bibitem{markov_convergence_rate}
P.~Br\'emaud, \emph{Markov Chains: Gibbs Fields, Monte Carlo Simulation, and Queues}.\hskip
  1em plus 0.5em minus 0.4em\relax New York, USA: Springer-Verlag, 1999.

\bibitem{doi:10.1080/01621459.1962.10482149}
G.~Bennett, ``Probability inequalities for the sum of independent random
  variables,'' \emph{Journal of the American Statistical Association}, vol.~57,
  no. 297, pp. 33--45, 1962.

\bibitem{doi:10.1080/03610926.2017.1367818}
S.~Zheng, ``An improved {B}ennett's inequality,'' \emph{Communications in
  Statistics - Theory and Methods}, vol.~47, no.~17, pp. 4152--4159, 2018.

\bibitem{Jabera2018}
T.~Jebara, ``A refinement of {B}ennett’s inequality with applications to
  portfolio optimization,'' in \emph{arXiv:1804.05454}, 2018.

\bibitem{Nedic_Ozdaglar_covergence_rate_delay}
A.~{Nedi\'c} and A.~{Ozdaglar}, ``Convergence rate for consensus with delays,''
  \emph{Journal on Global Optimization}, vol.~47, p. 437–456, 2010.

\end{thebibliography}
%\end{spacing}

% or


\begin{spacing}{0.85}
\bibliographystyle{IEEEtran}

\end{spacing}

%\vspace{-0.2in}
%\section*{Acknowledgment}

%The authors gratefully acknowledge partial funding support from NSF CAREER grant \#1845225.
%\vspace{-0.2in}

\ifCLASSOPTIONcaptionsoff
  \newpage
\fi

\end{document}